\documentclass{amsart}

\usepackage{rlepsf,amssymb}

\newtheorem{thm}{Theorem}[section]
\newtheorem{lem}[thm]{Lemma}
\newtheorem{cor}[thm]{Corollary}
\newtheorem{prop}[thm]{Proposition}
\newtheorem{rem}[thm]{Remark}

\newtheorem*{quest}{Question}
\newtheorem{conj}{Conjecture}
\newtheorem*{warning}{Warning}

\newtheoremstyle{definition}{7pt plus6.3pt minus6.3pt}{7pt plus3pt minus3pt}%
{\rm}{}{\bf}{}{0.75em}{\thmname{#1}\thmnumber{ #2}\thmnote{\sl\stdspace#3}}
\theoremstyle{definition}
\newtheorem{example}[thm]{Example}
\newtheorem{exercise}[thm]{\small Exercise}

\newcommand{\bbr}{\begin{rem}\em} 
\newcommand{\eer}{\end{rem}}

\newcommand{\bex}{\begin{example}} 
\newcommand{\eex}{\end{example}}

\newcommand{\bhw}{\begin{exercise}\small} 
\newcommand{\ehw}{\end{exercise}}

\newcommand{\beq}{\begin{equation}} 
\newcommand{\eeq}{\end{equation}}

\newcommand{\be}{\begin{enumerate}}
\newcommand{\ee}{\end{enumerate}}

\def\tb{\operatorname{tb}}
\def\C{\hbox{$\mathbb C$}}
\def\Z{\hbox{$\mathbb Z$}}
\def\Q{{\hbox{$\mathbb Q$}}}
\def\R{\hbox{$\mathbb R$}}
\def\K{\hbox{$\mathcal{K}$}}
\def\L{\hbox{$\mathcal{L}$}}
\def\T{\hbox{$\mathcal{T}$}}
\def\A{\hbox{$\mathcal{A}$}}
\def\CC{\hbox{$\mathcal{C}$}}
\def\M{\hbox{$\mathcal{M}$}}

\def\dfn#1{{\em #1}}

\begin{document}

\title{Legendrian and Transversal Knots}

\author{John B. Etnyre}
\address{University of Pennsylvania, Philadelphia, PA 19104}
\email{etnyre@math.upenn.edu}
\urladdr{http://www.math.upenn.edu/\char126 etnyre}

\thanks{Supported in part by NSF Grant \# DMS-0203941.}


\maketitle


\section{Introduction}
Contact structures on manifolds and Legendrian and transversal knots
in them are very natural objects, born over two centuries ago, in the
work of Huygens, Hamilton and Jacobi on geometric optics and work of
Lie on partial differential equations.  They touch on diverse areas of
mathematics and physics, and have deep connections with topology and
dynamics in low dimensions. The study of Legendrian knots is now a rich
and beautiful theory with many applications. This survey is an
introduction to, and overview of the current state of knot theory in
contact geometry. For a discussion of the driving questions in the
field see \cite{EtnyreNg} and for a more historical discussion of
contact geometry see \cite{Geiges01}.

This paper will concentrate on Legendrian and transversal knots in
dimension three where their theory is most fully developed and where
they are most intimately tied to topology. Moreover, in this dimension
one may use a predominately topological and combinatorial approach to
their study.  In Section~\ref{sec:hd} we will make our only excursion
into the study of higher dimensional Legendrian knots.

Throughout this survey we assume the reader is familiar with basic
topology at the level of \cite{Rolfsen}. We have tried to keep the
contact geometry prerequisites to a minimum, but it would certainly be
helpful to have had some prior exposure to the basics as can be found in
\cite{EtnyreIntro, GeigesIntro}. Some of the proofs in
Section~\ref{sec:class} rely on convex surface theory which can be
found in \cite{EtnyreIntro}, but references to convex surfaces can be
largely ignored without serious loss of continuity.

\tableofcontents
\vfill\eject

This paper is part of the forthcoming book:
\bigskip

\centerline{\large{{\bf Table of Contents for the Handbook of Knot
Theory}}}

\

\centerline{William W. Menasco and Morwen B. Thistlethwaite, Editors}%

\

\begin{enumerate}%

\item Colin Adams, {\it Hyperbolic knots}
\item Joan S. Birman and Tara Brendle {\it Braids: A Survey}
\item John Etnyre, {\it Legendrian and transversal knots}
\item Greg Friedman {\it Knot Spinning}
\item Jim Hoste {\it The enumeration and classification of knots and
links}
\item Louis Kauffman  {\it Knot Diagrammatics}
\item Charles Livingston  {\it A survey of classical knot concordance}
\item Lee Rudolph  {\it Knot theory of complex plane curves}
\item Marty Scharlemann  {\it Thin position in the theory of classical knots}
\item Jeff Weeks  {\it Computation of hyperbolic structures in knot
theory}

\end{enumerate}
\vfill\eject

\section{Definitions and examples}
A \dfn{contact structure} on an oriented 3--manifold $M$ is a completely non-integrable plane field
$\xi$ in the tangent bundle of $M.$ That $\xi$ is a plane field simply means that at each point
$x\in M,$ $\xi_x$ is a two dimensional subspace of $T_xM.$ Near each point one may always describe a
plane field as the kernel of a (locally defined) 1--form, {\em i.e.} there is a 1--form $\alpha$ so
that $\xi_x=ker (\alpha_x).$ The plane filed $\xi$ is completely non-integrable if any 1--form
$\alpha$ defining $\xi$ satisfies $\alpha\wedge d\alpha\not=0.$ We will require that $\alpha\wedge
d\alpha$ defines the given orientation on $M.$ This orientation compatibility has become a fairly
standard part of the definition of contact structure, but in the past such a contact structure was
called a positive contact structure. The condition $\alpha\wedge d\alpha\not=0$ implies that $\xi$
is not everywhere tangent to any surface. Intuitively one can see this in Figure~\ref{fig:fex}.
There one sees the plane fields twisting. This twisting prevents the planes from being everywhere
tangent to a surface. To make this ``not everywhere tangent'' condition precise one should consult
the Frobenius Theorem \cite{Conlon}.

\subsection{The standard contact structure on $\R^3$}
\begin{figure}[ht]
  \relabelbox 
  \small {\epsfxsize=4.5in\centerline{\epsfbox{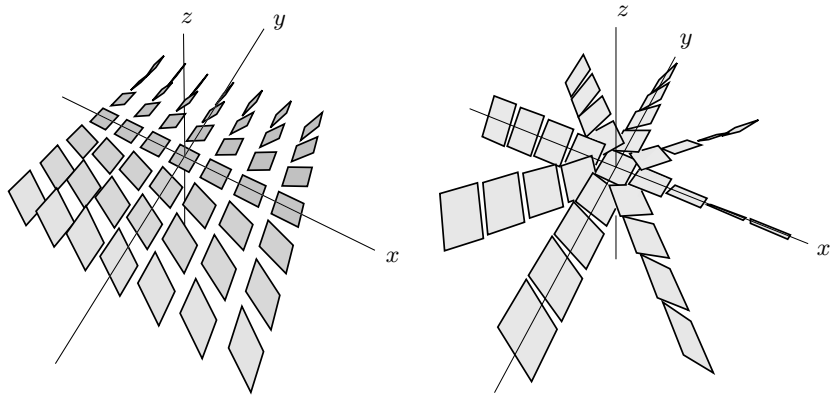}}} 
  \relabel {1}{$x$}
  \relabel {2}{$y$}
  \relabel {3}{$z$}
  \relabel {4}{$x$}
  \relabel {5}{$y$}
  \relabel {6}{$z$}
  \endrelabelbox
  \caption{The contact structure $\xi_{std}$ (left) and
    $\xi_{sym}$ (right) on $\R^3.$ (Figures courtesy of S.\ Sch\"onenberger.)}
  \label{fig:fex}
\end{figure}

The simplest example of a contact structures is on $\R^3$ (with Cartesian coordinates $(x,y,z)$) and
is given by
\[\xi_{std}=\text{span}\left\{\frac{\partial}{\partial y}, 
        \frac{\partial}{\partial x}+y\frac{\partial} {\partial z}\right\}.\] 
Clearly $\xi_{std}$ is the kernel of the 1--from $\alpha=dz-ydx$ and thus one may easily verify that
it is a contact structure. The plane field is indicated in Figure~\ref{fig:fex}.
We observe a few things about the plane fields that will be useful later. Note that the plane fields
along the $xz$-plane are all horizontal ({\em ie.} parallel to the $xy$-plane). Moving out along the
$y$ axis the planes start horizontal and ``twist'' around the $y$ axis in a left handed manner.
Traversing the $y$ axis from the origin ``to infinity'' the planes will make a $90^\circ$ twist.
Somehow it is this twisting that makes the plane field contact.

There are many other contact structures on $\R^3$ and all 3-manifolds have many contact structures,
but Darboux's theorem \cite{a, McDuffSalamon} says that all of them locally look like this one. By
this we mean that any point in a manifold has a neighborhood that is diffeomorphic to a neighborhood
of the origin in $\R^3$ by a diffeomorphism that takes the contact structure on the manifold to the
one described above. Throughout much of this article we will restrict attention to the contact
structure above. This is largely to simply the discussion. In a few places it will be essential that
we are in $(\R^3, \xi_{std})$ and we will make that clear then.

\subsection{Other contact structures}
We begin with a symmetric version of $\xi_{std}$ on $\R^3.$ 

\bex
Let $\alpha=dz+xdy-ydx$ or in cylindrical coordinates $\alpha = dz+r^2d\theta.$ The contact
structure $\xi_{sym}=\ker \alpha$ is as shown on the right hand side of Figure~\ref{fig:fex}. One
can find a diffeomorphism of $\R^3$ taking $\xi_{std}$ to $\xi_{sym}.$ Such a diffeomorphism is
called a \dfn{contactomorphism}. Thus in some sense they are the same contact structure, but for
various purposes one is sometimes easier to work with than the other.
\eex
 
\begin{figure}[ht]
        {\epsfxsize=2.5in\centerline{\epsfbox{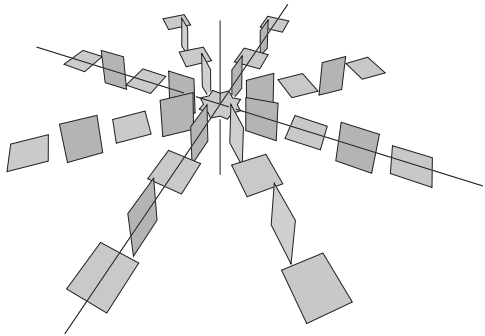}}}
        \caption{The contact structure $\xi_{ot}$ on $\R^3.$ 
          (Figures courtesy of S.\ Sch\"onenberger.)}
        \label{fig:otex}
\end{figure}

\bex\label{otex} 
In $\R^3$ with cylindrical coordinates consider $\alpha=\cos r dz+ r\sin r d\theta.$ The contact
structure $\xi_{ot} \ker \alpha$ is as shown in Figure~\ref{fig:otex}. Note the plane fields are
similar to the ones in $\xi_{sym}$ but as you move out along rays perpendicular to the $z$-axis the
planes twist around many times (infinitely often) whereas for $\xi_{sym}$ they twist only
$90^\circ.$
\eex 
Notice in this last example the disk $D=\{(r,\theta, z) | z=0, r\leq \pi\}$ is tangent to the
contact planes along the boundary. Such a disk is called an \dfn{overtwisted disk}. A contact
structure is called \dfn{overtwisted} if there is an embedded overtwisted disk, otherwise it is
called \dfn{tight}. From the discussion above we see that the plane field must twist to be a contact
structure, but $\xi_{ot}$ twists much more that is necessary, hence the name ``overtwisted''.
Clearly $\xi_{ot}$ is overtwisted and below we will see that $\xi_{std}$ and $\xi_{sym}$ are tight.
Moreover, we will also see that tight contact structures tell us interesting things about knots
related to them. For more discussion of tight and overtwisted contact structures see
\cite{Eliashberg89, EtnyreIntro}. Throughout this survey we will be mainly concerned with tight
contact structures with a real focus on $(\R^3, \xi_{std}).$

We end this section with a contact structure on a closed manifold.
\bex 
Let $S^3$ be the unit three sphere in $\R^4.$ Let
\[
\alpha=i^*\frac12(x_1dy_2-y_1dx_1+x_2dy_2-y_2dx_2)
\]
where $i:S^3\to R^4$ is the inclusion map. One may readily check that $\alpha$ is a contact form so
$\xi=\ker\alpha$ is a contact structure on $S^3.$ If one removes a point form $S^3$ the contact
structure $\xi$ is contactomorphic to $\xi_{std}$ on $\R^3.$ For a more thorough discussion of this
and similar examples see \cite{EtnyreIntro}.
\eex
For other examples of contact structures on closed manifolds see
\cite{a, EtnyreIntro, Gompf98, McDuffSalamon}.

\subsection{Legendrian knots}\label{sec:lk}
A \dfn{Legendrian knot} $L$ in a contact manifold $(M^3,\xi)$ is an embedded $S^1$ that is always
tangent to $\xi:$
\[
T_xL\in \xi_x,\quad x\in L.
\] 

For the rest of this section we restrict attention to $(\R^3,\xi_{std}).$ This is for two reasons.
First, according to Darboux's theorem all contact structures look locally like
$(\R^3,\xi_{\text{std}})$ and thus we are studying the ``local Legendrian knot theory'' in any
contact manifold. Secondly, when studying this contact structure we can use various projections to
help us understand the Legendrian knots, allowing us to ``get a feel for Legendrian knots''.

Below it will be convenient to have a parameterization of $L.$ So throughout this section
\[
\phi:S^1\to \R^3: \theta\mapsto (x(\theta), y(\theta), z(\theta))
\]
will be a parameterization of $L.$ Moreover, we will assume that $\phi$ is at least a $C^1$
immersion. Now the fact that $L$ is tangent to $\xi$ can be easily expressed by
\[
\phi'(\theta)\in \xi_{\phi(\theta)}
\]
or since $\xi=\ker (dz-ydx)$
\begin{equation}\label{mainleg}
z'(\theta)-y(\theta)x'(\theta)=0.
\end{equation} 

There are two ways to picture Legendrian knots in $(\R^3,\xi_{std}),$ that is via the {\em front
projection} and the {\em Lagrangian projection}. We begin with the front projection. Let
\[
\Pi:\R^3\to \R^2: (x,y,z)\to (x,z).
\] 
The image, $\Pi(L),$ of $L$ under the map $\Pi$ is called the front projection of $L.$ If $\phi$
above parameterizes $L$ then
\[
\phi_\Pi:S^1\to \R^2:\theta\mapsto (x(\theta), z(\theta))
\]
parameterizes $\Pi(L).$ While $\phi$ was an immersion (in fact an embedding), $\phi_\Pi$ will
certainly not be one. To see this note that Equation~\eqref{mainleg} implies that
$z'(\theta)=y(\theta)x'(\theta)$ thus anytime $x'(\theta)$ vanishes so must $z'(\theta).$ So if
$\psi_\Pi$ is to be an immersion $x'(\theta)$ must never vanish. But this implies that $\Pi(L)$ has
no vertical tangencies and of course any immersion of $S^1$ into $\R^2$ must have vertical
tangencies. (Here and below ``vertical'' means spanned by $\frac{\partial}{\partial z}$.) This
brings us to our first important fact about front projections.
\begin{itemize}
\item[{\em FF 1.}] Front projections $\Pi(L)$ have no vertical tangencies.
\end{itemize}
So how does $\psi_\Pi$ fail to be an immersion? Note the discussion above implies that $z'(\theta)$
always vanishes to at least the order of $x'(\theta).$ Thus we can always recover the $y$-coordinate
of $\phi$ from $\phi_\Pi$ by rewriting Equation~\eqref{mainleg} as 
\beq\label{eqforx}
y(\theta)=\frac{z'(\theta)}{x'(\theta)}, 
\eeq 
if $x'(\theta)$ is non-zero. If $x'$ is zero at some $\theta_0$ then we have 
\beq 
y(\theta_0)=\lim_{\theta\to \theta_0} \frac{z'(\theta)}{x'(\theta)}.
\eeq

\begin{figure}[ht]
  \relabelbox \small {\epsfxsize=2.5in\centerline{\epsfbox{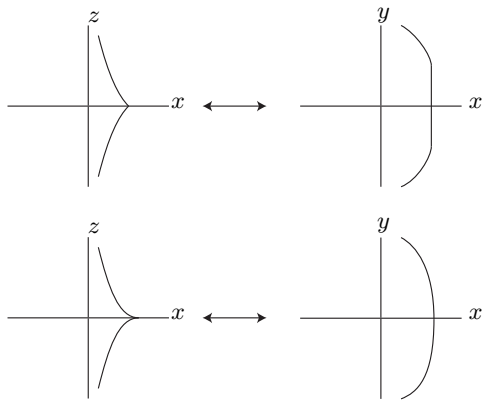}}} 
  \relabel {x1}{$x$}
  \relabel {x2}{$x$} 
  \relabel {x3}{$x$}   
  \relabel {x4}{$x$} 
  \relabel {y}{$y$} 
  \relabel {y6}{$y$}
  \relabel {z}{$z$} 
  \relabel {z8}{$z$} 
  \endrelabelbox
  \caption{Top row shows two projections of a Legendrian arc
    with $x'(\theta)=0$ on a open interval. Bottom row shows a near by Legendrian arc with
    $x'(\theta)=0$ only at one point.}
  \label{fig:good}
\end{figure}

There are Legendrian knots for which $x'(\theta)=0$ on open intervals, see top of
Figure~\ref{fig:good}, but this is not a stable phenomena. We can easily make the cusp point in the
projection ``sharper'', see bottom of Figure~\ref{fig:good}.
This leads to a Legendrian isotopic arc where $x'(\theta)=0$ only at one point. From this it is easy
to convince oneself that for a generic $C^1$ smooth Legendrian embedding in $\R^3,$ $x'(\theta)$ can
only vanish at isolated points. Moreover, at these isolated points there is a well defined tangent
line in the front projection. Thus we may assume that
\begin{itemize}
\item[{\em FF 2.}] Front projections may be parameterized by a map that is an immersion except at a
  finite number of points, at which there is still a well defined tangent line. Such points are
  called \dfn{generalized cusps}.
\end{itemize}
Actually this condition only guarantees the $y$ coordinate defined by Equation~\eqref{eqforx} is a
$C^0$ function. We need to add conditions on the second derivatives of $x$ and $z$ to get the $y$
coordinate to be $C^1.$ We will not concern ourselves here with this. It is interesting to note that
if we demand that all the coordinates be $C^\infty$ then generically our cusps must be ``semi-cubic
parabolas''. By this we mean that after a change of coordinates $z(\theta)=3\theta^3$ and
$x(\theta)= 2\theta^2.$

Our above discussion implies that {\em FF 1.\ }and {\em FF 2.\ }characterize front projections. In
particular, any map $f:S^1\to xz\text{--plane}:\theta\mapsto (x_f(\theta), z_f(\theta))$ that
satisfies {\em FF 1.\ }and {\em FF 2.\ }represents a Legendrian knot, since under these conditions
we can always define $y(\theta)$ by Equation~\eqref{eqforx}. Then the image of the map
$\phi(\theta)=(x_f(\theta), y(\theta), z_f(\theta))$ will be a Legendrian knot.

\begin{figure}[ht]
        {\epsfxsize=3in\centerline{\epsfbox{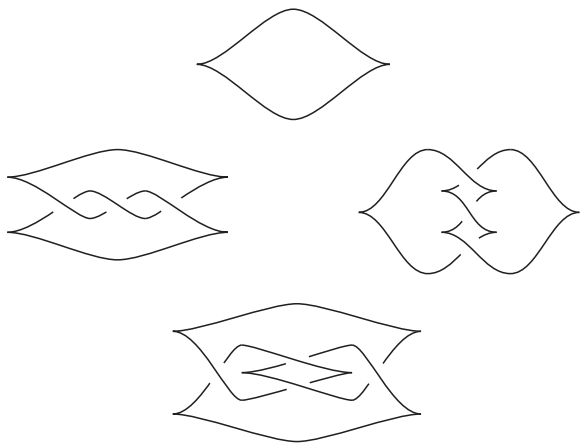}}}
        \caption{Legendrian knots realizing the unknot, right and left
        trefoil knots and the figure eight knot.}
        \label{fig:three}
\end{figure}

Interpreting this in terms of knot diagrams one sees that given any knot diagram
\begin{enumerate}
\item that has no vertical tangencies,
\item the only non-smooth points are generalized cusps and
\item at each crossing the slope of the overcrossing is smaller (that is, more negative) than the 
undercrossing
\end{enumerate}
represents the front projection of a Legendrian knot. To understand the condition on the
overcrossing recall that in order to have the standard orientation on $\R^3$ the positive $y$ axis
goes into the page. 
\bex 
In Figure~\ref{fig:three} we show front diagrams for Legendrian knots realizing the unknot, right and 
left trefoil knots and the figure eight knot.
\eex

Using these observations one can easily prove.
\begin{thm}\label{thm:approx}
  Given any topological knot $K$ there is a Legendrian knot $C^0$ close to it. In particular, there
  are Legendrian knots representing any topological knot type.
\end{thm}

\begin{proof}
Consider $(\R^3,\xi_{std}).$ First we show that any knot type can be represented by a Legendrian
knot. This is now quite simple, just take any diagram for the knot (Figure~\ref{fig:k2l}), make the
modifications shown if Figure~\ref{fig:conv}
\begin{figure}[ht]
  {\epsfxsize=3in\centerline{\epsfbox{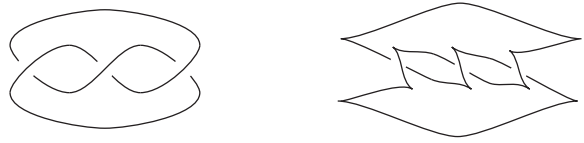}}}
  \caption{Converting a knot diagram (left) into a Legendrian
    front (right).}
  \label{fig:k2l}
\end{figure}
\begin{figure}[ht]
  \relabelbox 
  \small {\epsfxsize=2.8in\centerline{\epsfbox{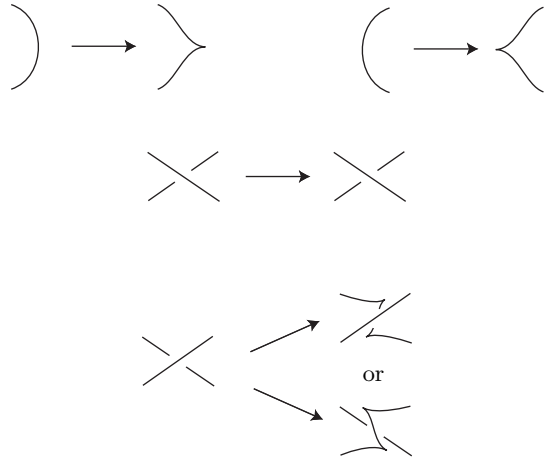}}}
  \relabel {or}{or}
  \endrelabelbox
  \caption{Realizing a knot type as a Legendrian knot.}
  \label{fig:conv}
\end{figure}
and then use the above procedure to recover a Legendrian knot (Figure~\ref{fig:k2l}).

Of course the problem with the construction of a Legendrian knot this way is that it is certainly
not $C^0$ close to the original knot since the difference between $y$-coordinates of the original
knot and the Legendrian knot determined by the front projection can be quite far apart. This problem
with the $y$-coordinate can be fixed with the idea illustrated in Figure~\ref{fig:arcapp} that any
arc may be $C^0$ approximated rel end points by a Legendrian arc.
\begin{figure}[ht]
        {\epsfxsize=1.8in\centerline{\epsfbox{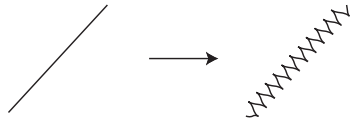}}}
        \caption{$C^0$ approximating an arc by a Legendrian arc.}
        \label{fig:arcapp}
\end{figure}

Using Darboux's theorem (that all contact structures are locally the same as
$(\R^3,\xi_{\text{std}})$) one may easily finish the proof for a general contact 3--manifold.
\end{proof}

Of course we could have used this last technique to show that any knot type has a Legendrian
representative thus avoiding the first part of the proof. But, in practice, if one is trying to
construct Legendrian representatives of a knot type (and not $C^0$ approximations of a specific
knot) one uses Figure~\ref{fig:conv}. This is because using Figure~\ref{fig:arcapp} introduces too
many ``zig-zags''. We will see below it is best to avoid these as much as possible.

\begin{figure}[ht]
  {\epsfxsize=2.4in\centerline{\epsfbox{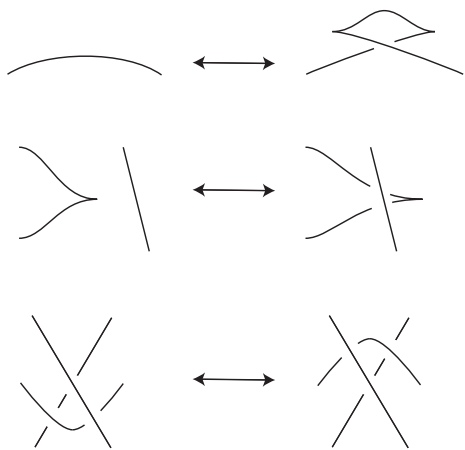}}}
  \caption{Legendrian Reidemeister moves. (Also need the
    corresponding figures rotated 180 degrees about all three
    coordinate axes.)}
  \label{fig:lreid}
\end{figure}

Just as there are Reidemeister moves for topological knot diagrams there is a set of
``Reidemeister'' moves for front diagrams too.
\begin{thm}[See \cite{Swiatkowski}]\label{thm:reidfront}
  Two front diagrams represent the Legendrian isotopic Legendrian knots if and only if they are
  related by regular homotopy and a sequence of moves shown in Figure~\ref{fig:lreid}.
\end{thm}

\begin{figure}[ht]
  {\epsfxsize=4in\centerline{\epsfbox{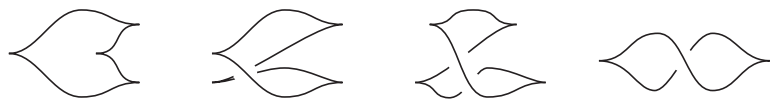}}}
  \caption{Various fronts of the same Legendrian unknot.}
  \label{fig:uniso}
\end{figure}

\begin{figure}[ht]
        {\epsfxsize=3in\centerline{\epsfbox{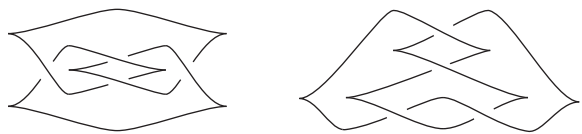}}}
        \caption{Two fronts of the same Legendrian figure eight knot.}
        \label{fig:two8}
\end{figure}

\bex 
If Figure~\ref{fig:uniso} we show that two different front diagrams for a Legendrian unknot are
Legendrian isotopic.
\eex 
To illustrate the difficulty in using these move one might want to try to show the two Legendrian
figure eight knots in Figure~\ref{fig:two8} are Legendrian isotopic.

We now discuss the \dfn{Lagrangian projection} of a Legendrian knot. Let 
\[
\pi:\R^3\to \R^2: (x,y,z)\mapsto (x,y).
\]
Then the Lagrangian projection of a Legendrian knot $L$ is $\pi(L).$ The terminology ``Lagrangian
projection'' comes from the fact that $d\alpha|_{xy-\text{plane}},$ which is a symplectic form,
vanishes when restricted to $\pi(L).$ This is very important when considering Legendrian knots in
$\R^{2n+1}$ but is irrelevant in our discussion. However, we keep the terminology for the sake of
consistency. If we again parameterize $L$ by $\phi$ (all notation is as above) then $\pi(L)$ is
parameterized by $\phi_\pi(\theta) =(x(\theta), y(\theta)).$ Unlike the front projection, the
Lagrangian projection is always parameterized by an immersion, since if $x'(\theta)=y'(\theta)=0$
for some $\theta$ then $z'(\theta)\not=0$ so the tangent vector to $L$ is pointing in the
$\frac{\partial}{\partial z}$ direction which does not lie in $\xi.$

The Legendrian knot $L$ can be recovered (up to translation in the $z$-direction) from the
Lagrangian projection as follows: pick some number $z_0$ and define $z(0)=z_0.$ Then define
\beq\label{recz} 
z(\theta)=z_0+\int_0^{\theta} y(\theta)x'(\theta) d\theta. 
\eeq 
Since a Legendrian knot satisfies Equation~\eqref{mainleg} we see this equation can be written
\[
z(\theta)=z_0+\int_0^{\theta} z'(\theta)d\theta,
\]
which is a tautology. So the only ambiguity in recovering $L$ is the choice of $z_0.$

\begin{figure}[ht]
  {\epsfxsize=1.7in\centerline{\epsfbox{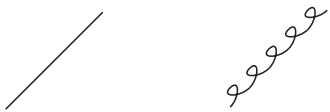}}}
  \caption{Using the Lagrangian projection to approximate
    topological arcs.}
  \label{fig:leglag}
\end{figure}

Let's observe a few restrictions on immersions $S^1\to \R^2$ that can be Lagrangian projections of a
Legendrian knot. First let $g:S^1\to \R^2:\theta\mapsto (x(\theta), y(\theta))$ be any immersion. If
we try to define $z(\theta)$ by Equation~\eqref{recz} then we run into problems. Specifically, if we
think of $\theta\in[0,2\pi]$ then $z(\theta)$ will be a well defined function on $S^1$ only if
$z(0)=z(2\pi).$ This condition can be written
\[
\int_0^{2\pi} y(\theta)x'(\theta) d\theta=0
\]
and of course is not satisfied for all immersions $g$. This is the only obstruction to lifting $g$
to an immersion $G:S^1\to \R^3$ whose image is tangent to $\xi.$ The image of $G$ will be a
Legendrian knot if $G$ is an embedding. The only way it can fail to be an embedding is if double
points in the image of $g$ lift to have the same $z$-coordinate. Thus an immersion $g$ lifts to a
Legendrian knot (well defined up to isotopy) if
\begin{enumerate}
\item $\int_0^{2\pi} y(\theta)x'(\theta) d\theta=0$ and
\item $\int_{\theta_0}^{\theta_1} y(\theta)x'(\theta) d\theta\not =0$ for all 
  $\theta_0\not= \theta_1$ with $g(\theta_0)=g(\theta_1).$
\end{enumerate}
Unfortunately these conditions are not easy to interpret diagrammatically, making working with
Lagrangian projections somewhat harder that working with front projections. None the less,
Lagrangian projections are still quite useful as we will see below. We also note that
Theorem~\ref{thm:approx} can be proven using Lagrangian projections. Specifically, given any
topological arc $\gamma:[0,1]\to \R^3$ we want to use the Lagrangian projection to prove
$\gamma([0,1])$ can be $C^0$ approximated by a Legendrian arc with the same end points. By
considering Figure~\ref{fig:leglag} one may clearly do this.

\begin{figure}[ht]
  {\epsfxsize=2in\centerline{\epsfbox{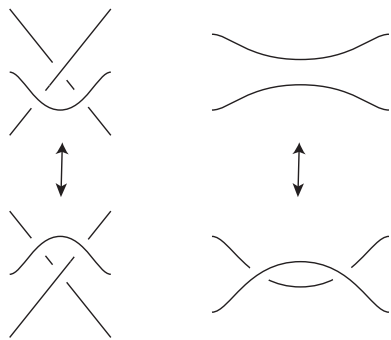}}}
  \caption{Legendrian Reidemeister moves in the Lagrangian
    projection. (Also need the corresponding figures rotated 180 degrees about all three coordinate
    axes.)}
  \label{fig:lagreid}
\end{figure}

We also have a weak Reidemeister type theorem.
\begin{thm}
  Two Lagrangian diagrams represent Legendrian isotopic Legendrian knots only if their diagrams are
  related by a sequence of moves shown if Figure~\ref{fig:lagreid}.
\end{thm}

Unfortunately these moves are not sufficient to guarantee Legendrian isotopy due to the integral
constraints discussed above.

In Section~\ref{FrontDGA} we will discuss converting a front projection into a Lagrangian
projection.

\subsection{Transverse knots}
A \dfn{transverse knot} $T$ in a contact manifold $(M^3,\xi)$ is an embedded $S^1$ that is always
transverse to $\xi:$
\[
T_xT\oplus\xi_x=T_xM,\quad x\in T.
\] 
If $\xi$ is orientable ({\em i.e.} if it is defined globally by a 1-form) then we can fix an
orientation on $\xi$. Recall our definition of contact structure requires a fixed orientation on
$M,$ thus if $T$ is transverse to $\xi$ we can orient $T$ so that it always intersects $\xi$
positively. When $T$ is so oriented we call $T$ a \dfn{positive transverse knot}. If we give $T$ the
opposite orientation then we call $T$ a \dfn{negative transverse knot}. If we do not put
positive/negative in front of the phrase ``transverse knot'' then positive is always implied. In
particular transverse knots are always oriented (unless $\xi$ is not orientable).

For the rest of this section we restrict attention to $(\R^3,\xi_{std}).$ Transverse knots are
usually studied via the analogy of the front projection: $\Pi:\R^3\to \R^2:(x,y,z)\mapsto (x,z).$
Unlike for Legendrian knots we cannot recover $T$ from $\Pi(T),$ however we can reconstruct $T$ up
to isotopy through transverse knots ({\em i.e.} we can recover the transverse isotopy class of $T$).
Thus it is reasonable to study $T$ through $\Pi(T).$

\begin{figure}[ht]
  {\epsfxsize=2in\centerline{\epsfbox{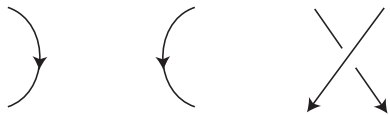}}}
  \caption{Segments excluded from projections of transverse
    knots.}
  \label{fig:tex}
\end{figure}

Let's examine $\Pi(T).$ Since the condition on a knot being transverse is an open condition in the
space of embeddings $S^1\to \R^3$ it is easy to check that for a generic transverse knot $\Pi(T)$
will be the image of an immersion. This immersion satisfies two obvious constraints:
\begin{enumerate}
        \item the immersion has no vertical tangencies pointing down and
        \item there are no double points as shown in Figure~\ref{fig:tex}.
\end{enumerate}
To check these conditions let $\phi:S^2\to \R^3:\theta\mapsto(x(\theta),y(\theta), z(\theta))$ be a
parameterization for $T.$ The fact that $T$ is a positive transverse knot implies that
\[
z'(\theta)-y(\theta)x'(\theta)>0.
\]
Now the projection of $T$ is parameterized by $\phi_\Pi(\theta)=(x(\theta), z(\theta)).$ At a
vertical tangency pointing down we have $x'=0$ and $z'<0,$ this contradicts the above equation thus
establishing condition (1). To see condition (2) note that
\[
y(\theta)< \frac{z'(\theta)}{x'(\theta)}.
\]
Thus the $y$ coordinate is bounded by the slope of the diagram in the $xz$-plane. (Once again recall
the positive $y$ axis points into the page.) One may now easily check that the picture on the right
side of Figure~\ref{fig:tex} cannot be the projection of a transverse knot.

\begin{figure}[ht]
  {\epsfxsize=2.8in\centerline{\epsfbox{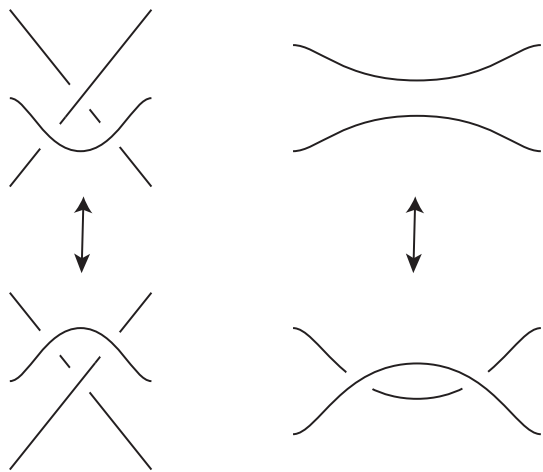}}}
  \caption{Transverse Reidemeister moves. Add arrows in all ways
    that don't violate the conditions in Figure~\ref{fig:tex}.
    (Also need the corresponding figures rotated 180 degrees about
    all three coordinate axes.)}
  \label{fig:tranreid}
\end{figure}

\begin{thm}[See \cite{EtnyreHonda01b, Swiatkowski}]
  Any diagram satisfying conditions (1) and (2) above can be lifted to a transversal knot in $\R^3$
  well defined up to isotopy through transversal knots. Moreover, two diagrams will represent the 
  same transverse isotopy class of transverse knots if any only if they are related by a sequence of
  move shown in Figure~\ref{fig:tranreid}.
\end{thm}

A second way to study transversal knots is by using closed braids. Recall a closed braid is simply a
knot (or link) in $\R^3$ (we will use cylindrical coordinates $(r,\psi, z)$ now) that can be
parameterized by a map $f:S^1\to \R^3:\theta\mapsto (r(\theta), \psi(\theta), z(\theta))$ for which
$r(\theta)\not=0$ and $\psi'(\theta)>0$ for all $\theta.$ For more on braids and closed braids see
\cite{BirmanBraid}.

To see the connection between braids and transverse knots we use the symmetric version of the
standard contact structure $(\R^3,\xi_{sym}).$ Since this contact structure is contactomorphic to
the standard one we can transfer any question about the contact structure $\xi_{std}$ to the contact
structure $\xi_{sym}.$ Now given a closed braid $B$ we can isotopy it through closed braid so that
it is far from the $z$-axis. Very far from the $z$-axis the planes that make up $\xi_{sym}$ are
almost vertical, that is close to the planes spanned by $\frac{\partial}{\partial z}$ and
$\frac{\partial}{\partial r}.$ Thus the closed braid type $B$ represents a transversal knot.

Bennequin proved the opposite assertion:
\begin{thm}[Bennequin 1983, \cite{Bennequin83}]
  Any transverse knot in $(\R^3, \xi_{sym})$ is transversely isotopic to a closed braid.
\end{thm}
The proof of this theorem is quite similar to the proof that knots can be braided
\cite{BirmanBraid}. One just needs the check that the ``braiding process'' can be done in a
transverse way. For details see \cite{Bennequin83, OS}.

\begin{figure}[ht]
  {\epsfxsize=1.2in\centerline{\epsfbox{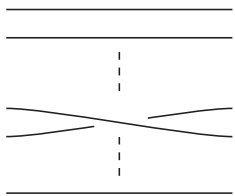}}}
  \caption{A generator $\sigma_i$ for the braid group $B_n.$}
  \label{fig:si}
\end{figure}

Recall, fixing $n$ points $p_i,$ in a disk $D^2$, an $n$-braid is an embedding of $n$ arcs
$\gamma_i:[0,1]\to D^2\times [0,1]$ so that $\gamma_i(t)\in D^2\times\{t\}$ and the endpoints of the
$\gamma_i$ corresponding to $0$ (resp. $1$) as a set map to $\{p_i\}_{i=1}^n$ in $D^2\times \{0\}$
(resp. $D^2\times\{1\}$). The set of all $n$-braids $B_n$ form a group. It is easy to see the group
is generated by $\sigma_i, i=1,\ldots n-1,$ where $\sigma_i$ is the $n$ braid with the $i$ and $i+1$
strands interchanging in a right handed fashion and the rest unchanged, see Figure~\ref{fig:si}.
The group $B_n$ naturally includes in $B_{n+1}.$ Given a braid $b$ in $B_n$ the \dfn{positive
stabilization} of $b$ is $b\sigma_n$ in $B_{n+1}.$
\begin{thm}[Orevkov and Shevchishin 2003, \cite{OS}; Wrinkle, \cite{Wrinkle}]
  Two braids represent the same transverse knot if and only if they are related by positive
  stabilization and conjugation in the braid group.
\end{thm}

\subsection{Types of classification}
In trying to classify Legendrian or transverse knots one could mean many different things. We
concentrate on Legendrian knots here but an analogous discussion holds for transverse knots too.

One can classify Legendrian knots up to {\em isotopy through Legendrian knots}. That is $L_0$ and
$L_1$ are Legendrian isotopic if there is a continuous family $L_t, t\in[0,1],$ of Legendrian knots
starting at $L_0$ and ending at $L_1.$ One can also classify Legendrian knots up to {\em ambient
contact isotopy}. Here we mean $L_0$ and $L_1$ are ambient contact isotopic if there is a one
parameter family $\phi_t:M\to M, t\in [0,1],$ (here $(M,\xi)$ is the ambient contact manifold) of
contactomorphisms of $M$ such that $\phi_0$ is the identity map and $\phi_1(L_0)=L_1.$
\begin{thm}
  The classification of Legendrian knots up to Legendrian isotopy is equivalent to the 
  classification up to contact isotopy. The analogous statement is true for transverse knots.
\end{thm}
\begin{proof}
The implication of contact isotopy to Legendrian isotopy is obvious. For the other implication
assume $L_t, t\in[0,1]$ is an isotopy through Legendrian knots in the contact manifold $(M,\xi).$
There is a family of diffeomorphisms $\phi_t:M\to M, t\in[0,1],$ such that $\phi_t(L_0)=L_t.$
Moreover, it is easy to arrange that $\phi_t^*(\xi|_{L_t})=\xi|_{L_0}.$ Let $\xi_t=\phi_t^*(\xi).$
This is a one parameter family of contact structures with $\xi_t=\xi_0$ along $L_0.$ Thus Gray's
Theorem, \cite{a, McDuffSalamon}, implies there is a family of diffeomorphisms $\psi_t$ such that
$\psi_t^*(\xi_t)=\xi_0$ and $\psi_t$ is the identity on $L_0.$ Now set $f_t=\phi_t\circ \psi_t.$
Note
\[
f_t^*(\xi_0)=\psi_t^*(\xi_t)=\xi_0.
\]
Thus $f_t$ are all contactomorphisms of $\xi_0=\xi.$ Moreover $f_t(L_0)=\phi_t(L_0)=L_t.$
\end{proof}

In $S^3$ with the standard contact structure there is another type of classification that is
equivalent to these.
\begin{thm}\label{compclass}
  In $(S^3,\xi_{std})$ (or $(\R^3, \xi_{std})$) two Legendrian knots are Legendrian isotopic if 
  and only if their compliments are contactomorphic.
\end{thm}
This theorem is not necessarily true in other contact manifolds.
\begin{proof}
A Legendrian knot $L$ has a canonical neighborhood $N(L).$ Denote by $M(L)$ the closure of the
complement of $N(L).$ If two Legendrian knots $L_0$ and $L_1$ in $(S^3, \xi_{std})$ have
contactomorphic complements, then let $\psi:M(L_0)\to M(L_1)$ be the contactomorphism. We can extend
$\psi$ over $N(L_0)$ so that it takes $N(L_0)$ to $N(L_1)$ and is a contactomorphism for
$(S^3,\xi_{std}).$ In \cite{Eliashberg92}, Eliashberg has shown that there is a unique tight contact
structure on $S^3.$ That coupled with Gray's Theorem implies there is a family of contactomorphisms
$\psi_t:S^3\to S^3, t\in[0,1]$ with $\psi_0$ the identity map and $\psi_1=\psi.$ Thus
$L_t=\psi_t(L_0)$ is an Legendrian isotopy from $L_0$ to $L_1.$
\end{proof}

\subsection{Invariants of Legendrian and transversal knots}\label{invtformula} 
Though it is not always essential we will always consider {\em oriented} knots.

\subsubsection{Classical invariants of Legendrian knots}
The most obvious invariant of a Legendrian knot is its underlying topological knot type, since any
Legendrian isotopy between two Legendrian knots is, among other things, a topological isotopy of the
underlying knots. Given a Legendrian knot $L$ we will denote its underlying topological knot
type by $k(L).$ We will use \K to denote a topological knot type ({\it i.e.} the set of topological knots
isotopic to a fixed knot) and the set of all Legendrian knots
$L$ with $k(L)\in\K$ will be denoted by $\L(\K).$

The next invariant of a Legendrian knot $L$ is the {\em Thurston--Bennequin invariant} which
intuitively measures the ``twisting of $\xi$ around $L$''. More rigorously this invariant is defined
by a trivialization of the normal bundle $\nu$ of $L.$
A fixed identification of $\nu$ with $L\times\R^2$ is called a
\dfn{trivialization} of $\nu$ or a \dfn{framing} of $L$.  A Legendrian
knot has a canonical framing: since $\xi_x$ and $\nu_x$ intersect
transversally ($T_xL\subset \xi_x$) one gets a line bundle
$l_x=\xi_x\cap \nu_x,$ for $x\in L.$ The line bundle $l$ gives a
framing of $\nu$ over $L.$ This framing is the
\dfn{Thurston--Bennequin framing} of $L$ and is denoted $tbf(L).$

If the normal bundle has a preassigned framing $\mathcal{F}$ then we can assign
a number to the Thurston--Bennequin framing of $L.$ This number $tw(L,\mathcal{F})$ is just the 
twisting of $l$ with respect to $\mathcal{F}$ and is called the twisting of $L$ with respect to
$\mathcal{F}.$ If $L$ is null homologous then $L$ has a framing given by a Seifert surface.
The twisting of $L$ with respect to this Seifert framing will be called the 
\dfn{Thurston-Bennequin invariant} of $L$ and is denoted $tb(L).$ One may alternately define $tb(L)$
as follows: let $v$ be a non-zero vector field along $L$ in $\nu\cap \xi$ and let $L'$ be a copy of $L$
obtained by pushing $L$ slightly in the direction of $v.$ Now define $tb(L)$ as the linking
of $L$ with $L',$ {\em i.e.} $tb(L)= lk(L,L').$ 
\bbr\label{trandeftb}
If $v'$ is a nonzero vector field along $L$ transverse to $\xi$ and $L''$ is obtained
from $L$ by pushing $L$ slightly in the direction of $v'$ then $tb(L)=lk(L,L'').$
\eer

The last ``classical'' invariant of a Legendrian knot $L$ is the {\em rotation number}
of $L.$ This invariant will only be defined for null homologous knots, so assume that $L=\partial 
\Sigma$ where $\Sigma$ is an embedded orientable surface. 
The contact planes when restricted to $\Sigma,$ $\xi|_\Sigma,$ form a trivial two dimensional
bundle (any orientable two plane bundle is trivial over a surface with boundary). This trivialization
of $\xi|\Sigma$ induces a trivialization $\xi|_L=L\times \R^2.$ Recall $L$ has an orientation, so
let $v$ be a non zero vector
field tangent to $L$ pointing in the direction of the orientation on $L.$ The vector field $v$ is
in $\xi|L=L\times \R^2$ and thus using this trivialization we can think of $v$ as a path of non zero
vectors in $\R^2,$ as such it has a winding number. This winding number $r(L)$ is the 
rotation number of $L.$  Note the rotation number depends on the orientation of $L$ and changes sign
if the orientation is reversed.
\bbr
One can also define $r(L)$ as the obstruction to extending $v$ to a non-zero vector field
in $\xi|\Sigma.$ 
\eer

The invariants $k(L), tb(L)$ and $r(L)$ will be called the {\em classical invariants} of $L.$ 

\subsubsection{Computation of the classical invariants via projection}
We now consider Legendrian knots $L$ in the standard contact structure on $\R^3$ and
interpret the classical invariants in the front and Lagrangian projection of $L$.
We begin with the rotation number. Let $w=\frac{\partial}{\partial y}.$ This is a nonzero section
of $\xi$ and thus can be used to trivialize $\xi|_L$ independent of finding a Seifert surface for 
$L.$ Now to compute the rotation number of $L$ we just need to see how many times a nonzero 
tangent vector field $v$ to $L,$ thought of as living in $\R^2$ via the trivialization given
by $w,$ winds around the origin of $\R^2.$ This is equivalent to making a signed count 
of how many times $v$ and $w$ point in the same direction, we call this an intersection of $v$ and $w.$ 
The ``sign of the intersection'' is determined by whether $v$ passes $w$ counterclockwise ($+1$) or
clockwise ($-1$). Now in the front projection $v$ will be pointing in the direction of 
$\pm w=\pm\frac{\partial}{\partial y}$ at the cusps. One may easily check that the intersection will be
positive when going down a cusp and negative when going up a cusp. Moreover, since this  counts
the number of times $v$ ``intersects'' $\pm w$ we need to divide by two in order to get $r(L).$ 
Thus in the front projection
\begin{equation}
r(L)=\frac12(D-U),
\end{equation}
where $U$ is the number of up cusps in the front projection and $D$ is the number of down cusps.
In the Lagrangian projection
$w$ projects to $\frac{\partial}{\partial y}$ and thus the rotation number is easily seen to
be the winding number of the tangent vectors to $\pi(L)$ 
\begin{equation}
r(L)= \text{ winding} (\pi(L)).
\end{equation}

\begin{figure}[ht]
  {\epsfxsize=1.3in\centerline{\epsfbox{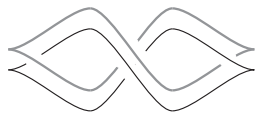}}}
  \caption{Knots $L$ (black) and $L'$ (grey) used to compute $tb(L).$}
  \label{fig:ftb}
\end{figure}

\begin{figure}[ht]
  \relabelbox 
  \small {\epsfxsize=1.8in\centerline{\epsfbox{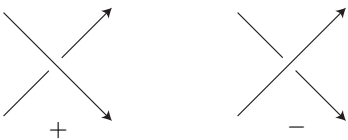}}}
  \relabel {+}{$+$}
  \relabel {-}{$-$} 
  \endrelabelbox
  \caption{Right handed crossing (on left) and a left handed crossing (on right). Note this
    figure is to indicate that right (left) crossings contribute $+$($-$) to the write of a
    knot. The picture on the right cannot occur in a front projection, but non the less left
    handed crossings can.}
  \label{fig:writhe}
\end{figure}

For the Thurston--Bennequin invariant. Let $v=\frac{\partial}{\partial z}$ then
for any Legendrian knot $L,$  $v$ is a vector field transverse to $\xi$ along $L$ (as
in Remark~\ref{trandeftb}). Thus $tb(L)$ is the linking of $L$ with a copy $L'$ of $L$ obtained by
shifting slightly in the $z$ direction. 
So in the front projection $L$ and $L'$ are as in Figure~\ref{fig:ftb}.
Now the linking number of $L$ and $L'$ is just one half
the signed count of the intersections between them. Where an intersection is positive if it is
right handed and negative if it is left handed (Recall $L$ and hence $L'$ are oriented). 
See Figure~\ref{fig:writhe}.
One can see that at each right (left) handed self crossing of the projection of $L$
there will be two right (left) handed crossings of $L$ and $L'$ and at a right or left cusp of $L$
there will be a left handed crossing of $L$ and $L'.$ Thus in the front projection
\begin{equation}
        tb(L)=\text{writhe}(\Pi(L))-\frac12(\text{number of cusps in $\Pi(L)$}).
\end{equation}
In the Lagrangian projection $L$ and $L'$ project to the same diagram, but we have
\begin{equation}\label{lagtb}
        tb(L)=\text{writhe}(\pi(L)).
\end{equation}
One may prove this formula by first
recalling $tb(L)$ is the linking number of $L$ and $L'.$ In the Lagrangian projection think
about trying to pull $L'$ straight up. If during this process $L'$ never intersects $L$ then
their linking is $0.$ Moreover, each crossing in the diagram for
$L$ contributes $\pm 1$ to the linking number.

\subsubsection{Classical invariants of transverse knots}
For transverse knots $T$ there are only two classical invariants, the topological knot type
$k(T)$ and the \dfn{self-linking number} (this is sometimes called the Bennequin number).
We denote transverse knots in a fixed topological knot type \K by $\mathcal{T}(\K).$ To define
the self-linking number of $T$ we assume it is homologically trivial. (There is a ``relative'' version
of the self-linking number if this is not true but we will not discuss this point here.) Thus
there is an orientable surface $\Sigma$ such that $\partial \Sigma=T.$ As above we know that 
$\xi|_{\Sigma}$ is trivial so we can find a nonzero vector field $v$ over $\Sigma$ in $\xi.$ Let
$T'$ be a copy of $T$ obtained by pushing $T$ slightly in the direction of $v.$ The self-linking 
number $sl(T)$ of $T$ is the linking of $T'$ with $T.$ 

The self-linking number also has an interpretation in terms of a relative Euler class. To see
this let $v$ be a non zero vector in $\xi\cap T\Sigma$ along $T$ that points out of $\Sigma.$ 
If this $v$ were to extend over $\Sigma$ then we could use it to get $T'$ above and the 
self-linking would be 0. If $v$ does not extend then the self-linking is not 0.
\bbr
The self-linking number $l(T)$ is precisely the obstruction to 
extending $v$ over $\Sigma$ to a non-zero vector field in $\xi.$ 
\eer

\subsubsection{Computations of the self-linking number}
To compute $sl(T)$ for a transverse knot in the standard contact structure on $\R^3$ note that
the vector $v=\frac{\partial}{\partial y}$ is always in $\xi_{std}$ and thus can be used to trivialize
$\xi_{std}$ independent of a Seifert surface for $T.$ Now let $T'$ be a copy of $T$ obtained by 
pushing $T$ along $v.$ The self-linking number of $T$ is the linking number of $T$ and $T'.$ If we 
consider $T$ and $T'$ in the ``front projection'' we see that
\beq
        sl(T)=\text{writhe}(\Pi(L)).
\eeq
The argument for this formula is exactly like the one for Equation~\eqref{lagtb}.

\begin{figure}[ht]
        {\epsfxsize=2in\centerline{\epsfbox{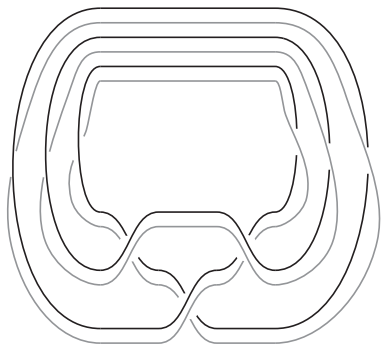}}}
        \caption{The braid $T$ in black and its push off by $v$ in grey.}
        \label{fig:bsl}
\end{figure}

We would also like to see a formula for $sl(T)$ in terms of a braid representation for $T.$
To do this recall we must use the cylindrically symmetric contact structure $\xi_1$ on $\R^3.$
For our non zero vector field we take $v=\frac{\partial}{\partial r}$ on $\{x\geq \epsilon\},$ 
$v=-\frac{\partial}{\partial r}$ for $\{x\leq -\epsilon\}$ and on $\{-\epsilon<x<\epsilon\}$ we 
have interpolate between these two choices by rotating clockwise in the contact planes. 
If we think of 
all the braiding as occurring in $\{x\geq \epsilon\}$ the we see  each generator
$\sigma^{\pm}$ contributes $\pm 1$ to the self-linking number. See Figure~\ref{fig:bsl}.
Moreover
each strand as it passes in and out of $\{x\leq \epsilon\}$ contributes $-1$ to the
self-linking number. Thus we have
\beq
        sl(T)=a(T)-n(T),
\eeq
where $n(t)$ is the number of strands in the braid $b$ representing $T$ and
$a(T)$ is the algebraic length of $b$ when written in terms of the generators
$\sigma_i$ (algebraic length is the sum of the exponents on the generators).
This useful equality was first observed by Bennequin \cite{Bennequin83}.

\subsection{Stabilizations}\label{section:stab}
\begin{figure}[ht]
  \relabelbox 
  \small {\epsfxsize=2.2in\centerline{\epsfbox{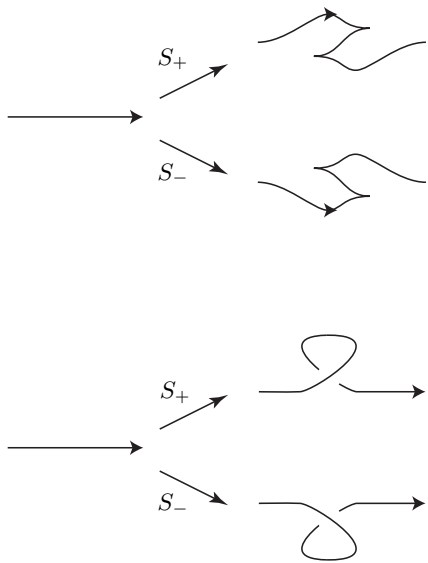}}}
  \relabel {Sp}{$S_+$}
  \relabel {Sn}{$S_-$}
  \relabel {Sp1}{$S_+$}
  \relabel {Sn1}{$S_-$}
  \endrelabelbox
        \caption{Stabilizations in the front projection (top) and the Lagrangian projection (bottom).}
        \label{fig:stab}
\end{figure}

Given a Legendrian knot $L$ there is a simple way to get another Legendrian knot in the
same topological knot type: stabilization. We describe \dfn{stabilization} in the standard contact 
structure on $\R^3.$ 
If a strand of $L$ in the front projection of $L$ is as on the left hand
side of Figure~\ref{fig:stab} 
then the stabilization of $L$ is obtained by removing the stand and replacing
it with one of the zig-zags shown on the right hand side of Figure~\ref{fig:stab}. If down cusps are
added then the stabilization is called positive, and denoted $S_+(L),$ 
and if up cusps are added then the stabilization is called negative, and denoted $S_-(L).$ 
Since stabilizations are done locally, this actually defines, via Darboux's 
Theorem, stabilizations for Legendrian knots
in any manifolds. Note
\begin{equation}
        tb(S_\pm(L))=tb(L)-1,
\end{equation}
and
\begin{equation}
        r(S_\pm(L))=r(L)\pm 1.
\end{equation}
It is important to observe that stabilization is a well defined operation, that is, it does
not depend at what point the stabilization is done. This is true in general 
\cite{EtnyreHonda01b} but in $\R^3$ it is somewhat easier to prove \cite{FuchsTabachnikov97}. 
One only needs to check that the ``zig-zags''
used to define the stabilization can be move past cusps and crossings in the
front projection. Using Theorem~\ref{thm:reidfront} this is a fun exercise.

One important fact concerning stabilizations is the following.
\begin{thm}[Fuchs and Tabachnikov 1997, \cite{FuchsTabachnikov97}]\label{stablist}
Given two Legendrian knots $L_1$ and $L_2$ in $(\R^3, \xi_{std})$ that are topologically isotopic then
after each has been stabilized some number of times they will be Legendrian isotopic.
\end{thm}
So the topological classification of knots is equivalent to the stable classification of
Legendrian knots. 
To prove this theorem one checks that
the topological Reidemeister moves can be performed via Legendrian Reidemeister
moves after sufficiently many stabilizations. 
It should not be difficult to prove this theorem in any contact 3-manifold,
but such a proof does not exist in the literature.

Note using stabilization it is easy to get Legendrian knots with arbitrarily negative Thurston-Bennequin
invariants. In any tight contact structure there is an upper bound on the Thurston-Bennequin invariant
(see Section~\ref{sec:bounds}). Thus while it is easy to stabilize a knot it is not necessarily easy to ``destabilize'' 
a knot. We say $L$ destabilizes if there is a Legendrian knot $L'$ such that $L=S_\pm(L').$ 

\begin{figure}[ht]
        {\epsfxsize=2.8in\centerline{\epsfbox{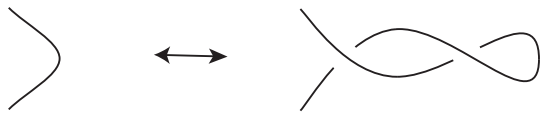}}}
        \caption{Transverse stabilization.}
        \label{fig:tstab}
\end{figure}

There is notion of transverse stabilization too. The stabilization of a transverse knot
$T$ is formed by taking an arc in the front projection as shown on the left hand side
of Figure~\ref{fig:tstab} and replacing it with the arc on the right hand side. Denote the resulting
knot by $S(T).$
Note there is only one type of stabilization here and
\begin{equation}
        sl(S(T))=sl(T)-2.
\end{equation}
Once again it is not hard to show that stabilization is well defined for transverse knots. 
We also have
\begin{thm}[Fuchs and Tabachnikov 1997, \cite{FuchsTabachnikov97}]
Given two transverse knots $T_1$ and $T_2$ in $(\R^3,\xi_{std})$ that are topologically isotopic then
after each has been stabilized some number of times they will be transversely isotopic.
\end{thm}

There is a notion of destabilization for transverse knots that is exactly analogous to destabilization
of Legendrian knots.

\subsection{Surfaces and the classical invariants}\label{sfcandinvt}
In this section we see how to compute the classical invariants of a Legendrian and transverse
knot in terms of the characteristic foliation on a surface bounded by the knot. We will also
see how convex surface theory \cite{Giroux91, Honda1} 
can be used to understand the classical invariants of Legendrian knots.

Let $T$ be a transverse knot bounding the surface $\Sigma.$ Orient $\Sigma$ so that
$T$ is its oriented boundary. With this orientation the characteristic foliation $\Sigma_\xi$
points out of $\partial \Sigma=T$ (recall the characteristic foliation is oriented).
Perturb $\Sigma$ so that $\Sigma_\xi$ is generic. In particular, we can assume that all the
singularities are isolated elliptic or hyperbolic points. Moreover, each singularity has
a sign depending on whether the orientation of $\xi$ and $T\Sigma$ agree at the singularity.
Let $e_\pm$ be the number of $\pm$ elliptic singularities in $\Sigma_\xi$ and let $h_\pm$
be the number of $\pm$ hyperbolic singularities.
We now have
\beq\label{lformula} 
-sl(T)= (e_+-h_+) - (e_--h_-).
\eeq
This formula follows easily by interpreting $sl(T)$ as a relative Euler class (that is let
$v$ be a vector field along $T$ tangent to $\Sigma_T$ and contained in $\xi$ and pointing
into $\Sigma_T,$ then $sl(T)$ is the obstruction to extending $v$ to a nonzero vector field 
on $\Sigma_T$).

\begin{figure}[ht]
  \relabelbox 
  \small{\epsfxsize=2in\centerline{\epsfbox{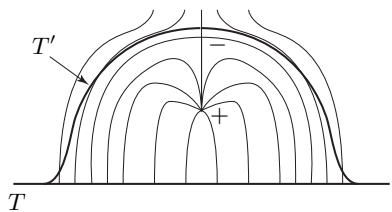}}}
  \relabel {T'}{$T'$}
  \relabel {T}{$T$}
  \relabel {+}{$+$}
  \relabel {-}{$-$}
  \endrelabelbox
  \caption{Using the characteristic foliation to recognize a destabilization.}
  \label{fig:tpass}
\end{figure}

By looking at the characteristic foliation one may frequently see how to destabilize a transverse
knot. In particular, if one sees Figure~\ref{fig:tpass}
in the characteristic foliation
of a surface with a transverse knot $T$ in its boundary, then one can let $T'$ be the knot shown
in the picture and it is clear that $sl(T')=sl(T)+2.$ With a little work one can, in fact, see that
$S(T')=T.$ Thus when one finds a negative hyperbolic point
on a surface whose unstable manifolds separate off a disk with one elliptic point (necessarily
positive), then the knot can be destabilized. Such a disk is called a \dfn{transverse bypass} and
first appeared in \cite{Etnyre00} and was used in the classification of some transverse torus knots in
\cite{Etnyre99}.

Now consider a Legendrian knot $L$ with Seifert surface $\Sigma.$ 
It is relatively easy to see we 
may isotopy $\Sigma$ relative to its boundary so that the singularities in the characteristic
foliation along
the boundary alternate in sign. When this is arranged it is easy to see
$|\tb(L)|$ is simply half the number of singularities along the boundary.
We now observe:
\begin{lem}[Eliashberg and Fraser 1998, \cite{EliashbergFraser}]
Assume the singularities in the characteristic foliation along $\partial \Sigma=L$ alternate in sign.
If $tb(L)\leq 0$ then positive/negative singularities along the boundary will be 
sources/sinks (along the boundary). If $tb(L)>0$ the positive/negative singularities
along the boundary will be sinks/sources (along the boundary).
\end{lem}
The idea for this lemma is simply to examine how the characteristic 
foliation inherits an orientation from the contact planes and the Seifert surface.
A singularity along the boundary is positive or negative depending on whether the
contact planes are twisting past $\Sigma$ in a right handed or left handed fashion.

Recall that elliptic singularities in the characteristic foliation are sources (sinks)
when the singularity is positive (negative). Thus when $tb(L)>0$ we see from this lemma
that all the singularities must be hyperbolic! So the characteristic foliation can be 
standardized along the boundary when $tb(L)>0.$ 

When $\tb(L)\leq 0$ it takes a little more work to normalize to characteristic foliation near
the boundary. In addition, the foliation can be normalized in various ways. In particular,
one can arrange that all the singularities are hyperbolic, or that all the singularities
are elliptic, or that they alternate between hyperbolic and elliptic \cite{EliashbergFraser}.

If $tb(L)\leq 0$ then we can isotope $\Sigma$ relative to $L$ so that it is convex.
Denote the dividing curves of $\Sigma$ by $\Gamma$ and the $\pm$ regions of $\Sigma\setminus \Gamma$
by $\Sigma_\pm.$ Then we have
\beq\label{tbforconvex}
tb(L)=-\frac12 (L\cap \Gamma),
\eeq
and
\beq
r(L)=\chi(\Sigma_+)-\chi(\Sigma_-).
\eeq
Call a properly embedded arc $a\subset \Gamma$ 
in $\Sigma$ boundary parallel if the closure of one of the components of 
$\Sigma\setminus a$ is a disk that contains no components of $\Gamma$ in its interior. Such
a boundary parallel dividing curve is frequently called a \dfn{bypass}, because it allows 
one to ``bypass'' some twisting as the following lemma shows. Bypasses were first introduced in
\cite{Honda1}.
\begin{lem}
If the dividing curves $\Gamma$ of $\Sigma$ contain a boundary parallel arc and $\tb(L)<1$ or
$\Sigma$ has genus greater than 0, then $L$ can be destabilized.
\end{lem}
\begin{proof}
The Legendrian Realization Principle \cite{Giroux91, Honda1} says that we can alter the characteristic
foliation to any other singular foliation as long as $\Gamma$ ``divides'' this foliation.
Thus we may realize the foliation 
shown in Figure~\ref{fig:lpass} 
\begin{figure}[ht]
  \relabelbox 
  \small{\epsfxsize=2.8in\centerline{\epsfbox{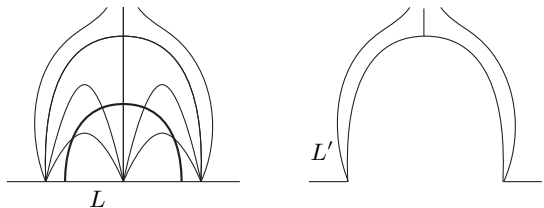}}}
  \relabel {L'}{$L'$}
  \relabel {L}{$L$}
  \endrelabelbox
  \caption{Characteristic foliation associated to a boundary parallel dividing curve (left).
    The destabilization $L'$ of $L$ (right).}
  \label{fig:lpass}
\end{figure}  
near the boundary parallel arc. We then let $L'$ be the Legendrian
knot indicated on the right hand side of Figure~\ref{fig:lpass} (after the corners are
smoothed). It is not hard to show that $L=S_\pm(L').$ 
\end{proof}
\bbr
A bypass is technically the disk cobounded by $L$ and $L'$ but using this lemma we see that
this is essentially equivalent to the definition above.
\eer

\subsection{Relation between Legendrian and transversal knots}
\begin{figure}[ht]
  \relabelbox 
  \small{\epsfxsize=2.8in\centerline{\epsfbox{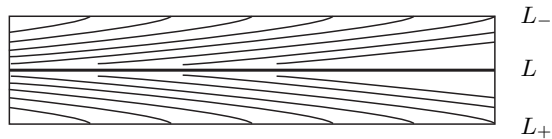}}}
  \relabel {m}{$L_-$}
  \relabel {L}{$L$}
  \relabel {p}{$L_+$}
  \endrelabelbox
  \caption{Annulus involved in transverse push-off of $L$.}
  \label{fig:tpushoff}
\end{figure}

Let $L$ be a Legendrian knot in a contact manifold $(M,\xi).$ Let $A=S^1\times[-1,1]$ 
be an embedded annulus in $M$ such that $S^1\times\{0\}=L$ and $A$ is transverse to $\xi.$
If $A$ is sufficiently ``thin'' then the characteristic foliation on $A$ is as shown in 
Figure~\ref{fig:tpushoff}.
Given such an $A$ the curve $L_+=S^1\times\{+ \frac12\}$ (resp. $L_-=S^1\times\{- \frac12\}$) is a positively
(resp. negatively) transverse knot. The knots $L_\pm$ are called the positive (resp. negative)
transverse push off of $L.$ It is easy to check that any two transversal push-offs are transversely
isotopic.
\begin{warning}
The definition of the positive and negative push-off is not standard in the literature. Some
authors ({\em eg.} Bennequin \cite{Bennequin83}) reverse the naming of positive and negative push-offs.
\end{warning}

\begin{figure}[ht]
        {\epsfxsize=3.5in\centerline{\epsfbox{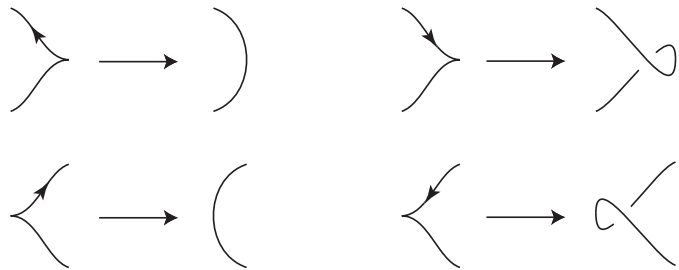}}}
        \caption{Positive transverse push-off in the front projection.}
        \label{fig:ltot}
\end{figure}

In the standard contact structure on $\R^3$ we can see this transverse push off in the front projection.
Let $L$ be an oriented Legendrian knot and consider its front projection. Away from the cusps, any
strand oriented to the right should be pushed in the $-y$-direction slightly to get $L_+$ any strand
oriented to the left should be pushed in the $y$ direction slightly. Thus way from the cusps the
front projection of $L$ and $L_+$ agree. A careful analysis near the cusps shows that they change as
shown in Figure~\ref{fig:ltot}.

Given a Legendrian knot $L$ we can construct a positive and a negative transversal 
knot $T_\pm(L)$ the relation between their classical invariants is described in the following 
lemma.
\begin{lem} 
        The invariants of Legendrian knots and their transverse push offs are related by
        \beq\label{invtrel} sl(T_\pm(L))= tb(L) \mp r(L).\eeq
\end{lem}
\begin{proof}
Given $L$ a Legendrian knot and $\Sigma$ a Seifert surface for $L$ we let 
$\tau$ be a nonzero vector in $\xi|_L$ that extends to a nonzero vector in $\xi|_\Sigma.$
Thus if $v$ is a nonzero vector field tangent to $L$ then $r(L)$ is the twisting of $v$ relative
to $\tau$ in $\xi$ which we denote: $r(l)=t(v,\tau,\xi|_L).$ If $\nu$ is the normal bundle to $L,$ 
$s$ is the outward pointing normal to $\Sigma$ and $w$ is a vector field in $\xi|_L$ transverse to $v$ then
$tb(L)=t(w,s,\nu).$ 

Let $v_+, w_+, \tau_+$ be the vector fields in $\xi|_{L_+}$ coming naturally from the definition of $L_+.$
Note we can think of $s$ as the same as $s_+$ since $L$ and $L_+$ are the same topological knots.
Now keeping in mind $\nu=\xi|_{L_+}$ we have 
\begin{align*}
sl(L_+)&= t(\tau_+, s, \nu)=t(\tau_+, w_+, \nu)+t(w_+,s,\nu)\\
&=t(\tau, w, \xi)+ t(w, s, \nu)= -r(L)+tb(L).\end{align*}
The equation for $L_-$ is similar except that $\nu=-\xi$ so $t(\tau_+, w_+, \nu)=-t(\tau, w, \xi).$

For Legendrian knots in $(\R^3,\xi_{std})$ there is a direct diagrammatic proof of these
formulas using the formulas derived in the previous sections and Figure~\ref{fig:ltot}.
\end{proof}

There is also a Legendrian push off of a transversal knot. Let $T$ be a transverse knot. We want
to construct a standard model for a neighborhood of $T.$ To this end consider $(\R^3,\xi_{sym}).$
Recall $\xi_{sym}=\ker(dz+r^2d\phi).$ Now let $M$ be $\R^3$ modulo the action $z\mapsto z+1.$ The
contact structure $\xi_{sym}$ clearly induces a contact structure on $M$ (also denoted $\xi_{sym}$).
Now $M=S^1\times \R^2$ and $T_a=\{r=a\}\subset M$ is a torus that bounds the solid torus
$S_a=\{r\leq a\}\subset M.$ The characteristic foliation of
$T_a$ is by lines of slope $a^2.$ A standard application of Moser's technique 
shows that the transverse knots $T$ has a neighborhood $N$ contactomorphic to $S_a$ for some $a.$
Thus in a neighborhood of $T$ there is a torus $T_b$ with $b<a$ and $b^2=\frac1n$ for some positive
integer $n.$ The characteristic foliation on $T_b$ is by lines of slope $b^2=\frac1n.$ Let $T_l$ be
a leaf in the foliation of $T_b.$ It is easy to see that $T_l$ is a Legendrian knot
topologically isotopic to $T.$ The Legendrian knot $T_l$ is called a \dfn{Legendrian push off} of $T.$
By considering the obvious annulus between $T$ and $T_l$ it is easy to see that
the positive transverse push off of $T_l$ is $T.$ 

Notice that $T_l$ is called {\em a} Legendrian push off and not {\em the} Legendrian push off. The reason
for this is that there are many other Legendrian push offs. Note that we can find a $c<b<a$ such that
$c^2=\frac{1}{n+1}.$ The the Legendrian push off obtained using $T_c,$ denote it $T_l',$ 
is related to $T_l$ by a negative stabilization \cite{EtnyreHonda01b}: $T_l'=S_-(T_l).$ 
We say two Legendrian
knots are \dfn{negatively stably isotopic} if they are Legendrian isotopic
after each has been negatively stabilized some number of times. Thus the Legendrian push off of 
a transverse knot gives a well defined negative stable isotopy class of Legendrian knots.
\begin{thm}[In $(\R^3, \xi_{std}),$ Epstein, Fuchs and Meyer 2001, \cite{EFM}: In
a general contact manifold, Etnyre and Honda 2001, \cite{EtnyreHonda01b}]\label{nslist}
Two Legendrian knots are negatively stably isotopic if and only if their transverse push offs
are transversely isotopic.
\end{thm}

It is quite easy to show the ``only if'' part of this theorem.
The ``if'' part is much more difficult.

Since we know all transverse knots are the transverse push off of some Legendrian
knot, the classification of transversal knots is equivalent to the classification of Legendrian
knots up to negative stabilization.


\section{Tightness and bounds on invariants}\label{sec:bounds}
In the standard contact structure on $\R^3,$ or more generally in a tight contact
structure, there are many bounds on the classical invariants of a Legendrian or
transversal knot. In fact, the existence of some of these bounds
is at the heart of the nature of tightness. In this section we discuss these inequalities. 

\subsection{Bennequin's inequality}
Part of the origins of modern contact geometry can be found in the following theorem.
\begin{thm}[Eliashberg 1992, \cite{Eliashberg92}]\label{ben}
        Let $(M,\xi)$ be a tight contact 3-manifold. Let $L$ be a Legendrian
        knot in $M$ with Seifert surface $\Sigma_L$ and 
        $T$ be a transverse knot with Seifert surface $\Sigma_T.$ Then
        \begin{equation}\label{slbbound}
                sl(T)\leq -\chi(\Sigma_T),
        \end{equation} 
        and
        \begin{equation}\label{tbbound}
                tb(L)+|r(L)|\leq -\chi(\Sigma_L).
        \end{equation} 
\end{thm}
\begin{proof}[Sketch of Proof]
Clearly the second inequality follows from the first and Equation~\eqref{invtrel}. To prove the first 
inequality recall the discussion in Section~\ref{sfcandinvt}: if the foliation on $\Sigma$ is generic
then
\[-sl(T)= (e_+-h_+)-(e_--h_-).\]
Moreover one may easily show
\[\chi(\Sigma_T)=(e_++e_-)-(h_++h_-).\]
Thus 
\[sl(T)+\chi(\Sigma_T)=2(e_--h_-).\]
So we can prove Inequality~\ref{slbbound} by showing that we may isotop $\Sigma_T$ (rel boundary)
so that $e_-=0.$ This is done by showing that the tightness of $\xi$ implies that 
each positive elliptic point is connected to
a positive hyperbolic point and then canceling this pair. 
This argument is discussed in \cite{Eliashberg92, EtnyreIntro}.  
\end{proof}

This bound is not sharp for many knot types. This is easily seen once other inequalities 
are established, so we defer this discussion until later in this section.

The inequalities in Theorem~\ref{ben} are called the Bennequin inequalities. However, the
theorem as stated is due to Eliashberg \cite{Eliashberg92}. 
What Bennequin actually showed was that for any transverse
knot in $(\R^3,\xi_{std})$ Equation~\eqref{slbbound} holds \cite{Bennequin83}. 
One reason for this is that the notion of
tight vs.~overtwisted was unknown at the time. Bennequin was trying to show that $\R^3$ 
has more than one contact structure. He did this by showing that in the standard contact
structure Equation~\eqref{slbbound} is true while it is not in the contact structure $\xi_{ot}$
defined in Example~\ref{otex}.
Specifically the Legendrian unknot
$L=\{z=0, r=\pi\}$ in $(\R^3,\xi_{ot})$ is easily seen to have $tb(L)=0.$  

This was one of the first indications of the existence of a tight vs.~overtwisted dichotomy.
In fact, tightness can be characterized in terms of knots.
\begin{thm}
        A contact structure $\xi$ is overtwisted if and only if there is a Legendrian unknot
        with Thurston--Bennequin invariant equal to 0 if and only if there is a transverse 
        unknot with self-linking number equal to 0.
\end{thm}
Using Theorem~\ref{ben} the only non-trivial (but still easy) part of this 
theorem is that overtwisted implies
the existence of a transverse unknot with self-linking 0.
%

Since Equation~\ref{tbbound} gives an upper bound on the Thurston--Bennequin invariant
of a knot $L$ in a tight contact structure we can make the following definition
\[\overline{tb}(\K)=\text{max}\{tb(L)| L\in \L(\K)\}.\]
So $\overline{tb}(\K)$ is the maximal Thurston-Bennequin invariant for Legendrian knots
in the knot type $\K.$ Clearly this is an invariant of the topological knot type.

\subsection{Slice genus}\label{sec:sg}

There has been a refinement of Bennequin's inequality for knots in Stein 
fillable contact structures. For a discussion of Stein fillable contact structures 
see \cite{EtnyreIntro, Gompf98}. The important point here is that the standard
tight contact structure on $S^3$ is Stein fillable by $B^4$ (with its standard complex
structure).
\begin{thm}[Lisca and Matic 1998, \cite{LiscaMatic98}]
Let $X$ be a Stein manifold with boundary $M$ and $\xi$ the induced contact
structure. If $L$ is a Legendrian knot in $M$ then
\beq\label{sliceg}
tb(L)+|r(L)|\leq -\chi (\Sigma),
\eeq
where $\Sigma$ is an embedded oriented surface in $X$ with $\partial \Sigma=L.$ 
\end{thm}

Given a knot type $\K$ in $M=\partial X,$
let $g_s(\K)$ denote the minimal genus of an 
embedded orientable surface $\Sigma$ in $X$ with $\partial \Sigma=\K.$ This
is called the \dfn{slice genus} of $\K$ in $(X,M).$  Equation~\eqref{sliceg} can be
interpreted in terms of the slice genus:
\[\frac12(tb(L)+|r(L)|+1)\leq g_s(\K)\]
for all $L\in\L(\K).$  Clearly the
genus of a knot is an upper bound on the slice genus, thus this inequality implies the
Bennequin inequality.

For Legendrian knots in the standard contact structure in $S^3$ this theorem is due
to Rudolph \cite{Rudolph}. In full generality is appears in \cite{LiscaMatic98}. Currently,
the most direct proof of this theorem involves a standard adjunction type inequality in
Seiberg-Witten theory and an embedding theorem for Stein manifolds proved
in \cite{LiscaMatic97}.

While much of the deep mathematical content to the above theorem is Seiberg-Witten theory, the language
of Legendrian knots provides a user friendly ``front end'' for the Seiberg-Witten theory. In particular, in
many examples it is easy to simply draw pictures to see knots cannot be slice or to get bounds on their
slice genus without ever having to explicitly invoke Seiberg-Witten theory. We will see applications of this
in Section~\ref{app:concord}. Here we
illustrate the usefulness of this theorem with a simple example.
Let $\K$ be a knot type and $\overline{\K}$ be its mirror image. It is well known \cite{Rolfsen} that
$\K\#\overline{\K}$ is slice. Thus the maximal Thurston--Bennequin invariant satisfies
$\overline{tb}(\K\#\overline{\K})\leq -1.$
But later we will see (Section~\ref{sec:consum}) that $\overline{tb}(\K\#\overline{\K})=\overline{tb}(\K)+
\overline{tb}(\overline{\K})+1.$ Thus 
\[\overline{tb}(\overline{\K})\leq-\overline{tb}(\K)-2.\]
By drawing pictures (see Section~\ref{sec:torus}) 
one may easily see that for $p,q$ positive and relatively prime,
a $(p,q)$-torus knot $\K_{(p,q)}$
has $\overline{tb}(\K_{(p,q)})=pq-p-q.$ Thus 
\[\overline{tb}(\K_{(-p,q)}) = \overline{tb}(\overline{\K}_{(p,q)}) \leq p+q-pq-2,\]
and we see the Bennequin inequality is not sharp. (Later we will see that $\overline{tb}(\K_{(-p,q)})=-pq$.)

\subsection{Other inequalities in $(\R^3,\xi_{std})$}

\begin{figure}[ht]
        {\epsfxsize=2.5in\centerline{\epsfbox{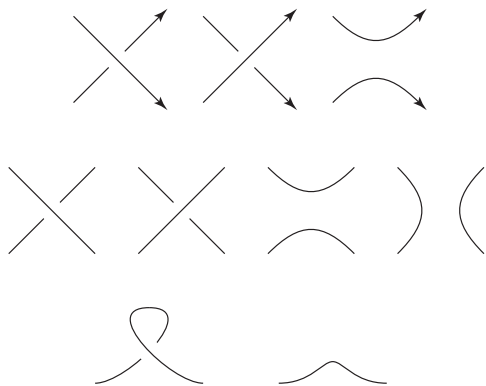}}}
        \caption{The relation between $K_+, K_-, K_0$ (top), $K_+,K_-,K_0,K_\infty$
        (middle) and $K_r, K$ (bottom).}
        \label{fig:skien}
\end{figure}

We begin by recalling the definition of several knot polynomials. Consider a projection of
an oriented knot $K.$  
To this projection we associate the HOMELY (Laurent) polynomial (aka, two variable
Jones polynomial) $P_K(v,z)$ which is defined by the
skien relation 
\[\frac{1}{v}P_{K_+} - vP_{K_-}=zP_{K_0},\]
where $K$ is either $K_+$ or $K_-$ and 
$K_+, K_-, K_0$ are related to each other as shown in the top row of Figure~\ref{fig:skien}.

Let $D$ be a regular projection
of a (non-oriented) topological knot $K.$ 
To this projection we associate a two variable (Laurent) polynomial
$L_K(a, x)$ defined by
\[L_U=\frac{a-a^{-1}}{x}+1,\]
\[L_{K_+}-L_{K_-}=x(L_{K_0}-L_{K_\infty}),\]
\[L_{K_r}= a L_K,\]
where $K_+, K_-, K_0, K_\infty$ are as related in the middle row of Figure~\ref{fig:skien}, $K_r, K$ are as
related in the bottom row of Figure~\ref{fig:skien} and $U$ is the diagram for the unknot with no crossings.
The polynomial $L$ is an invariant of $D$ through regular
homotopies. The polynomial $F_K(a, x)=a^{\text{writhe}(D)} H(a, x)$ is the Kauffman polynomial
of $K.$  It is an invariant of the topological knot $K.$

\begin{thm}\label{otherbounds} 
Let $\K$ be a topological knot type.
For any $L\in\L(\K)$ we have
\begin{enumerate}
\item $tb(L)+|r(L)|\leq d_{P_K}$  and
\item $tb(K)\leq d_{F_K}.$
\end{enumerate}
Where $d_{P_K}$ is the lowest degree in $v$ of $P_K(v,z)$ and $d_{F_K}$ is the
lowest degree in $a$ of $F_K(a,x).$
\end{thm}

The first inequality, observed by Fuchs and Tabachnikov \cite{FuchsTabachnikov97},
follows from work of Franks and Williams \cite{FranksWilliams87} and Morton \cite{Morton86} 
who showed that for an $n$ braid with algebraic length $a,$ 
\[a-n\leq d_{P_K},\]
and the fact we observed in Section~\ref{invtformula} that $sl=a-n.$ 
Thus this inequality gives a bound on the self-linking numbers of
transverse knots. The relation between the invariants of a Legendrian knot and its
transverse push off complete the proof of the first inequality in the theorem. 

The second inequality is originally due to Rudolph \cite{Rudolph90}.
Both inequalities can be proved using state models for the polynomials, see \cite{Tabachnikov97}.
Both inequalities have been extended to Legendrian knots in $J^1S^1$ the one jet space
of $S^1,$ see \cite{ChmutovGoryunov97}. For a more detailed discussion of the history of these inequalities
see \cite{Ferrand02}.

\bex
For the left handed trefoil we have 
\[P_K(v,z)=v^{-3}( \frac{v-v^{-1}}{z})((2v-v^{-1}-vz^2)\]
and 
\[F_K(a,x)=a^{-3}( 1+\frac{a-a^{-1}}{x})(2a-a^{-1}+z-a^{-2}z+az^2-a^{-1}z^2).\]
Thus we get the following bounds on Legendrian left handed trefoils $L$
\[tb(L)+|r(L)|\leq d_{P_K}=-5\]
and
\[tb(L)\leq -6.\]
We will see below that the maximal $tb$ for left handed trefoils is $-6$ thus the second
inequality seems sharper, but the trefoils with $tb=-6$ must have $r=\pm1.$ Thus both 
inequalities are sharp in this case. Note that in this case the Bennequin inequality is
not sharp, as it only give an upper bound of 1.
\eex

While no known bound on $tb$ is always sharp the Kauffman polynomial $F_K$ frequently
seems to be very good. For example 
\begin{thm}[Ng 2001, \cite{Ng01}]
If \K is a 2-bridge knot, then the Kauffman bound is sharp: $\overline{tb}(\K)=d_{F_{\K}}.$
Moreover, the Kauffman bound is sharp for all knot types with eight or fewer crossings
except the $(4,-3)$ torus knot ($8_{19}$ in \cite{Rolfsen}).
\end{thm}

\section{New invariants}\label{newinvt}
In this section we discuss two invariants of Legendrian knots that have been 
defined in recent years. They are considerably more complicated than the ``classical
invariants'' of Legendrian knots and currently we do not know exactly what they
tell us about Legendrian knots. But they do show us that the classical invariants
do not completely determine the Legendrian isotopy type of Legendrian knots.

\subsection{Contact homology (aka Chekanov-Eliashberg DGA)}\label{chintro}
In \cite{Eliashberg00} Eliashberg first outlined the general theory of contact homology. It is
a sophisticated algebraic tool for keeping track of ``Reeb cords'' and ``pseudo-holomorphic
curves'' in symplecizations of contact manifolds. Subsequently the more elaborate framework
of Symplectic Field Theory has been developed in \cite{EGH}. The foundational parts of these
grand theories (especially the parts of interest here) 
are currently being developed. In certain circumstances one can translate the general ideas into
simpler terms and rigorously establish a theory based in these simpler terms. The 
differential graded algebra (DGA) of Chekanov \cite{Chekanov} is such a theory for Legendrian knots
in $\R^3.$   We will concentrate on
this simpler, but still deep, theory here. For a discussion of how the theory described here
fits into the ``big picture'' see \cite{EtnyreNgSabloff}.
In \cite{Sabloff} a similar theory was established for Legendrian knots in circle bundles
(with contact structures transverse to the circle fibers).

\begin{figure}[ht]
  {\epsfxsize=2.8in\centerline{\epsfbox{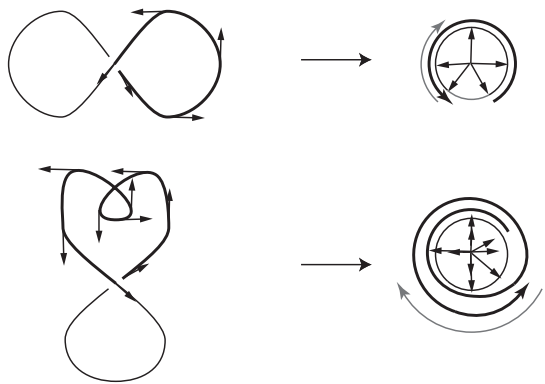}}}
  \caption{On the right are two crossings, the thicker arc is $\gamma_c.$ On the right is 
    the image of the Gauss map for the arc. The grey arrow indicates the completion of $g_c$ to a closed loop.
    When the Gauss map is projectivized we see the top example has grading 1 and the bottom example
    has grading 3.}
  \label{fig:grad}
\end{figure}

\begin{figure}[ht]
  \relabelbox 
  \small{\epsfxsize=.8in\centerline{\epsfbox{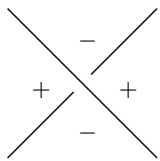}}}
  \relabel {+}{$+$}
  \relabel {-}{$-$}
  \relabel {p}{$+$}
  \relabel {m}{$-$}
  \endrelabelbox
  \caption{Quadrants near a crossing in $\pi(L).$}
  \label{fig:cornerasy}
\end{figure}

\begin{figure}[ht]
  \relabelbox 
  \small{\epsfxsize=4in\centerline{\epsfbox{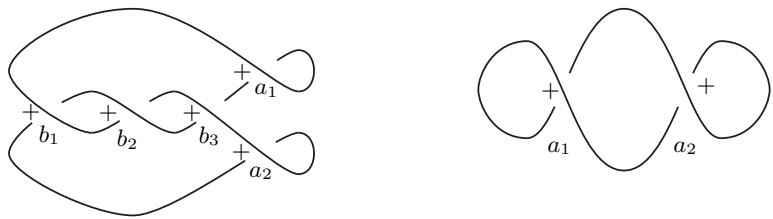}}}
  \relabel {b1}{$b_1$}
  \relabel {b2}{$b_2$}
  \relabel {b3}{$b_3$}
  \relabel {a1}{$a_1$}
  \relabel {c}{$a_2$}
  \relabel {aa1}{$a_1$}
  \relabel {d}{$a_2$}
  \relabel {1}{$+$}
  \relabel {2}{$+$}
  \relabel {3}{$+$}
  \relabel {4}{$+$}
  \relabel {5}{$+$}
  \relabel {6}{$+$}
  \relabel {7}{$+$}
  \endrelabelbox        
  \caption{The Lagrangian projection of a Legendrian trefoil and unknot. 
    One corner of each crossings has
    been labeled with a $+$ to help determine the signs of all the quadrants near the crossing.}
  \label{fig:chex1}
\end{figure}

Let $L$ be such a Legendrian knot. To $L$ we want to assign a graded algebra and a differential
on the algebra.  There are various levels of
sophistication one can use in defining the algebra and its differential. We begin with
the simplest description and discuss generalizations in Section~\ref{exten}.\hfill\break
{\bf The Algebra.} Denote the double points in $\pi(L)$ (recall
this is the Lagrangian projection of $L$) by $\mathcal{C}.$ We assume $\CC$ is a finite set 
of transverse double points (for
a generic Legendrian this is clearly true).
Let $\mathcal{A}_{\Z_2}$ be the free associative unital algebra over $\Z_2$ generated
by $\mathcal{C}.$ If we want to make the generators explicit we write $\A_{\Z_2}(\CC).$\hfill\break
{\bf The Grading.} To each crossing $c\in\CC$ there are two points $c^+$ and $c^-$ in $L\subset\R^3$ that
project to $c.$ We let $c^+$ be the point with larger $z$ coordinate.
Choose a map $\gamma_c:[0,1]\to L$ that parametrizes an arc running from $c^+$ to $c^-$ 
(note there are two such arcs, we can choose either). Let $P$ be the projectivized unit circle in 
the $xy$-plane ({\em i.e.} identify antipodal points). We get a map $g_c:[0,1]\to P$ by $g_c(t)=\frac
{\gamma_c'(t)}{|\gamma_c'(t)|}.$ Since $c$ is a transverse double point $g_c(0)\not=g_c(1).$ We 
extend $g_c$ to a map $g_c:S^1\to P$ by rotating $g_c(1)$ clockwise until it agrees with
$g_c(0).$ Now define the grading on $c$ to be
\[|c|=\text{ degree } (g_c).\]
See Figure~\ref{fig:grad}.
Note the grading is not well-defined! If we chose the other arc from $c^+$ to $c^-$ 
in $L$ then we could have gotten a different number. Note the union of the two arcs is all of $L.$ From this
it is easy to see that the difference between the two possible gradings on $c$ is the twice the rotation
number of $L$ (recall the rotation number is the degree of the Gauss map for the immersion $\pi(L)$). Thus
$|c|$ is well defined modulo $2r(L).$ Thus $\A$ is a graded algebra with grading in $\Z_{2r(L)}.$\hfill\break
{\bf The Differential.} We will define the differential $\partial$ on $\A_{\Z_2}$ by defining it on
the generators of $\A$ and then extending by the signed Leibniz rule:
\[
\partial ab= (\partial a)b+(-1)^{|a|}a\partial b.
\]
(Note that since we are working over $\Z_2$ we do not really need the $(-1)^{|a|}$ in this formula, but
when working over other rings/fields the sign is important so we keep it in the formula.)
The neighborhood of a crossing $c$ in $\R^2$ is divided into four quadrants by $\pi(L)$ two are labeled
$+$ quadrants and two are labeled $-$ quadrants as shown in Figure~\ref{fig:cornerasy}. 
Let $a\in \CC$ be a generator of $\A$ and let $b_1\ldots b_k$ be a word in the ``letters'' \CC.
Let $P_{k+1}$ be a $k+1$ sided polygon with vertices labeled counterclockwise $v_0,\ldots, v_k$ 
and set 
\[\M^a_{b_1\ldots b_k}=
\{u:(P_{k+1},\partial P_{k+1})\to (\R^2, \pi(L)) | \text{ $u$ satisfies 1.-- 3. below}\}/ \text{
reparam.}\]
where $\R^2$ denotes the $xy$-plane and the conditions are
\begin{itemize}
\item[1.] $u$ is an immersion 
\item[2.] $u(v_0)=a$ and near $a$ the image of a small neighborhood of $v_0$ under $u$ covers a $+$
quadrant.
\item[3.] $u(v_i)=b_i, i=1,\ldots, k,$ and near $b_i$ the image of a 
small neighborhood of $v_i$ under $u$ covers a $-$ quadrant.
\end{itemize}
We can now define
\[\partial a=\sum_{b_1\ldots b_k} (\#_2\M) b_1b_2\ldots b_k,\]
where the sum is taken over all words in the letters \CC and $\#_2$ denotes the modulo two count of
elements in $\M.$ It is not hard to show that $\M$ is finite and the sum is finite \cite{Chekanov}, thus
the differential is well-defined.
\begin{lem}[Chekanov 2002, \cite{Chekanov}]
The map $\partial$ is a differential:
\[\partial\circ \partial =0,\]
and lowers degree by one.
\end{lem}
\bex\label{sunknot}
Here we compute the graded algebra and the differential for the Legendrian unknot with $tb=-2$ (and $r=\pm1$).
The Lagrangian projection is shown on the right hand side of Figure~\ref{fig:chex1}. 
We have two generators 
of the algebra $a_1, a_2.$ Each has grading 1. Moreover, $\partial a_i=1$ for $i=1,2.$
\eex
\bex\label{trefoilex}
Consider the Legendrian right handed trefoil with $tb=1.$ The Lagrangian projection is shown on the left of
Figure~\ref{fig:chex1}. The algebra has five generators $a_1, a_2, b_1, b_2, b_3.$ Their gradings are
\[|a_i|=1,\quad |b_i|=0.\]
One easily computes
\begin{align*}
\partial a_1&=1+b_1+b_3+ b_1b_2b_3\\
\partial a_2&=1+b_1+b_3+ b_3b_2b_1\\
\partial b_i&= 0.
\end{align*}
\eex

To discuss the invariance of $(\A_{\Z_2},\partial)$ under Legendrian isotopy we need a few preliminary 
definitions. An automorphism of $\A_{\Z_2}(c_1,\ldots, c_n)$ of the form 
\begin{equation} 
    \phi(c_i) = 
    \begin{cases}
      c_i, & i \neq j \\
      \pm c_j + u, & u \in
      \A(c_1,\ldots,c_{j-1},c_{j+1},\ldots,c_n),\ i=j.
    \end{cases}
  \end{equation}
for some fixed $j$
is called \dfn{elementary}. A composition of such automorphism is called a \dfn{tame automorphism}. 
A tame isomorphism
is an identification of the generators of two algebras followed by a tame automorphism. (Since our algebras
are all graded our automorphisms are also assumed to be graded automorphisms).

An \dfn{index $i$ stabilization}
of a DGA $(\A_{\Z_2}(c_1, \ldots, c_n),\partial)$ is the graded algebra $\A_{\Z_2}(c_1, \ldots, c_n, a, b)$
where $|a|=|b|-1=i$ and the differential on this algebra agrees with the original differential on
the $c_i$'s and $\partial a=0,$ and $\partial b=a.$  If the index of a stabilization is unimportant 
then it is simply referred to as a stabilization instead of an index $i$ stabilization.

Two DGA's are \dfn{stably tame isomorphic} if after stabilizing each of the DGA's some number of times
they become tame isomorphic.

\begin{thm}[Chekanov 2002, \cite{Chekanov}]
The differential graded algebra $(\A_{\Z_2},\partial)$ associated to a Legendrian knot changes 
by stable tame isomorphisms under Legendrian isotopy, moreover the homology 
\[CH_*(L)=\frac{\ker \partial}{\hbox{\em image } \partial}\]
is unchanged under Legendrian isotopy.
\end{thm}

\bex\label{compstab}
The contact homology for the Legendrian unknot in Example~\ref{sunknot} is
\[CH_*(L)=0,\]
since $\partial a_1=1.$ Indeed $\ker(\partial)$ is spanned by $a_1-a_2$ but 
$\partial (a_1(a_1-a_2))=a_1-a_2.$
\eex
\bbr
In general it is easy to show that anytime $\partial a=1$ for some element in the algebra then $CH_*=0.$
\eer
\bex
The contact homology for the Legendrian trefoil in Example~\ref{trefoilex} is
\[CH_*(L)= CH_0(L)=\langle b_1, b_2, b_3 | b_1b_2b_3= b_3b_2b_1=0\rangle\]
\eex

\begin{prop}[Chekanov 2002, \cite{Chekanov}]
The DGA associated to a stabilized knot is stably tame isomorphic to the ``trivial'' algebra, that is the
algebra generated by one element $a$ such that $\partial a=1.$ In particular, the contact homology is 0.
\end{prop}
The main observation in the proof of this proposition is that a stabilization in the Lagrangian projection is
achieved by adding a small loop (see Figure~\ref{fig:stab}). Since this loop is small a generator $a$ has been added to 
the algebra for which $\partial a=1$ by the following lemma.
\begin{lem}[Chekanov 2002, \cite{Chekanov}]\label{lem:area}
If $u:P_{k+1}\to \R^2$ is a polygon in $\M^a_{b_1\ldots b_k}$ then
\[h(a)-\sum_{i=1}^k h(b_i)=\int_{P_{k+1}} u^*(dx\wedge dy)\]
where $h(c)$ is the difference in the $z$-coordinates of the two points $c^+, c^-$ in $L$ lying above a 
double point $c.$ 
\end{lem}
Thus we easily see the contact homology is 0 as in
Example~\ref{compstab}. With a little more work one can see the algebra is as claimed in the proposition.

\subsection{Linearization}

Since our algebra is non-commutative it is quite hard in general to look at the 
homology of two DGA's and determine if they are the
same or not. In this section we describe two methods to extract more ``computable'' invariants
out of the DGA associated to a Legendrian knot.

Given a DGA $(\A,\partial)$ we can filter $\A$ by ``word length''. That is, define
$\A_n$ to be the vector space generated by words in the generators of $\A$ of length
less than or equal to $n.$ Note $\A_n$ is not an algebra, but merely a vector space.
We call the DGA \dfn{augmented} if for each generator $a\in\A,$ 
$\partial a$ contains no constant term. Said another way the ``the differential is non-decreasing on word length''.
It is easy to see that if $(\A,\partial)$ is an augmented DGA then $\partial$ induces a differential
$\partial_n$ on $\A_n$ by applying $\partial$ to a generator of $\A_n$ and projecting the result to
$\A_n.$  If $(\A,\partial)$ is the DGA associated to a Legendrian knot $L$ and
it is augmented then we call the homology of $(\A_n,\partial_n)$ the \dfn{order $n$ contact homology}
of $(\A,\partial)$ and is denoted $L_nCH_*(L).$ The first order contact homology is frequently referred
to as the \dfn{linear (or linearized) contact homology}. The Poincar\'e polynomial 
\[P_n(\lambda)=\sum_{i\in\Z_{2r(L)}} \dim(L_nCH_i(L)) \lambda^i\]
of $L_nCH_*(L)$ is called the \dfn{order $n$ Chekanov-Poincar\'e polynomial} of $(\A,\partial).$
It is interesting to observe (see \cite{Chekanov}) 
that one can compute the Thurston--Bennequin invariant of a Legendrian knot
by evaluating $P_1(\lambda)$ at $-1:$
\[tb(L)=P_1(-1).\]
\hfill\break
{\bf Important Facts:} First, not all DGA's are augmented, and second, the 
order $n$ contact homology  and order $n$ Chekanov-Poincar\'e polynomial
of $(\A,\partial)$ is not an invariant of the stable tame isomorphism class of $(\A,\partial)$!

These two facts seem to indicate these notions are useless for Legendrian knots and that our notation
is not so good. It turns out we can still make use of these order $n$ approximations to the contact
homology. If $\A$ is generated by $a_1,\ldots, a_k$ then set 
\[G(\A)=\{g \hbox{ a graded automorphism of $\A$ of the form $a_i\mapsto a_i+c_i$ where $c_i\in\Z_2$}.\}\]
Every element $g\in G(\A)$ gives a tame isomorphism of DGA's from $(\A,\partial)$ to $(\A,\partial^g),$
where $\partial^g=g\partial g^{-1}.$ Now set 
\[G_a(\A)=\{g\in G(\A) : (\A,\partial^g) \hbox{ is augmented}\}.\]
\begin{thm}[Chekanov 2002, \cite{Chekanov}]
The sets 
\[I_n(L)=\{\hbox{Chekanov-Poincar\'e polynomial of order $n$ for $(A,\partial^g)$} :
g\in G_a(\A)\}\]
and
\[L_nCH_*(L)=\{L_nCH_*(\A,\partial^g): g\in G_a(\A)\}\]
are invariants of the Legendrian knot $L.$
\end{thm}
It is relatively easy to show these sets do not change under tame isomorphisms. Thus to prove this theorem
one only needs to check what happens to $G_a, I_n$ and $L_nCH_*$ under stabilization.

A convenient way to find $G_a(\A)$ is via augmentations. Let $(\A,\partial)$ be a DGA. A map 
\[\epsilon: \A\to \Z_2\]
is called an \dfn{augmentation} of $\A$ if $\epsilon(1)=1,$ $\epsilon\circ \partial=0$ and
$\epsilon$ vanishes on any element of non-zero degree.
\begin{lem}
The augmentations of $(A,\partial)$ are in one to one correspondence with $G_a(\A)$.
\end{lem}
This lemma is simple to prove once one observes that given an augmentation $\epsilon$ the 
tame automorphism $g(a_i)= a_i+\epsilon (a_i)$ is in $G_a(\A).$ 

\begin{figure}[ht]
  \relabelbox 
  \small{\epsfxsize=3.7in\centerline{\epsfbox{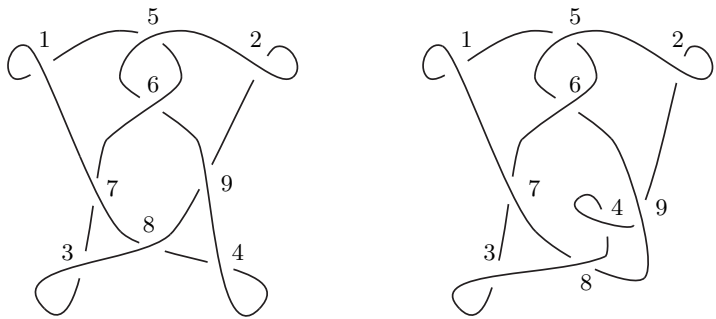}}}
  \relabel {1}{$1$}
  \relabel {2}{$2$}
  \relabel {3}{$3$}
  \relabel {4}{$4$}
  \relabel {5}{$5$}
  \relabel {6}{$6$}
  \relabel {7}{$7$}
  \relabel {8}{$8$}
  \relabel {9}{$9$}
  \relabel {a}{$1$}
  \relabel {b}{$2$}
  \relabel {12}{$3$}
  \relabel {13}{$4$}
  \relabel {14}{$5$}
  \relabel {15}{$6$}
  \relabel {16}{$7$}
  \relabel {17}{$8$}
  \relabel {18}{$9$}
  \endrelabelbox      
  \caption{The knot $L$ (left) and $L'$ (right).}
  \label{fig:chex}
\end{figure}

\bex
In this example we consider the two ``Chekanov-Eliashberg knots''. These were the first two Legendrian
knots with the same $tb, r$ and knot type that were shown to be non Legendrian isotopic. The 
Lagrangian projection of these knots is shown if Figure~\ref{fig:chex}. 
The generators for $L$ are $a_i, i=1,\ldots 9.$ Their gradings are
\begin{align*}
|a_i|&= 1, \quad i=1\ldots 4,\\
|a_5|&= 2,\\
|a_6|&= -2,\\
|a_i|&= 0, \quad i=7\ldots 9
\end{align*}
and the boundary map is
\begin{align*}
\partial a_1&= 1+ a_7 + a_7a_6a_5,\\
\partial a_2&= 1+ a_9+ a_5a_6a_9,\\
\partial a_3&= 1+ a_8a_7,\\
\partial a_4&= 1+ a_9a_8,\\
\partial a_i&=0\quad i\geq 5.
\end{align*}

The generators for $L'$ are $b_i, i=1,\ldots 9.$ Their gradings are 
\begin{align*}
|b_i|&= 1, \quad i=1\ldots 4,\\
|b_i|&= 0, \quad i=5\ldots 9
\end{align*}
and the boundary map is
\begin{align*}
\partial b_1&= 1+ b_7 + b_7b_6b_5 +b_5+ b_9b_8b_5,\\
\partial b_2&= 1+ b_9+ b_5b_6b_9,\\
\partial b_3&= 1+ b_8b_7,\\
\partial b_4&= 1+ b_9b_8,\\
\partial b_i&=0\quad i\geq 5.
\end{align*}

Though these differentials look different it is difficult at this stage to see that the algebras are
not stability tame isomorphic. We now compute the linearized contact homology of $L.$ We begin by looking
for augmentations. These must be maps $\epsilon:\A\to \Z_2$ that send $a_i$ to $0$ for $i\leq 6$ and
send $a_i$ to $a_i+c_i$ for $i\geq 7$ where $c_i$ is $0$ or $1.$ The equation $\epsilon\circ \partial=0$ implies
 \begin{align*}
1+c_7+c_7c_6c_5&=0,\\
1+c_9+c_5c_6c_9&=0,\\
1+c_8c_7&=0,\\
1+c_9c_8&=0.
\end{align*}
The only solution to these equations is $c_7=c_8=c_9=1.$ The differential $\partial_\epsilon$ associated to this
augmentation is
\begin{align*}
\partial_\epsilon a_1&= a_7 + a_6a_5+ a_7a_6a_5,\\
\partial_\epsilon a_2&= a_9+ a_5a_6+ a_5a_6a_9,\\
\partial_\epsilon a_3&= a_8+a_7+a_8a_7,\\
\partial_\epsilon a_4&= a_9+a_8+ a_9a_8,\\
\partial_\epsilon a_i&=0\quad i\geq 5.
\end{align*}
Thus one easily computes the linearized homology is generated by $a_1+a_2+a_3+a_4, a_5, a_6.$ It has
dimension one in grading $-2, 1,2$ so $P_1(\lambda)=\lambda^{-2}+\lambda+\lambda^2.$

For $\A(L')$ 
one finds three augmentations $c_i=0, i=1\ldots 4,$ $c_7=c_8=c_9=1$ and
$c_5c_6=0.$ No matter which augmentation is used one finds the linearized homology generated by $b_1+b_2+b_3+b_4,
b_5, b_6.$ Thus it has dimension one in grading 1 and dimension two in grading 0.

Hence one can distinguish $L$ and $L'$ using the linearized contact homology.
\eex

\subsection{The characteristic algebra}
Let $(\A,\partial)$ be the DGA generated by $a_1, \ldots, a_k$. Let $I$ be the two sided ideal
in \A generated by $\langle \partial a_i \rangle_{i=1}^k.$ The \dfn{characteristic algebra}
of $(\A,\partial)$ is 
\[\CC(\A)=\A/I.\]
Two characteristic algebras $\A_1/I_1$ and $\A_2/I_2$ are \dfn{tame isomorphic} if after adding
some generators to $A_i$ and the same generators to $I_i, i=1,2,$ then there is a tame 
isomorphism of $\A_1$ to $\A_2$ taking $I_1$ to $I_2.$ It is important to note that the definition
of tame isomorphism for characteristic algebras actually requires the pair $\A,I$ not just $\A/I.$
So one should probably define the characteristic algebra to be the pair $(\A,I)$ but we will
hold with tradition and define $\CC(\A)=\A/I,$ keeping in mind the need to remember both
$\A$ and $I.$ We say $\A_1/I_1$ and $\A_2/I_2$ are \dfn{equivalent} if after adding some 
(possibly different) number of generators to both $A_1$ and $A_2$ (but not the ideals $I_i$) 
the algebras $\A_1/I_1$ and $\A_2/I_2$ are tame isomorphic.
\begin{thm}[Ng 2001, \cite{Ng01}]
Legendrian isotopic  knots have equivalent characteristic algebras.
\end{thm}

\begin{figure}[ht]
  \relabelbox 
  \small{\epsfxsize=3.5in\centerline{\epsfbox{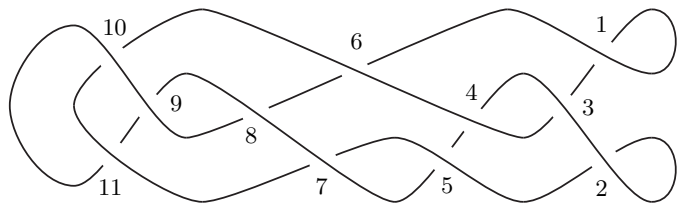}}}
  \relabel {1}{$1$}
  \relabel {2}{$2$}
  \relabel {3}{$3$}
  \relabel {4}{$4$}
  \relabel {5}{$5$}
  \relabel {6}{$6$}
  \relabel {7}{$7$}
  \relabel {8}{$8$}
  \relabel {9}{$9$}
  \relabel {10}{$10$}
  \relabel {a}{$11$}
  \endrelabelbox
  \caption{The Legendrian knots $L$.}
  \label{fig:mex}
\end{figure}

We will illustrate the efficacy of the characteristic algebra by addressing Ng's solution to
the ``Legendrian mirror problem''
\cite{FuchsTabachnikov97}. 
That is: given a Legendrian knot $L$ is it isotopic to its image
under the map $(x,y,z)\mapsto (x,-y,-z)$? The answer is sometimes NO.
\bex
Consider the Legendrian knot $L$ whose Lagrangian projection is shown if Figure~\ref{fig:mex}. 
The algebra has
11 generators with gradings $|a_1|,=|a_2|=|a_7|= |a_9|=|a_{10}|=1,$ $|a_3|=|a_4|=0$ and 
$|a_5|=|a_6|=|a_8|=|a_{11}|=-1.$ The differential is
\begin{align*}
\partial a_1 = 1 + a_{10}a_5a_3 \quad \quad & \partial a_6= a_{11}a_8\\
\partial a_2 = 1 + a_3(1+a_6a_{10}+a_{11}a_7) \quad \quad& \partial a_7= a_8a_{10}\\
\partial a_4 = a_{11} +(1+a_6a_{11} + a_{11}a_7)a_3 \quad \quad& \partial a_9=1+a_{10}a_{11}\\
\partial a_5=\partial a_8=\partial a_{11}=\partial a_3=0.
\end{align*}
Let $I$ be the ideal generated by these relations and set $\CC(\A)=\A/I.$ The characteristic algebra $\CC(\A'),$ 
the algebra associated to the Legendrian mirror of $L,$ is the same except all the terms in $I$ are reversed.

One may show \cite{Ng01} that 
\[\CC(\A)= \langle a_1,\ldots, a_10\rangle/\langle 1+ a_{10}a_5a_3, 1+a_3a_{10}a_5, 1+
a_{10}^2a_5^2, 1+a_{10}a_5+a_6a_{10}+ a_{10}a_5^2a_7\rangle\]

We note that $\CC(\A)$ does not have elements $v$ and $w$ of grading $-1$ and $1,$ respectively, 
that satisfy $vw=1.$
The reason for this is if there were such elements then we could consider the further quotient 
algebra $\CC'$ obtained
by setting $a_3=1, a_1=a_2=a_6=a_7=a_9=0.$ This algebra has presentation $\langle a_5, a_{10}\rangle/\langle 
1+a_{10}a_5\rangle.$ Moreover the elements $v$ and $w$ would map to such elements in $\CC'.$ 
But in $\CC'$ one may
easily see that such elements do not exist.

Now in $\CC(\A')$ we see the elements $v=a_{10}$ and 
$w=a_5a_3$ satisfy $vw=1.$ Thus $L$ is not Legendrian isotopic to its mirror.
\eex

Let $\gamma_i$ be the number of generators of degree $i$ of the DGA $(\A,\partial).$ The \dfn{degree
distribution} of $\A$ is the map 
\[\gamma:\Z_{2r(L)}\to \Z_{\geq0}: i\mapsto \gamma_i\]
where $r(L)$ is the rotation number of $L.$

\begin{thm}[Ng 2001, \cite{Ng01}]
The first and second order Poincar\'e-Chekanov polynomials for all possible augmentations of
the DGA $(\A,\partial)$ are determined by the characteristic algebra $\CC(\A)$ and  the degree
distribution of $\A$
\end{thm}

Ng has found examples that the characteristic algebra distinguishes that have the same linearized contact
homology. Thus the characteristic algebra is a strictly stronger invariant that the linearized contact homology.

\subsection{Lifting the DGA to $\Z[t,t^{-1}]$}\label{exten}

\begin{figure}[ht]
  \relabelbox 
  \small{\epsfxsize=1.9in\centerline{\epsfbox{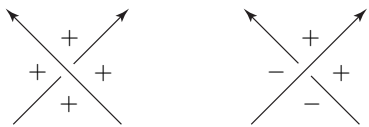}}}
  \relabel {1}{$+$}
  \relabel {2}{$+$}
  \relabel {3}{$+$}
  \relabel {4}{$+$}
  \relabel {5}{$+$}
  \relabel {6}{$+$}
  \relabel {7}{$-$}
  \relabel {8}{$-$}
  \endrelabelbox
  \caption{Signs $s(v_i)$ associated to a corner.}
  \label{fig:os}
\end{figure}

It is fairly simple to describe the lift of the DGA to $\Z[t,t^{-1}].$ Once again
we start with the Lagrangian projection $\pi(L)$ of the Legendrian knot $L$ and assume
the set of double point \CC is finite and consists of transverse double points. Now let
$\A$ be the free associative unital algebra over $\Z[t,t^{-1}]$ generated
by $\mathcal{C}.$ We grade the generators as above and set $|t|=2r(L).$ The grading in this
case is over $\Z$ (and not just $\Z_{2r(L)}$).
Recall to define the grading on $c\in\CC$ we chose paths 
$\gamma_c$ from the upper point $c^+\in L$ above $c$ to the lower point $c^-\in L.$ 
To define the differential we want to assign to each $u\in \M^a_{b_1\ldots b_k}$ (for the 
notation see Section~\ref{chintro}) an integer $n_u$ and a sign $s_u.$ We begin with the integer.
Given $u\in \M^a_{b_1\ldots b_k}$ note that the image of the boundary $u(\partial P_{k+1})$
is a union of arcs in $\pi(L)$ that can all be uniquely lifted to arcs $\alpha_i$ in $L.$ 
Let $\Gamma_u=\cup \alpha_i\cup \gamma_a\cup \gamma_{b_i}.$ This is a closed loop in $L.$ We orient
$\Gamma_u$ so that all the arcs $\alpha_i$ are oriented as the boundary of the polygon $u(P_{k+1}).$
Now $\Gamma_u$ is an oriented loop in the oriented $L=S^1.$ As such it corresponds to an element
in $H_1(L,\Z)=\Z.$ This integer is denoted $n_u.$ 
Now for the sign. To each corner $v_i$ of $P_{k+1}$ we assign a sign $s(v_i)$ as shown in 
Figure~\ref{fig:os}.
The sign is $u$ is defined to be 
\[s_u=\prod_{i=0}^k s(v_i).\]
Now the boundary map is defined on a generator $a\in\CC$ by 
\[\partial a=\sum_{b_1\ldots b_k} \sum_{u\in\M^a_{b_1\ldots b_k}} s_u b_1b_2\ldots b_k t^{n_u},\]
where the first sum is taken over all words in the letters \CC.

Once again we have 
\begin{thm}[Etnyre, Ng and Sabloff 2002, \cite{EtnyreNgSabloff}]
The map $\partial$ is a differential
\[\partial\circ \partial =0.\]
The stable tame isomorphism class of $(A,\partial)$ is invariant under Legendrian isotopy.
\end{thm}

\begin{figure}[ht]
  \relabelbox 
  \small{\epsfxsize=2in\centerline{\epsfbox{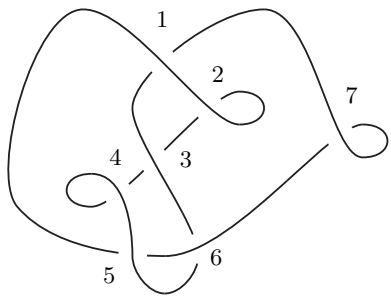}}}
  \relabel {1}{$1$}
  \relabel {2}{$2$}
  \relabel {3}{$3$}
  \relabel {4}{$4$}
  \relabel {5}{$5$}
  \relabel {6}{$6$}
  \relabel {7}{$7$}
  \endrelabelbox
  \caption{The Legendrian knots $L$.}
  \label{fig:8com}
\end{figure}

\bex
Here we compute the contact homology of the Legendrian figure eight knot with $tb=-3.$ 
From Figure~\ref{fig:8com}
we see the algebra has seven generators with $|a_2|=|a_4|=|a_5|=|a_7|=1,$ $|a_1|=|a_3|=0$ 
and $|a_6|=-1.$ We compute
the differential
\begin{align*}
\partial a_1 &= a_6-a_6a_3 -ta_6a_3a_5a_6,\\
\partial a_2 &= t^{-1} +a_1a_3 -a_6a_3a_4,\\
\partial a_4 &= 1-a_3+ta_5a_6a_3,\\
\partial a_7 &=t^{-1}+a_3-ta_3a_6a_3a_5,\\
\partial a_3&= \partial a_5 = \partial a_6=0.
\end{align*}
\eex

It is not currently known if this refined contact homology is a stronger invariant of Legendrian knots than the
original contact homology. One can still define the order $n$ contact homology, and in particular the linearization 
(though augmentations are much harder to find), and the characteristic algebra of this enhanced contact homology.

One potential benefit of including orientations when we count disks in $\M^a_{b_1\ldots b_k}$ is that we can
abelianize. In particular we can define $\A_\Q$ to be the supercommutative associative algebra with unit
over $\Q[t,t^{-1}]$ generated by $\CC,$ the crossings in a Lagrangian projection of a Legendrian knot. We grade
the elements of $\A_\Q$ as we did before and the differential is defined as above as well.
\begin{thm}[Etnyre, Ng and Sabloff 2002, \cite{EtnyreNgSabloff}]
The stable tame isomorphism class and homology of 
$(\A_\Q, \partial)$ are  invariants of the Legendrian isotopy class of
the Legendrian knot $L.$ 
\end{thm} 
The homology of $(\A_\Q,\partial)$ is called the \dfn{abelianized contact homology} of $L.$ 
It is surprising but not known if one can abelianize the contact homology over $\Z$ or $\Z_2.$ In the proof
of this theorem it is simple to show the stable tame isomorphism class of the
differential algebra is an invariant of $L$ 
but to show its homology is
invariant one writes down a chain map that involves division by integers. Thus if the abelianized contact homology
over $\Z$ or $\Z_2$ is invariant the proof will have to be considerably different from the proof of the above
theorem.

\subsection{DGA's in the front projection}\label{FrontDGA}

\begin{figure}[ht]
        {\epsfxsize=2.7in\centerline{\epsfbox{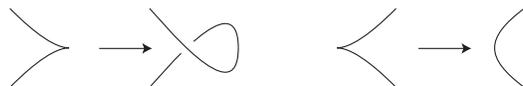}}}
        \caption{Transition from a front projection to a Lagrangian projection.}
        \label{fig:ftol}
\end{figure}

Though it is more natural to define Chekanov's DGA in the Lagrangian projection it
is somewhat difficult to work with Lagrangian projections of Legendrian knots. In this
section we discuss Ng's description of the DGA in the front projection \cite{Ng01}. This 
front description is based on the following observation.
\begin{thm}\label{ftol}
Given the front projection of a Legendrian knot $L$ one can obtain a diagram isotopic
to the the Lagrangian projection by altering the cusps in the
front diagram as shown in Figure~\ref{fig:ftol}.
\end{thm}
\begin{proof}[Sketch of Proof]
The idea is to stretch any Lagrangian projection in the $x$-direction so that it is very long and flat.
Then it is easy to arrange that except near crossings and cusps each strand of the projection has
constant $z$-value. Now, in the complement of the crossings and the cusps, tilt the strands so that
they each have a distinct slope $\epsilon z,$ where $z$ is the $z$ coordinate of the line
before it is tilted. Form the crossings by arcs whose slopes lie between the
slopes of the arcs involved and form the left cusps in a similar fashion. See Figure~\ref{fig:mification}.
\begin{figure}[ht]
        {\epsfxsize=2.7in\centerline{\epsfbox{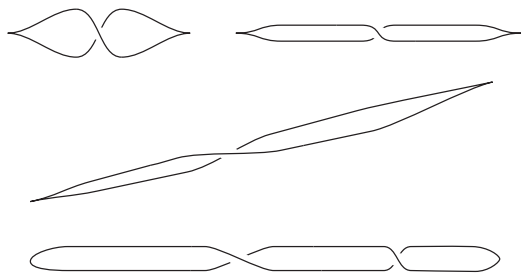}}}
        \caption{Top left: the front projection of a Legendrian unknot. Top right:
        the stretched projection. Middle: the tilted projection that is easy to 
        convert into the Lagrangian projection (bottom).}
        \label{fig:mification}
\end{figure}
It is easy to see the Lagrangian projection of these constant slope lines. With a
moments thought the crossings and cusps project as indicated in the statement of the theorem. 
For full details see \cite{Ng01}.
\end{proof}

\begin{figure}[ht]
        {\epsfxsize=3.3in\centerline{\epsfbox{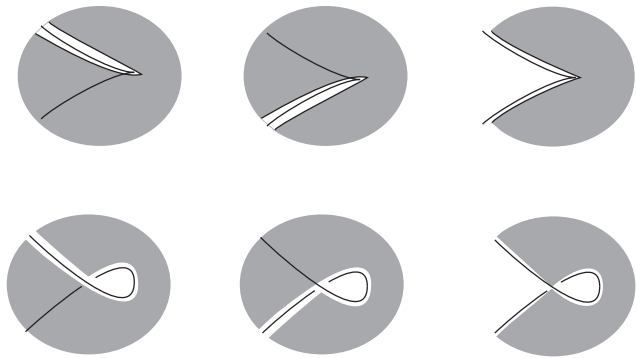}}}
        \caption{The ways a disk can have a negative corner at a right cusp (top) and the
        corresponding disk in the Lagrangian projection (bottom). The gray is the image of
        the disk. The left two diagrams indicate a branch point at the cusp. In this case we count the
        cusp as a corner. The right hand picture indicates a corner at the cusp. In this case we count
        the ``corner'' as a double corner. That is we count it twice when reading off the work that
        goes into the differential.
        }
        \label{fig:cuspch}
\end{figure}

\begin{figure}[ht]
        {\epsfxsize=2.8in\centerline{\epsfbox{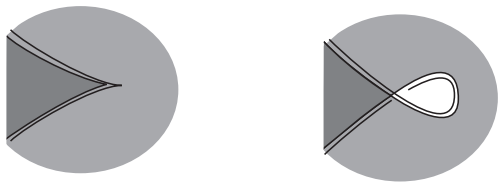}}}
        \caption{A polygon can cross a right cusp as shown here (left). The corresponding
        disk in the Lagrangian projection is also shown (right).
        }
        \label{fig:cusppass}
\end{figure}

With this theorem in mind it is simple to translate the elements used in the DGA into the
front projection.\hfill\break
{\bf The Algebra.} Let \CC be the double points and right cusp points in the projection. Let
$\A$ be the free associative unital algebra over $\Z_2$ generated
by $\mathcal{C}.$ (One can of course translate the algebra over $\Z[t, t^{-1}]$ too, but we leave
this to the reader, or see \cite{EtnyreNgSabloff}.)\hfill\break
{\bf The Grading.} Each right cusp point has grading 1. To each double point $c$ choose a capping
path $\gamma_c$ as we did in Section~\ref{chintro} 
for double points in the Lagrangian projection. The grading 
of the crossing is 
\[|c|=D_c-U_c,\]
where $D_c$ (respectively $U_c$) is the number of cusps traversed downward (respectively upward)
along $\gamma_c.$\hfill\break
{\bf The Differential.} 
Label the quadrants near a double point as in Figure~\ref{fig:cornerasy}. We will think of
a right cusps as a corner (or quadrant) and put a $+$ there.
Left cusps are not thought of as corners.
The differential is exactly as in Section~\ref{chintro} for crossings, except that it may have a corner or
branch point at a right vertex, see Figure~\ref{fig:cuspch}. 
Branch point are counted as corners and a corner at a
right cusp counts as two corners. The reason for this can be seen by transforming the
front projection into the Lagrangian projection as described in Theorem~\ref{ftol}. Moreover
a polygon can pass by a right cusp as shown in Figure~\ref{fig:cusppass}.
The differential of a right cusp is 1 plus the differential defined 
above for a crossing.

One may easily use Theorem~\ref{ftol} to see that this DGA is exactly the one you would get
using the Lagrangian projection of $L.$

\begin{figure}[ht]
  \relabelbox 
  \small{\epsfxsize=2.2in\centerline{\epsfbox{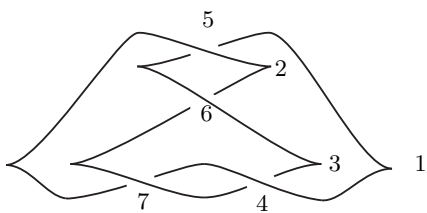}}}
  \relabel {1}{$1$}
  \relabel {2}{$2$}
  \relabel {3}{$3$}
  \relabel {4}{$4$}
  \relabel {5}{$5$}
  \relabel {6}{$6$}
  \relabel {7}{$7$}
  \endrelabelbox
  \caption{The front projection of the Legendrian figure eight knot with $\tb=-3$.}
  \label{fch}
\end{figure}

\bex
Consider the front projection of the Legendrian figure eight knot in 
Figure~\ref{fch}. There are seven generators for the algebra, three from cusps and 4 from crossings.
Their gradings are $|a_i|=1, i=1,\ldots 4,$ $|a_i|=0, i=5,6,7$ and $|a_7|=-1.$ 
The boundary map is 
\begin{align*}
\partial a_1 &= 1+ a_6 a_6a_4a_6a_7+(1+a_6a_5)a_3a_6a_7 +a_6a_2(1+a_6+a_7a_4a_6)a_7,\\
\partial a_2 &= 1+a_5a_6,\\
\partial a_3 &= 1+a_6 + a_6a_7a_4,\\
\partial a_4 &=\partial a_5=\partial a_6 =\partial a_7=0.
\end{align*}
\eex

\bbr
Note that many of the disks in the above computation of difficult to find. This is because
of the right cusps. They allow disks to be quite complicated. To eliminate this problem
one can use the Legendrian Reidemeister moves (Theorem~\ref{thm:reidfront}) to move them all to the far
right of the diagram. This will usually increase the number of crossings in the diagram, but
the computation of the boundary map will usually be considerably easier.
\eer

 \subsection{Decomposition Invariants}\label{sec:di}
Here we define another new invariant of Legendrian knots. This invariant was 
also discovered by Chekanov \cite{Chekanov01, ChekanovPushkar} and can be used to solve Arnold's famous
four cusp conjecture, see Section~\ref{sec:cusp}.

\begin{figure}[ht]
  \relabelbox 
  \small{\epsfxsize=2.7in\centerline{\epsfbox{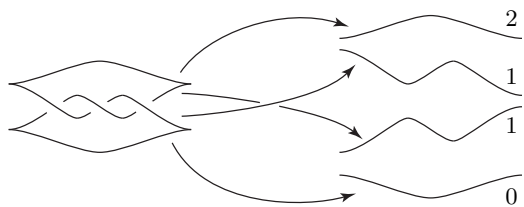}}}
  \relabel {0}{$0$}
  \relabel {1}{$1$}
  \relabel {2}{$2$}
  \relabel {e}{$1$}
  \endrelabelbox
  \caption{The four strands associated (right) associated with the front (left). The numbers 
    are a Maslov potential.}
  \label{decompmp}
\end{figure}

Let $F$ be the front projection of a Legendrian knot $L.$  Let $C(F)$ denote the cusp points in $F$
and let $S(F)$ be the smooth components of $F\setminus C(F).$ These are called the set 
of ``strands'' in $F.$  A \dfn{Maslov potential} for $F$
is a map
\[\mu:S(F)\to \Z_{2r(L)}\]
whose value on the upper stand at a cusp is one larger than its value on the lower
strand. 
Note the Maslov potential is well defined up to adding a constant.
A double point in $F$ is called \dfn{Maslov} if the value of the Maslov potential evaluated on the
two strands crossing there are are the same. See Figure~\ref{decompmp}

\begin{figure}[ht]
        {\epsfxsize=1.9in\centerline{\epsfbox{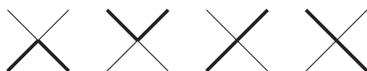}}}
        \caption{The left two figures are switching crossings the right two are non-switching crossings.
        (Different line weights represent different components in the decomposition.)}
        \label{fig:switch}
\end{figure}

Now let $D(F)$ be the double points in $F$ and let $A(F)$ be the closure of the components in 
$F\setminus (D(F)\cup C(F)).$ We call $A(F)$ the arcs in $F.$ A \dfn{decomposition} of $F$ 
is the collection of
simple closed curves $D_F=\{F_1\ldots F_k\}$  such that each $F_i$ is the union of arcs in $F,$ each 
arc is used in only one $F_i,$ and all arcs are in some $F_i.$ 
A crossing in $F$ is called \dfn{switching} or \dfn{non-switching} as indicated in Figure~\ref{fig:switch}.
A decomposition is called \dfn{admissible}
if 
\begin{enumerate}
\item Each $F_i$ bonds a disk $B_i$
\item The ``vertical slice'' $B_{i,x}=B_i\cap \{x=\text{ constant}\}$ of each disk is an interval,
a cusp point or empty
\item Near a switching crossing the intervals $B_{i,x}$ and $B_{j,x}$ are disjoint or one is 
contained within the other, where $F_i$ and $F_j$ are two components of the decomposition coming
together at the crossing.
\item Each switching crossing is Maslov.
\end{enumerate}

\begin{figure}[ht]
        {\epsfxsize=3.4in\centerline{\epsfbox{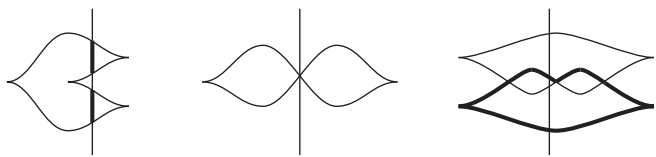}}}
        \caption{Left two figures are not admissible by Item (2). The right figure is not
        admissible by Item (3).}
        \label{fig:badex}
\end{figure}

Note Item (1) rules out trivial decompositions (except for the unknot with $tb=-1$). Item (2) 
rules out the left two decompositions shown in Figure~\ref{fig:badex}. (In fact, a stabilization in 
a diagram will prevent it from having any admissible decompositions because of Item (2).)
Item (3) rules out things like that shown if Figure~\ref{fig:badex}.

Let $\text{Adm}(F)$ be the set of admissible decompositions of $F.$ Define the function 
\[\theta:\text{Adm}(F)\to \Z\]
by $\theta(D_F)=\#(D_F)-S(D_F)$ where $S(D_F)$ is the number of switching crossings that occur in the
decomposition $D_F.$
\begin{thm}[Chekanov and Pushkar, \cite{ChekanovPushkar}] 
Let $L, L'$ be two Legendrian knots in $\R^3$ with generic front projections $F, F'.$  If $L$ and 
$L'$ are Legendrian isotopic then there is a one to one correspondence 
$f:\text{Adm}(F)\to \text{Adm}(F')$
such that $\theta(f(D_F))=\theta(D_F).$ 
\end{thm}
Thus we have two new invariants of Legendrian knots! The cardinality of $\text{Adm}(\Pi(L))$ and the
map $\theta.$ 

\begin{figure}[ht]
  \relabelbox 
  \small{\epsfxsize=4in\centerline{\epsfbox{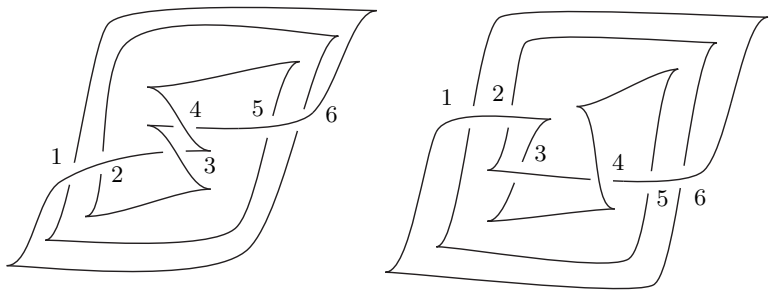}}}
  \relabel {1}{$1$}
  \relabel {2}{$2$}
  \relabel {3}{$3$}
  \relabel {4}{$4$}
  \relabel {5}{$5$}
  \relabel {6}{$6$}
  \relabel {a}{$1$}
  \relabel {b}{$2$}
  \relabel {c}{$3$}
  \relabel {d}{$4$}
  \relabel {e}{$5$}
  \relabel {f}{$6$}
  \endrelabelbox
  \caption{Legendrian projections of the two ``Chekanov knots'' from Figure~\ref{fig:chex}.
    $L$ is on the left and $L'$ is on the right.}
  \label{fig:chd}
\end{figure}

\bex
We can once again distinguish the Chekanov-Eliashberg knots. 
Here we draw the Chekanov examples somewhat
differently, see Figure~\ref{fig:chd}.
In this projection of $L$ one may easily see in any admissible decomposition crossings $2,3,4,5$ must
be switching crossings (by property (1) and (2)). Moreover, crossings 1 and 6 cannot be switching
since they are not Maslov. Thus one sees that this projection of $L$ has $|\text{Adm}(\Pi(L))|=1.$
Similarly in the projection of $L'$ crossings $2,3,4,5$ must be switching and either crossings $1,6$
are both switching or both non-switching. Thus  $|\text{Adm}(\Pi(L'))|=2.$
\eex

It is interesting to note that, just as the contact homology did, 
this invariant also vanishes on stabilized Legendrian knots since
there are no admissible decompositions of a stabilized knot diagram.

The existence of admissible decompositions has an interesting relation to augmentations of the
Chekanov-Eliashberg DGA. In particular Fuchs has shown the following theorem.
\begin{thm}[Fuchs, \cite{Fuchs??}]
If the front projection of a Legendrian knot has an admissible decomposition then
its DGA has an augmentation.
\end{thm}
Very recently Fuchs and Ishkhanov, and independently Sabloff, have found a converse to the above theorem.
\begin{thm}[Fuchs and Ishkhanov \cite{FI}; Sabloff \cite{Sabloff?}] 
If the DGA of a Legendrian knot $L$ has a (graded) augmentation, 
then any front projection of $L$ has a ruling.
\end{thm}
An interesting corollary of this result is the following.
\begin{cor}
If the DGA of a Legendrian knot $L$ has a (graded) augmentation then its rotation number is zero. 
\end{cor}

\section{Classification results}\label{sec:class}
In considering the classification of Legendrian knots it is convenient to consider the map
\[\Psi: {\widetilde{\L}} \to \widetilde{\K}\times\Z\times \Z:(L)\mapsto(k(L),tb(L),r(L)),\]
where $\widetilde{\L}$ is the space of {\em all} oriented Legendrian knots and $\widetilde{\K}$ 
is the space of {\em all} 
oriented topological knots and $k(L)$ is the topological knot type underlying
$L.$ The main questions one would like to answer are: What is the image of $\Psi$? and is
$\Psi$ injective? The answer to the second question is NO as the Chekanov-Eliashberg examples in 
Figure~\ref{fig:chex} show.

We can refine the above questions by considering the map $\Psi$ restricted to $\L(\K)$ for some
knot type $\K.$ This map is denoted $\Psi_{\K}.$ 
A knot type \K is called \dfn{Legendrian simple} if $\Psi_{\K}$ is injective. So for Legendrian simple
knot types a classification of Legendrian knots amounts to the identification of the range of 
$\Psi_{\K}.$

\begin{figure}[ht]
  \relabelbox 
  \small{\epsfxsize=2in\centerline{\epsfbox{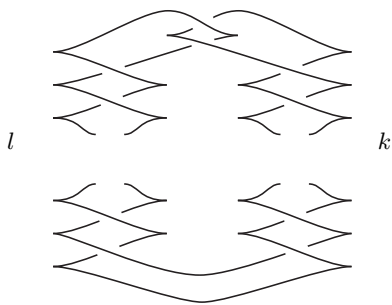}}}
  \relabel {l}{$l$}
  \relabel {k}{$k$}
  \endrelabelbox
  \caption{The Legendrian knot $E(k,l)$. There are $l$ crossings on the right
    hand side and $k$ on the left.}
  \label{fig:gench}
\end{figure}

One can also consider the map
\[\Phi: \widetilde{\L} \to \widetilde{\T}\times \Z,\]
where $\widetilde{\T}$ is the space of all transverse knots and $\Phi(L)=(L_+, tb(L)).$ Recall $L_+$ is
the (positive) transverse push-off of $L.$ Note we did not include the rotation number in the range
of $\Phi$ since $r(L)=tb(L)-sl(L_+).$ This map is potentially a better invariant than $\Psi.$ 
It turns out this map is not injective either.
\bex
Let $E(k,l)$ be the Legendrian knots shown in Figure~\ref{fig:gench}. 
Note the topological knots underlying $E(k,l)$ are determined by $k+l.$ 
In \cite{EFM} it was shown
that $E(k,l)$ is Legendrian isotopic to $E(k',l')$ if and only if the unordered pairs $\{k,l\}$ $\{k',
l'\}$ are the same. This proof involves computations of the linearized contact homology.
However if $l$ is odd then the transverse push-offs $E(k,l)_+$ and $E(k-1, l+1)_+$ 
are transversely isotopic.
\eex

Finally we note that all the new invariants described in Section~\ref{newinvt} vanish on stabilized knots
thus cannot distinguish non Legendrian isotopic knots if the knots are stabilized. Below in 
Section~\ref{sec:consum} we will show there are stabilized Legendrian knots that are not Legendrian 
isotopic. Thus contact homology and the decomposition invariants are not complete invariants.

A standard procedure for trying to classify Legendrian knots in a knot type
\K is to show 
\begin{itemize}
\item All Legendrian knots in $\mathcal{L}(\mathcal{K})$ can be destabilized except
        those with maximal Thurston-Bennequin invariant.
\item Classify Legendrian knots in  $\mathcal{L}(\mathcal{K})$ with maximal Thurston-Bennequin invariant.
\item Understand when Legendrian knots become the same under stabilization.
\end{itemize}
Recall, a Legendrian knot $L$ destabilizes if there is a Legendrian knot $L'$ such that
$L=S_\pm(L').$
This procedure does not always work (see the sections on connected sums and cablings) but
it is frequently a useful strategy.

\subsection{The unknot}
The main theorem concerning unknots is the following.
\begin{thm}[Eliashberg and Fraser 1995, \cite{EliashbergFraser}]
In any tight contact three manifold,
Legendrian unknots are determined by their Thurston-Bennequin invariant and rotation
number. All Legendrian unknots are stabilizations of the unique one
with $tb=-1$ and rotation $r=0,$ see Figure~\ref{fig:three}. 
\end{thm}

The bound on the Thurston Bennequin invariant follows from Bennequin's inequality~\eqref{tbbound}. The
theorem then follows from two lemmas.
\begin{lem}
Any Legendrian unknot with $tb<-1$ destabilizes.
\end{lem}
\begin{proof}
Let $D$ be a disk that $L$ bounds.
If $tb(L)<-1$ then we may make $D$ convex and it will have at least two
dividing curves by Equation~\eqref{tbforconvex}. Since $D$ is a disk one of these dividing
curves must be boundary parallel (recall this means that it separates off a disk
containing no other dividing curves). As discussed in Section~\ref{sfcandinvt} we may use
this dividing curve to destabilize $L.$ 
\end{proof}

\begin{lem}
There is a unique Legendrian unknot with $tb=-1.$
\end{lem}
\begin{proof}
Let $L_i$ be two Legendrian unknots with $tb=-1.$ 
Let $D_i$ be a disk that $L_i$ bounds. We can make $D_i$ convex.
Using the Legendrian Realization Principle \cite{Giroux91, Honda1} we can arrange that the
characteristic foliation on $D_i$ has two elliptic point on $\partial D_i$ and the foliation on the
rest of $D_i$ is by arcs. Thus we may Legendrian
isotop $L_i$ to an arbitrarily small neighborhood of an arc $A_i.$  One may easily show
that any two Legendrian arcs are Legendrian isotopic. Thus we can take $A_1=A_2.$ It is
now an exercise that any two ``thickenings'' of $A=A_i$ to a Legendrian unknot are Legendrian
isotopic. 
\end{proof}

\subsection{Torus knots}\label{sec:torus} 
Let $T$ be the boundary of a neighborhood of an unknot in $M^3.$ A knot that can be
isotoped to sit on $T$ is called a torus knot. 
\begin{thm}[Etnyre and Honda 2001, \cite{EtnyreHonda01b}]\label{maintorus}
        In any tight contact three manifold,
        Legendrian torus knots
        are determined up to Legendrian isotopy by their knot type, Thurston-Bennequin invariant
        and rotation number.
\end{thm}
To complete the classification of Legendrian torus knots we need to identify the invariants 
that are realized. For this we need to more carefully describe torus knots.
Let $V$ be an embedded unknoted solid torus in $M.$ Let $\mu$ be the meridian to $V$ and $\lambda$
be the longitude. Now all torus knots can be isotoped onto $\partial V$ so they can
be expressed by $p\mu+q\lambda$ where $|p|>q>0$ and $(p,q)=1.$ 
If $p>0$ then the torus knot is called \dfn{positive} otherwise it is called \dfn{negative}.
\begin{thm}[Etnyre and Honda 2001, \cite{EtnyreHonda01b}]\label{ptorus}
All Legendrian positive $(p,q)$-torus knots are stabilizations of the unique one
with $tb=pq-p-q$ and rotation $r=0,$ see Figure~\ref{fig:postor}. 
\begin{figure}[ht]
  \relabelbox 
  \small{\epsfxsize=2.3in\centerline{\epsfbox{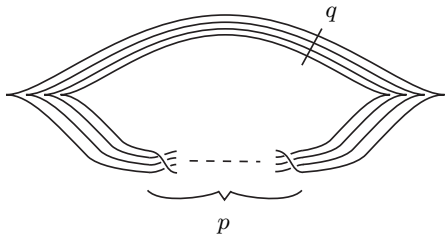}}}
  \relabel {p}{$p$}
  \relabel {q}{$q$}
  \endrelabelbox
  \caption{Positive $(p,q)$-torus knot.}
  \label{fig:postor}
\end{figure}
\end{thm}
\begin{thm}[Etnyre and Honda 2001, \cite{EtnyreHonda01b}]\label{ntorus}
All Legendrian negative $(p,q)$-torus knots are stabilizations of one with maximal $tb$
which is equal to $pq.$ Moreover, $(p,q)$-torus knots with maximal $tb$ are classified
by their rotation number and the set of realized rotation number is
\[\{\pm(p+q+n2q)| 0\leq n\leq \frac{2(p+q)}{q}\},\]
See Figure~\ref{fig:nleg}.
\begin{figure}[ht]
  \relabelbox 
  \small{\epsfxsize=3.9in\centerline{\epsfbox{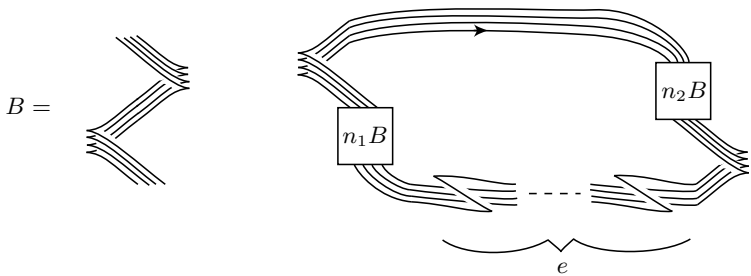}}}
  \relabel {B}{$B =$}
  \relabel {e}{$e$}
  \relabel {1}{$n_1B$}
  \relabel {2}{$n_2B$}
  \endrelabelbox
  \caption{Negative $(-p,q)$-torus knot. Here $p=(n_1+n_2+1)q+e.$ So for each
    such pair of positive integers $n_1,n_2$ there is a maximal $tb,$ $(-p,q)$-torus knot.}
  \label{fig:nleg}
\end{figure}
\end{thm}

\begin{figure}[ht]
  \relabelbox 
  \small{\epsfxsize=4.3in\centerline{\epsfbox{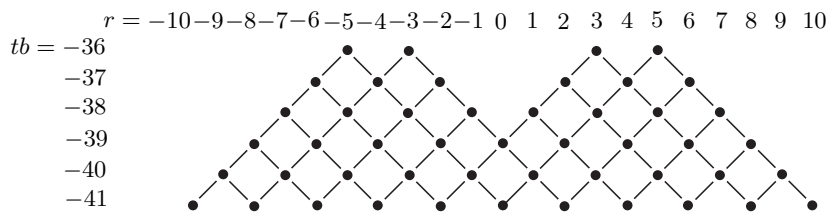}}}
  \relabel {-10}{$-10$}
  \relabel {-9}{$-9$}
  \relabel {-8}{$-8$}
  \relabel {-7}{$-7$}
  \relabel {-6}{$-6$}
  \relabel {-5}{$-5$}
  \relabel {-4}{$-4$}
  \relabel {-3}{$-3$}
  \relabel {-2}{$-2$}
  \relabel {-1}{$-1$}
  \relabel {0}{$0$}
  \relabel {1}{$1$}
  \relabel {2}{$2$}
  \relabel {3}{$3$}
  \relabel {4}{$4$}
  \relabel {5}{$5$}
  \relabel {6}{$6$}
  \relabel {7}{$7$}
  \relabel {8}{$8$}
  \relabel {9}{$9$}
  \relabel {10}{$10$}
  \relabel {tb = -36}{$tb=-36$}
  \relabel {-37}{$-37$}
  \relabel {-38}{$-38$}
  \relabel {-39}{$-39$}
  \relabel {-40}{$-40$}
  \relabel {-41}{$-41$}
  \relabel {r=}{$r=$}
  \endrelabelbox
  \caption{Possible $tb$ and $r$ for the $(-9,4)$-torus knot.}
  \label{fig:geog}
\end{figure}

In Figure~\ref{fig:geog} we see the possible $tb$ and $r$ of Legendrian $(-9,4)$-torus 
knots. The points indicate a unique 
Legendrian knot with corresponding invariants. 
The edges indicate either positive or negative stabilizations (depending
on whether the edge has a positive or negative slope).  This picture will be called the
\dfn{mountain range} associated to the knot type. We have not discussed mountain
ranges up till now because all previous knot types had a particularly simple 
mountain range.

\begin{proof}[Proof of Theorem~\ref{maintorus}]
We restrict our attention to $(S^3,\xi_{std})$ for convenience.  Let $T$ be a standardly 
embedded convex 
torus in $(S^3,\xi_{std})$ with two dividing curves. We will need two
facts \cite{Giroux00, Honda1} concerning this situation:
\begin{enumerate}
\item The slope $s$ of the dividing curves $\Gamma_T$ on $T$ is less than 0.
\item Any other slope in $(-\infty, 0)$ can be realized by a standardly embedded
convex torus parallel (and disjoint from) $T$ with two dividing curves. Moreover,
slopes larger than $s$ are realized by tori on one side of $T$ and slopes smaller 
than $s$ are realized on the other.
\end{enumerate}

We also need the following topological fact: the framing on the normal bundle
of a $(p,q)$-torus knot given by its Seifert surface differs from the one induced 
by the standard torus $T$ on which it sits by $pq.$ 

Now consider a Legendrian positive $(p,q)$-torus knot $L.$ 
The Bennequin inequality says that $tb(L)\leq pq-p-q.$ Thus $tw(L,T)\leq -p-q<0,$ 
so we can make $T$ convex without moving $L.$ Throughout this section we will
assume $T$ has only two dividing curves. It is easy to deal with the more
general case but this is left to the reader. Suppose the dividing curves on
$T$ have slope $-\frac{n}{m},$ then $tw(L,T)=-(pm+qn).$ For all positive, relatively
prime $n,m$ this number is always less than or equal to $-(p+q)$ with equality only
when $m=n=1.$  Thus if $tw(L,T)<-(p+q)$ (that is $tb(L)<pq-p-q$) then there is a 
convex torus $T'$ whose dividing curves have slope $-1$  and cobounds with $T$
a toric annulus $N=T^2\times[0,1].$ Let $A=L\times[0,1]$ in $N,$ here $L=L\times\{0\}.$ 
We can assume that $\partial A$ is Legendrian and $A$ is convex. All the dividing curves 
on $A$ are properly embedded arcs and they intersect $L\times\{0\}\subset \partial A,$
$2(pm+qn)$ times and $L\times\{1\}\subset \partial A,$ $2(p+q)$ times. From this one
may easily conclude that there is a dividing curve on $A$ that has boundary on
$L=L\times\{0\}$ and separates off a disk from $A.$ Thus as discussed in Section~\ref{sfcandinvt} 
we can use this dividing curve to destabilize $L.$ 

Using almost the same argument it is easy to show that a negative $(p,q)$-torus 
knot destabilized to $pq.$  This is left as an exercise, or see \cite{EtnyreHonda01b}. The fact that
this is the maximal $tb$ for such a Legendrian knot is more difficult to establish. In particular,
when $q$ is odd
it does not follow from any of the inequalities in Section~\ref{sec:bounds}. See \cite{Fuchs??}.
\end{proof}

\begin{proof}[Proof of Theorem~\ref{ptorus}]
Let $L_1$ and $L_2$ be two Legendrian positive $(p,q)$-torus knots with maximal $tb.$ 
From the proof of the previous theorem we know that $L_i$ sits on a convex torus
$T_i$ with two dividing curves of slope $-1,$ for $i=1,2.$ Let $V_i$ and$V_i'$ be the
two solid tori into which $T_i$ breaks $S^3.$ Using the
Legendrian realization principle we may isotop $T_i$ relative to $L_i$ so
that the characteristic foliation on $T_1$ is the same as the one on $T_2.$ Let
$\phi$ be a diffeomorphism $T_1$ to $T_2$ that preserves the characteristic foliation and
takes $L_1$ to $L_2.$ 
By the classification of contact structures on solid tori we know $\phi$ extends, as
a contactomorphism, over $V_1$ and $V_1'.$ Thus $\phi$ is a contactomorphism of $S^3$ to
itself that takes $L_1$ to $L_2.$ Applying Theorem~\ref{compclass} we see that $L_1$ and $L_2$ are
Legendrian isotopic.
\end{proof}

The proof of Theorem~\ref{ntorus} is similar to the above proof in spirit but somewhat
more involved. We refer the reader to \cite{EtnyreHonda01b}.

\subsection{Figure eight knot}
The main theorem concerning Legendrian figure eight knots is the following. 
\begin{thm}[Etnyre and Honda 2001, \cite{EtnyreHonda01b}]
In any tight contact three manifold,
Legendrian figure eight knots are determined by their Thurston-Bennequin invariant and rotation
number. All Legendrian figure eight knots are stabilizations of the unique one
with $tb=-3$ and $r=0,$ see Figure~\ref{fig:three}. 
\end{thm}
The bound on $tb$ comes from the inequality in Theorem~\ref{otherbounds}. Thus the following two lemmas
establish the theorem.
\begin{lem}
Any Legendrian figure eight knot with $tb<-3$ destabilizes.
\end{lem}
\begin{lem}
There is a unique Legendrian figure eight knot with $tb=-3.$
\end{lem}
To prove any Legendrian figure
eight knot $L$  with $tb<-3$ destabilizes we consider a minimal genus
Seifert surface $\Sigma$ for $L.$ So $\Sigma$ is a punctured torus and we can
make it convex. We would like to find a bypass  for $L$ on $\Sigma$ so 
we can destabilize $L.$ We might not always be able to do this. Thus to prove $L$ destabilizes
we need to isotop $\Sigma$ until we do see a bypass for $L$ on $\Sigma.$ 
To accomplish this we use
the fibration on the complement of the figure eight knot.
Specifically, let $X=S^3\setminus N$ where $N$ is a standard neighborhood of a 
$L.$ Then $X$ is a $\Sigma$ bundle over $S^1.$ The monodromy for $X$ is 
\[\Psi=\begin{bmatrix} 2 & 1 \\ 1 & 1 \end{bmatrix}.\]
If we cut $X$ open along $\Sigma$ then we get $\Sigma\times [0,1].$ Because of the
monodromy map the dividing curves on $\Sigma\times\{0\}$ and $\Sigma\times \{1\}$ will
be different. We can then take an annulus $A=S^1\times[0,1]$ where $S^1\subset \Sigma,$ 
make $\partial A$ Legendrian and $A$ convex. If we choose $S^1$ correctly we will be
able to find a bypass on $A$ for a boundary component of $A.$ 
We can use this to
alter the dividing curves on $\Sigma,$ \cite{Honda1}. This will frequently produce a bypass on $\Sigma$
for $L.$ The rigorous argument is somewhat involved. The reader is referred to \cite{EtnyreHonda01b} 
for details.  
 
To prove there is a unique Legendrian figure eight knot with $tb=-3$ we use the
above ideas to show that the dividing curves on $\Sigma$ in this situation can
always be arranged to look a certain way and that there is a unique tight contact
structure on $X$ that is a subset of $S^3$ and has the given dividing curves. The uniqueness
proof is then finished as in the proof of Theorem~\ref{ptorus}.

\subsection{Connected sums}\label{sec:consum}
There is a standard way to take a contact connect sum of manifolds. If $(M_i,\xi_i)$ 
are two contact three manifolds, then let $B_i$ be a ball with convex boundary 
contained in a neighborhood
of a point. The connect sum $M_1\# M_2$ is obtained from 
$\overline{M_i\setminus B_i}$ by gluing their boundaries together by a contactomorphism.
This construction preserves tightness \cite{Colin97}. If $L_i$ is a Legendrian knot in $M_i$ we can
arrange that $L_i\cap B_i$ looks like the intersection of the $x$-axis with a
unit ball about the origin in $(\R^3,\xi_{std}).$ We can then define the \dfn{connect
sum} $L_1\#L_2$ in $M_1\# M_2$ to be the closure of $L_1\setminus B_1 \cup L_2\setminus B_2$
(here we of course need to arrange that the gluing map for $M_1\# M_2$ sends $\partial 
(L_1\cap B_1)$ to $\partial (L_2\cap B_2)$).

\begin{figure}[ht]
  \relabelbox 
  \small{\epsfxsize=4in\centerline{\epsfbox{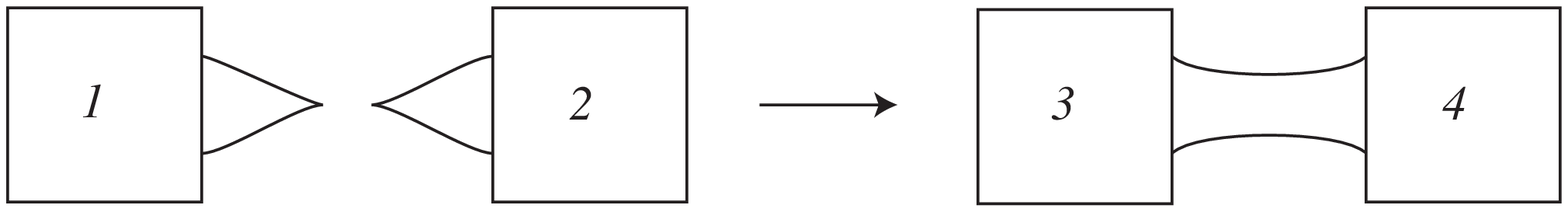}}}
  \relabel {1}{$L_1$}
  \relabel {2}{$L_2$}
  \relabel {3}{$L_1$}
  \relabel {4}{$L_2$}
  \endrelabelbox
  \caption{The connect sum in the front projection.}
  \label{fig:sumd}
\end{figure}

In $(\R^3,\xi_{std})$ the connect sum has a diagrammatic interpretation in the front
projection (note $\R^3\# \R^3$ is not $\R^3$ so we are somewhat abusing notation, but
we may think of one of the $\R^3$'s as $S^3$ with its standard contact structure).
This diagrammatic connect sum is shown in Figure~\ref{fig:sumd}.
It is a good exercise to show that this is the same a the connect sum defined above
and that it is well defined. Moreover using this diagrammatic interpretation it is easy
to see 
\[tb(L_1\# L_2)= tb(L_1) +tb(L_2) +1\]
and
\[r(L_1\# L_2)= r(L_1) +r(L_2).\]

While all the theorems concerning the connect sum work with the more general definition
we restrict attention here to the diagrammatic version in $\R^3.$
The following theorem says that the prime decomposition of Legendrian knots is
unique up to shifting stabilization from one summand to the other and topological
symmetries.
\begin{thm}[Etnyre and Honda 2003, \cite{EtnyreHonda03}]\label{thm:con}
Let $\K=\K_1\#\ldots \#\K_n$ be topological connected sum knot type in $\R^3,$
with $\K_i$ prime. The map
\[\frac{\L(\K_1)\times\ldots \times\L(\K_n)}{\sim} \to \L(\K_1\# \ldots\#\K_n)\]
given by connect sum is a one to one correspondence where $\sim$ is generated by
\begin{enumerate}
\item $(\ldots, S_\pm(L_i),\ldots, L_j, \ldots)\sim(\ldots, L_i,\ldots, S_\pm(L_j), \ldots)$ and
\item $(L_1,\ldots, L_n)\sim(L_{\sigma(1)},\ldots, L_{\sigma(n)})$ where $\sigma$ is permutation
        of $1, \ldots, n$ such that $\K_i=\K_{\sigma(i)}.$ 
\end{enumerate}
\end{thm}
The idea for the proof of this theorem is quite simple. One just uses the
spheres implicated in the definition of connect sum. Make these spheres convex, cut the
manifold along the spheres and glue in standard contact balls with a standard Legendrian
arc in them. One then analyses what happens if one decomposes using a different 
convex sphere. The details can be found in \cite{EtnyreHonda03}.

While the proof of this theorem is fairly straightforward it has some very interesting
consequences.
\begin{cor}
Given any integers $n$ and $m$ there is a topological knot type $\K$
and Legendrian knots $L_1,\ldots, L_n\in\L(\K)$ with the same $tb$ and $r$ that
are not Legendrian isotopic. Moreover, they stay non isotopic even after $m$
stabilizations of any type.
\end{cor}
These are the first examples that cannot be distinguished by the Chekanov-Eliashberg DGA or
the decomposition invariants. In fact, there is no known invariant 
distinguishing many of them.

Recall Theorem~\ref{stablist} says that any two Legendrian knots in the same knot type 
become Legendrian isotopic after some number of positive and negative stabilizations.
There are many theorems in topology concerning some type of equivalence under stabilization
({\em eg.} for Heegaard splittings of three manifolds, smooth structures on four
manifolds, \ldots), but this is the first time one can show that more than one
stabilization is required in the equivalence.

We indicate the proof of this corollary in the case $n=2$ and $m=0.$ To this
end consider the $(-5, 2)$ torus knots $\K.$ In $\L(\K)$ there are four Legendrian
knots with maximal $tb$ (this maximum is $-10$) having rotation numbers
$-3, -1, 1, 3.$ Denote them by $K_r$ where the $r$ indicates their rotation number.
Now set $L_1=K_{-3}\# K_3$ and $L_2=K_{-1}\# K_1.$ Clearly both $L_1$ and $L_2$ have
$tb=-19$ and $r=0.$ However they are not Legendrian isotopic since their summands
are not related by the equivalence relation in Theorem~\ref{thm:con}. Using
the fact that $S_+(K_{-3})=S_-(K_{-1})$ and $S_-(K_3)=S_+(K_1)$ one may easily check
that $S_+(L_1)=S_+(L_2)$ and $S_-(L_1)=S_-(L_2).$ Thus our examples only 
satisfy $m=0$ in the corollary. The reason for this is that the valleys in the 
mountain range for $\K$ (see Section~\ref{sec:torus} for the terminology) are only one
deep. Thus to find examples with $m>0$ we just need to find torus knots having
very deep valleys. This is a simple exercise, or see \cite{EtnyreHonda03}.

\subsection{Cables}
We now discuss Legendrian knots in cabled knot types. To do this we need
some preliminary definitions. We say a solid torus $S$ represents a knot type
\K if its core curve is in $\mathcal{K}.$ Define the width of a knot type as
\beq
w(\K)=\sup \frac{1}{\text{slope}(\Gamma_{\partial S})}, 
\eeq
where the supremum is taken over all solid tori with convex boundary representing \K, and
$\Gamma_{\partial S}$ are the dividing curves on the boundary of $S.$ 
This is clearly an invariant of topological knots.
Since the standard
neighborhood of a Legendrian knot with $tb=n$ has dividing curves on its boundary with
slope $\frac1n$ we can see that
\[\overline{tb}(\K)\leq w(\K) \leq \overline{tb}(\K)+1.\]

A knot type \K is called \dfn{uniformly thick} if 
\begin{enumerate}
\item $w(\K)=\overline{tb}(\K)$ and
\item any solid torus $S$ representing \K is contained in the interior of a solid torus
that is the standard neighborhood for a Legendrian knot in $\L(\K)$ with maximal $tb.$  
\end{enumerate}
For the unknot one may easily check that $w=0$ while $\overline{tb}=-1,$ thus 
the unknot is not uniformly thick. We will see below that there are
many knot types that are uniformly thick, but first we indicate the importance of
being uniformly thick.

Let $\gamma$ be a $(p,q)$ curve on $\partial (S),$ where $S=S^1\times D^2.$ Given a knot
type $K$ let $S'$ be a solid torus representing \K and choose a diffeomorphisms $S\to S'$ 
sending $S^1\times\{pt\}$ to a longitude of $S'$ and $\{pt\}\times \partial D^2$ to a
meridian. Then the knot type $\K_{(p,q)}$ determined by $\gamma$ under this diffeomorphism
is called the \dfn{$(p,q)$-cable of \K}.

\begin{thm}[Etnyre and Honda, \cite{EtnyreHonda??}]
Let \K be a Legendrian simple, uniformly thick knot type. Then
the knot type $\K_{(p,q)}$ is also Legendrian simple.
\end{thm}

Thus we cannot say the cables of Legendrian simple knots are Legendrian simple, but this
is true if the knot type is uniformly thick.

\begin{thm}[Etnyre and Honda, \cite{EtnyreHonda??}]
We have the following:
\begin{enumerate}
\item Negative torus knots are uniformly thick.
\item If \K is uniformly thick and $\frac{p}{q}<w(\K)$ then $\K_{(p,q)}$ is uniformly thick.
\item If $\K_1$ and $\K_2$ are uniformly thick then $\K_1\#\K_2$ is uniformly thick.
\end{enumerate}
\end{thm}

The proofs of these last two theorems are in principle quite similar to the proofs in the 
proceeding sections, especially the section on torus knots. 

\begin{figure}[ht]
  \relabelbox 
  \small{\epsfxsize=3in \centerline{\epsfbox{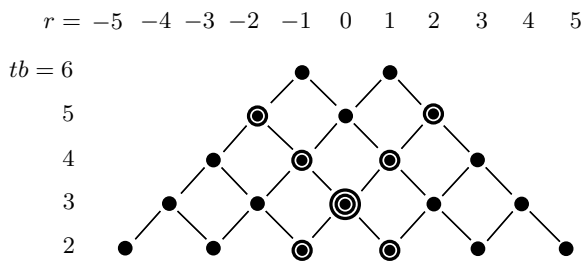}}}
  \relabel {- 5}{$-5$}
  \relabel {- 4}{$-4$}
  \relabel {- 3}{$-3$}
  \relabel {- 2}{$-2$}
  \relabel {- 1}{$-1$}
  \relabel {0}{$0$}
  \relabel {1}{$1$}
  \relabel {2}{$2$}
  \relabel {3}{$3$}
  \relabel {4}{$4$}
  \relabel {5}{$5$}
  \relabel {tb = 6}{$tb=6$}
  \relabel {a}{$5$}
  \relabel {b}{$4$}
  \relabel {c}{$3$}
  \relabel {d}{$2$}
  \relabel {r}{$r=$}
  \endrelabelbox
  \caption{Classification of Legendrian $(2,3)$-cables of $(2,3)$-torus 
    knots.  Concentric circles indicate multiplicities, {\em ie}, the number of 
    distinct isotopy classes with a given $r$ and $\tb$.}   
  \label{iterate}
\end{figure}

There are knot types that are not uniformly thick. As mentioned above the unknot is not
uniformly thick, but it acts somewhat as if it is. For example its cables (torus knots)
are Legendrian simple. It would appear that all positive torus knots are not Legendrian
simple, but this seems difficult to prove in general. We have the following:
\begin{thm}[Etnyre and Honda, \cite{EtnyreHonda??}]\label{cableex}
If $\K'$ is the $(2,3)$-cable of the $(2,3)$-torus knot, then 
$\L(\K')$ is classified as in Figure~\ref{iterate}.    This entails 
the following:
\be
\item There exist exactly two maximal Thurston-Bennequin 
representatives $K_\pm\in \L(\K')$.  They satisfy $\tb(K_\pm)=6$, 
$r(K_\pm)=\pm 1$. 
\item There exist exactly two non-destabilizable representatives 
$L_\pm\in \L(\K')$ which have non-maximal Thurston-Bennequin 
invariant. They satisfy $\tb(L_\pm)=5$ and $r(L_\pm)=\pm 2$. 
\item Every $L\in \L(\K')$ is a 
stabilization of one of $K_+$, $K_-$, $L_+$, or $L_-$. 
\item $S_+(K_-)=S_-(K_+)$,  $S_-(L_-)= S_-^2(K_-)$, and 
$S_+(L_+)=S_+^2(K_+).$
\item $S_+^k(L_-)$ is not (Legendrian) isotopic to $S_+^k S_-(K_-)$ 
and $S_-^k(L_+)$ is not isotopic to $S_-^k S_+(K_+)$, for all positive 
integers $k$.    Also, $S_-^2(L_-)$ is not isotopic to $S_+^2(L_+)$.
\ee \end{thm}

This example is quite interesting since it is the first known example where there 
are Legendrian knots which do not have the maximal $tb$ for the knot type and yet
they do not destabilize. Moreover, this is the first knot type where there are Legendrian
representatives with the same invariants that do not become Legendrian isotopic after
some number of stabilizations of a fixed sign. See the section on transverse knots
for more discussion of this knot type.

\subsection{Links}
\begin{figure}[ht]
{\epsfxsize=4.5in \centerline{\epsfbox{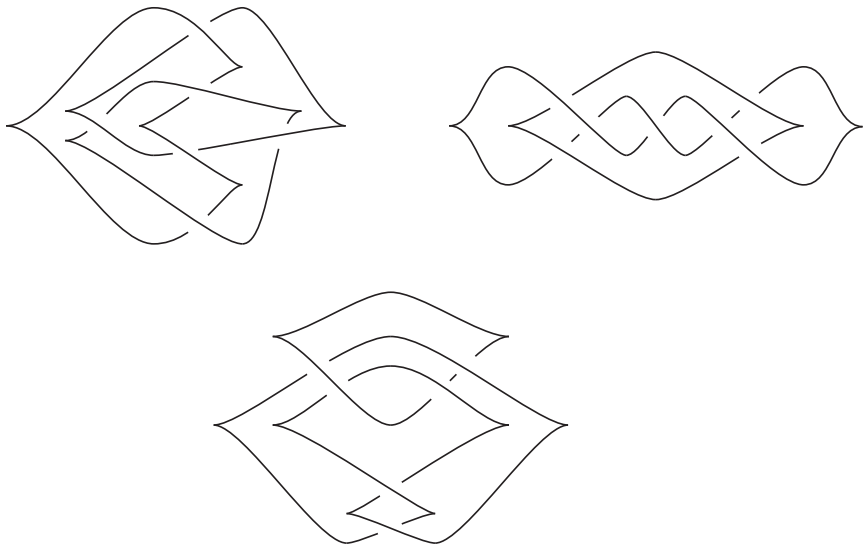}}
\caption{Legendrian Borromean rings with $tb=-4$ (top left), Legendrian Whitehead link with $tb=-5$
(top right) and Legendrian Whitehead mirror realized by two $tb=-1$ unknots (bottom).}
\label{llink}} 
\end{figure}

There has not been much study of Legendrian links.  One of the first result about links appeared
in \cite{Mohnke01} and addressed the question: Given a link type is there a Legendrian realization
such that each component has maximal Thurston--Bennequin invariant in its knot type? The answer is NO.
For a link $L=\cup L_i$ define its Thurston--Bennequin invariant to be the sum of the 
Thurston--Bennequin invariants of 
its components $tb(L)=\sum tb(L_i).$  Then the inequality for the Thurston--Bennequin invariant
involving the HOMFLY polynomial
in Theorem~\ref{otherbounds} still holds. For the Borromean rings $B$ and the Whitehead link $W$
a simple computations of the HOMFLY polynomial yield
\[tb(B)\leq -4\]
and 
\[tb(W)\leq -5.\]
Since $B$ had 3 components and $W$ has 2 we clearly cannot realized them with components consisting
of unknots with maximal Thurston--Bennequin invariant. These bounds are realized, see Figure~\ref{llink}.
Also in this figure it is shown that the mirror of the Whitehead link can be realized by Legendrian unknots 
of $tb=-1,$ once again illustrating the ``cirality'' of Legendrian knots.
In general, very little is known concerning realizing Legendrian links with given Legendrian knot components and
the restrictions on the Thurston--Bennequin invariants.

\begin{figure}[ht]
{\epsfxsize=4in \centerline{\epsfbox{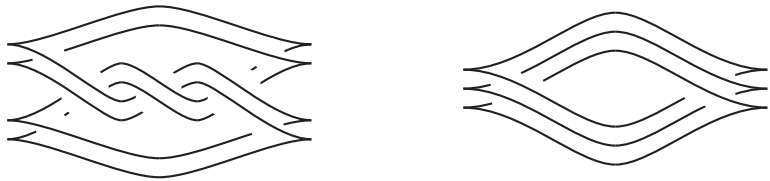}}
\caption{A 2-copy of the right handed trefoil (right) and a 3-copy of the unknot (left).}
\label{Ncopy}} 
\end{figure}

The $N$-copy of a Legendrian knot has also been studied. Given a Legendrian knot $L$ in
$(\R^3,\xi_{std}),$ let $L_N$ be the Legendrian link obtained by talking $N$ copies of $L$ each shifted 
slightly in the $z$ direction. See Figure~\ref{Ncopy}. 
The link $L_N$ is a Legendrian link in the topological link type of the
$(Ntb(L),N)$-cable of $L.$ If $L$ had maximal $tb$ for Legendrian knots in its knot type then $L_N$ is
a Legendrian realization of its link type with each component having maximal $tb$ in its knot type.

Topologically one can always realize any arbitrary permutation of the components of this link via 
isotopy. Can this be done for the components of the Legendrian $N$-copy.
We have the following result.
\begin{thm}[Michatchev 2001, \cite{Michatchev01}]
If $L$ is a Legendrian unknot, then only cyclic permutations of the components of $L_N$ are possible via 
Legendrian isotopy.
\end{thm}
One may easily check that cyclic permutations are possible. To show other permutations are not possible
one must use an enhanced version of contact homology see \cite{Michatchev01}.

One can study other cables as well. The following result follows easily from \cite{EtnyreHonda01b}.
\begin{thm}
Let $p>q>0$ be relatively prime. There is a unique Legendrian link with maximal Thurston--Bennequin
invariant in the $(np,nq)$-torus link type. Moreover,
each component of the link has maximal $tb$ in its knot type and any permutation of its components can be
realized by a Legendrian isotopy.
\end{thm}

\subsection{The homotopy type of the space of Legendrian knots.}
\begin{figure}[ht]
  \small{\epsfxsize=4in \centerline{\epsfbox{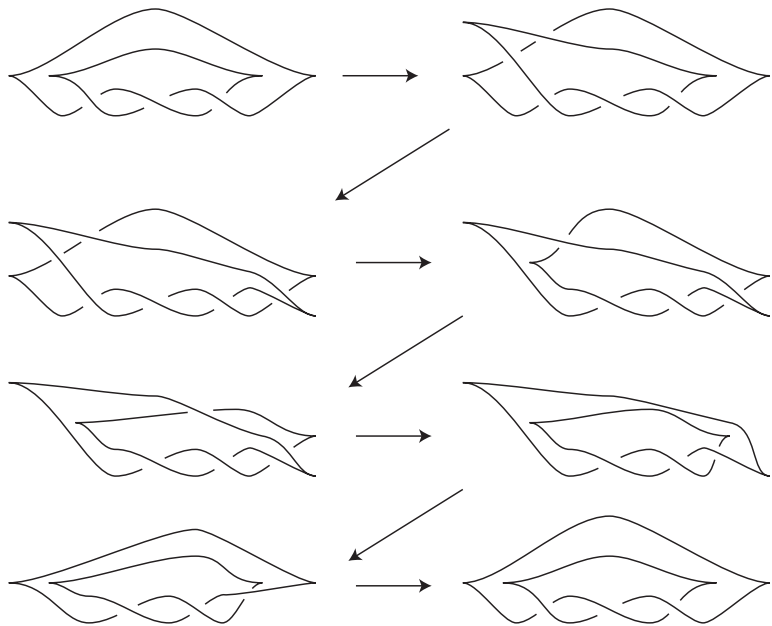}}}
  \caption{A loop of Legendrian trefoils.}
  \label{loop}
\end{figure}

Recall $\K$ denotes a topological knot type. That is \K is a space of topological knots isotopic to
a fixed knot.
Also recall
$\L(\K)$ denotes the set of all Legendrian knots in $\K.$  
All the
classification results above can be thought of in terms of identifying the kernel of the map
\[i_*:\pi_0(\L(\K))\to \pi_0(\K),\]
where $i:\L(\K)\to\K$ is the natural inclusion map. In general, studying the homotopy groups
of $\L(\K)$ is an interesting problem. To this end, fix a component $\L$ of $\L(\K)$ and consider
the inclusion $i:\L\to \K.$ Recently K\'alm\'an \cite{Kalman1, Kalman2} has shown the induced
map on $\pi_1$ need not be injective. Specifically, if $\L$ is the set of maximal Thurston-Bennequin invariant,
Legendrian right handed trefoil knots then consider the loop $\Omega$ in $\L$ illustrated
in Figure~\ref{loop}.
\begin{thm}[K\'alm\'an \cite{Kalman1, Kalman2}]
The loop $\Omega^6$ is non-trivial in $\pi_1(\L)$ but $i(\Omega^6)$ is trivial in $\pi_1(\K).$
\end{thm}
The key ingredient in proving this theorem is the following result.
\begin{thm}[K\'alm\'an \cite{Kalman1}]
If $L$ is a fixed (generic) Legendrian knot in $\L$ then there is a multiplicative homomorphism 
\[\mu : \pi_1(\L, L)\to Aut(CH(L))\]
defined by continuation on the Chekanov-Eliashberg contact homology. 
\end{thm}
Now showing that $\Omega^6$ is a non-trivial loop amounts to a computation that $\mu$ of
this loop is non-trivial. This is carried out in \cite{Kalman1}. Moreover, in \cite{Kalman2}
non-trivial loops are found for all maximal Thurston-Bennequin invariant positive torus knots.

\subsection{Transverse knots}

Recall a knot type \K is transversely simple if transverse knots in the knot type are determined by
their self-linking numbers. By Theorem~\ref{nslist} we know that the negative stable classification of 
Legendrian knots
implies the transverse classification. Thus from the results above we know the following knot types are 
transversely simple
\begin{enumerate}
\item the unknot (\cite{Eliashberg93}),
\item torus knots (positive torus knots in \cite{Etnyre99}, all torus knots in \cite{EtnyreHonda01b,Menasco01, Menasco?}), and
\item the figure eight knot (\cite{EtnyreHonda01b}).
\end{enumerate}
The theorem on Legendrian connect sum implies that the connect sum of Legendrian simple knot types need
not be Legendrian simple, but this is not the case for transverse knots.
\begin{thm}[Etnyre and Honda 2003, \cite{EtnyreHonda03}]
If $\K_1$ and $K_2$ are transversely simple then $\K_1\#\K_2$ is also transversely simple.
\end{thm}
This theorem easily follows from Theorem~\ref{thm:con} on 
Legendrian connect sums and the observation 
if a knot
type is transversely simple then its $tb$-$r$ mountain range 
(see Figure~\ref{fig:geog}) is connected via negative stabilizations. The proof is
an easy exercise or see \cite{EtnyreHonda03}.

\begin{figure}[ht]
  \relabelbox 
  \small{\epsfxsize=3in \centerline{\epsfbox{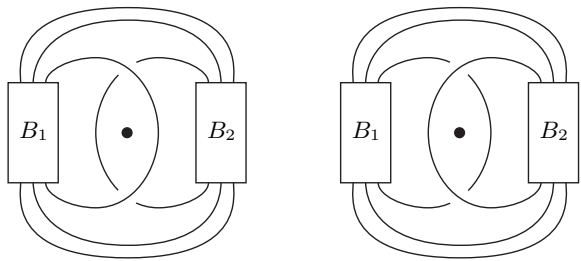}}}
  \relabel {1}{$B_1$}
  \relabel {2}{$B_2$}
  \relabel {3}{$B_1$}
  \relabel {4}{$B_2$}
  \endrelabelbox
  \caption{Exchange move.}
  \label{exmove}
\end{figure}

Once again there is an interesting relation between transverse knots and braids. We call a knot type \K
\dfn{exchange reducible} if any braid representing \K can be reduced to a unique minimal braid index braid
for \K by a sequence of braid destabilizations and exchange moves. An exchange move is shown in Figure~\ref{exmove}.
\begin{thm}[Birman and Wrinkle 2000, \cite{BirmanWrinkle00}]
If \K is exchange reducible then it is transversely simple.
\end{thm}
To prove this theorem one needs to observe that there are two types of braid destabilizations. A positive one
and a negative one. The positive one can be realized by a transverse isotopy and the negative one corresponds
to a transverse destabilization of the transverse knot. Moreover, an exchange move can also be done by a 
transverse isotopy of the knot. Thus being exchange reducible implies that any transverse knot destabilizes
to a unique maximal self-linking number transverse knot and thus the knot type is transversely simple. 
One can easily reprove that the unknot is transversely simple with this theorem \cite{BirmanMenasco92}, see
also \cite{Menasco01, Menasco?} for a discussion of torus knots and cables from this perspective.

The search for non-transversely simple knot types was long and difficult, due in part to the lack of
non classical invariants for transverse knots. But we now know there are non-transversely simple knot types.
\begin{thm}[Birman and Menasco, \cite{BirmanMenasco??}]
The transverse knots represented by the closures of the braids 
$\sigma_1^3\sigma_2^4\sigma_1^2\sigma_2^{-1}$ 
and $\sigma_1^3\sigma_2^{-1}\sigma_1^2\sigma_2^4$ are not transversely isotopic, but
have the same self-linking number.
\end{thm}
Many other
examples arise from Birman and Menasco's analysis, see \cite{BirmanMenasco??, BirmanMenasco93}, 
but they are all three braids. They prove the knots are
not transversely isotopic by classifying all closed three braids and understanding how to pass from one
three braid to another representing the same knot type. (They can also show that there must exist non transversely 
simple knots types that are not three braids, but the specific knot types cannot identified at this time.)

The classification of Legendrian knots in the knot type \K of the $(2,3)$-cable of 
the $(2,3)$-torus knot in 
Theorem~\ref{cableex} implies that \K is 
not transversely simple. In particular we get the first classification of 
transverse knots in a non-transversely simple knot type.
\begin{thm}[Etnyre and Honda, \cite{EtnyreHonda??}]
Let \K be the $(2,3)$-cable of the $(2,3)$-torus knot. Each odd integer less than or equal to seven 
is the self-linking number of a transverse knot in $\T (\K).$ Moreover, the self-linking number determines
the transverse knot if $sl\not= 3$ and there are precisely two transverse knots with $sl=3.$ 
\end{thm}
One should be able to construct many more such examples of cables of positive torus knot
using the ideas from \cite{EtnyreHonda??}.

\subsection{Knots in overtwisted contact structures}
Overtwisted contact structures have been classified by Eliashberg in \cite{Eliashberg89} 
and are determined by their
homotopy type of plane field. Since \cite{Eliashberg89} overtwisted structures have not been studied
much. In particular,
Legendrian knots in overtwisted contact structures have been somewhat ignored. There are however a few
things we do know.

If $L$ is a Legendrian knot in an overtwisted contact manifold $(M,\xi)$ then call $L$ \dfn{loose} if
the contact structure on $M\setminus L$ is overtwisted. 
While it has been claimed by several authors that 
if a Legendrian knot $L$ violates the Bennequin inequality then it is loose, this is not actually
true as can be seen from the example below (this was first observed in \cite{Dymara??}). Concerning the classification of loose knots we have the
following theorem.
\begin{thm}[Eliashberg and Fraser 1995, \cite{EliashbergFraser}: and Dymara 2001, \cite{Dymara01}]
Fix an overtwisted contact structure on $S^3,$ then two loose Legendrian knots in the same knot type and having the same
$tb$ and $r$  are contactomorphic. Moreover, if there is a fixed overtwisted disk in the complement of both 
Legendrian knots then they are Legendrian isotopic.
\end{thm}
The proof of this theorem relies on the classification of overtwisted contact structures in \cite{Eliashberg89}. 
In particular, $tb$ and $r$ determine the homotopy type
of the contact structure on the complement of $L.$ Thus two loose knots with the same invariants are related by
a contactomorphism of the ambient contact manifold. One gets a global contact isotopy if there is a fixed
overtwisted disk in the complement of both Legendrian knot. 

In \cite{EliashbergFraser} it was asked if non-loose 
knots exist in overtwisted contact structures. In fact they do, as shown
in \cite{Dymara01}. 
The example given there goes as follows. The standard tight contact structure on $S^3$ can be thought
to be the planes orthogonal (in the standard round metric) to the Hopf fibration. 
Now let $\xi$ be the contact
structure obtained by a ``half-Lutz twist.'' To understand this write the 
Hopf fibration as $\pi:S^3\to S^2.$
Let $A\subset S^2$ be an embedded annulus $A=S^1\times [0,1].$ 
Let $N=\pi^{-1}(A)$. Thus $N=A\times S^1.$ where the $S^1$ are the
Hopf fibers. The standard contact structure intersects the tori $T_t=S^1\times\{t\}\times S^1$ in a 
linear foliation.
On $N$ alter the contact structure by adding an extra half twist, that is leave the contact structure
fixed at the boundary of $N$ but as $t$ move from $0$ to $1$ make the slopes of the 
foliations induced on the
$T_t$ change more rapidly so that the angles swept out by the foliations are $\pi$ more that then angles
swept out by the standard contact structure. Call this new contact structure $\xi.$ Note there will be a $t$
such that the characteristic foliation on $T_t$ is parallel the the fibers in the Hopf fibration. 
Let $\gamma$ 
be a leaf in this characteristic foliation. In \cite{Dymara01} 
it was shown that the overtwisted contact structure
$\xi$ is tight when restricted to $S^3\setminus \gamma.$ Intuitively this is because $\gamma$ intersects all
the obvious overtwisted disks. To prove the complement is tight one needs to examine the universal cover,
$\R^3,$ of
the complement and prove that the contact structure pulls back the standard contact structure on $\R^3.$  
Note $\gamma$ is a Legendrian unknot and one may easily compute (in a manner similar to the computation
for torus knots) that $tb=1$ and $r=0.$ 
It is easy to find another Legendrian $\gamma'$ with $tb=1$ and $r=0$ so that the complement is still 
overtwisted
(just take the connected sum the boundaries of two disjoint overtwisted disks). Now $\gamma$ and $\gamma'$
cannot be Legendrian isotopic since their complements are not contactomorphic. 

From this example it seems the classification of Legendrian unknots in $(S^3,\xi)$ is more complicated the
the classification of Legendrian unknots in tight contact structures. It would be quite interesting  (and 
probably not too hard) to completely classify Legendrian unknots in this contact structure, 
or more generally
all overtwisted contact structures on $S^3.$ It seems studying Legendrian unknots in overtwisted contact 
structures
might be very interesting. But this must wait for future work.

\section{Higher dimensions}\label{sec:hd}
A contact structure on a $2n+1$ dimensional manifold $M$ is a hyperplane field $\xi$ that
can be locally given as the kernel of a 1-form $\alpha$ such that $\alpha\wedge(d\alpha)^n\not=0.$
As in dimension three one should think of such a plane field as ``maximally nonintegrable''.
The standard contact structure $\xi$ on $\R^{2n+1}$ is given as the kernel of 
\begin{equation}
        \alpha=dz-\sum_{j=1}^n y_j dx_j,
\end{equation}
where $x_1,y_1,\ldots, x_n,y_n, z$ are Euclidean coordinates on $\R^{2n+1}.$ In this 
section we restrict attention to $(\R^{2n+1},\xi).$ 
A \dfn{Legendrian knot} (or Legendrian
submanifold) of $\R^{2n+1}$ is an $n$ dimensional submanifold $L\subset \R^{2m+1}$ that is tangent to
$\xi$ at each point.  We will frequently use $L$ to
denote the Legendrian submanifold and the domain of a Legendrian embedding.
We also recall that 
the standard symplectic structure on $\C^n$ is given by 
\begin{equation}
        \omega=\sum_{j=1}^n dx_j\wedge dy_j,
\end{equation}
and a Lagrangian submanifold is an $n$ dimensional submanifold $L\subset \C^n$ for which 
$\omega(v_1, v_2)=0$ for all vectors $v_1,v_2\in TL.$ 

\subsection{Legendrian knots in $\R^{2n+1}$}\label{highknot}
As in three dimensions there are two standard 
projections that are useful in studying Legendrian knots. The first is
the \dfn{Lagrangian projection} (sometimes called the complex projection) 
which projects out the $z$ coordinate:
\begin{equation}
        \pi:\R^{2n+1}\to \C^n: (x_1,y_1,\ldots,x_n,y_n,z)\mapsto (x_1,y_1,\ldots,x_n,y_n).
\end{equation}
The second is called the \dfn{front projection} and it projects out the $y_j$'s:
\begin{equation}
        \Pi:\R^{2n+1}\to \R^{n+1}: (x_1,y_1,\ldots,x_n,y_n,z)\mapsto (x_1,\ldots,x_n, z).
\end{equation}

We begin with the Lagrangian projection.
If $L$ is an embedded Legendrian knot then
$\pi (L)$ is an immersed Lagrangian submanifold of $\C^n.$ 
If $L$ is generic, that is in a $C^\infty$ dense subset of all Lagrangian embeddings, 
then $\pi (L)$ will have a finite number of isolated
double points. 
The embedding of $L$ in $\R^{2n+1}$ can be recovered (up to rigid translation
in the $z$ direction) from $\pi(L)$ as in the $n=1$ case described above. Specifically,
pick a point $p\in \pi(L)$ and choose any $z$ coordinate
for $p$ then the $z$ coordinate of any other point $p'$ is determined by 
\begin{equation}\label{1recoverz}
        \sum_{j=1}^n \int_\gamma y_j dx_j, 
\end{equation}
where $\gamma$ is any path in $\pi(L)$ from $p$ to $p'.$ 
Furthermore, given any Lagrangian immersion in $\C^n$ 
with isolated double points, if the integral in Equation~\ref{1recoverz} is independent
of the path $\gamma$ then we obtain a Legendrian immersion
in $\R^{2n+1}.$ Note that we will get an embedding as long as the above integral is
not zero for paths connecting the double points. The integral in Equation~\eqref{1recoverz} will be
independent of $\gamma$ when we have
exact Lagrangian immersions $i:L\to\R^{2n},$ {\em i.e.}  ones for which  $i^*(\sum_{j=1}^n x_j dy_j)$ is
exact. In particular, if $H^1(L)=0,$ as is the case for $S^n$ if $n\not=1,$ 
then all Lagrangian immersions are exact. Thus an isotopy $L'_t$ of exact Lagrangian manifolds in $\C^n$
with transverse double points will lift to an isotopy $L_t$ of Legendrian manifolds in
$\R^{2n+1}$ (again there is no guarantee that the manifolds will be embedded). 
\bex\label{1basicexample}
Here we consider the most basic example of an exact Lagrangian immersion of a sphere in $\C^n,$
thus providing a Legendrian sphere in $\R^{2n+1}.$ Let $S^n=\{(x,y)\in\R^n\times\R: |x|^2+y^2=1\}$ and
then define
\begin{equation}
        f(x,y):S^n\to \C^n:(x,y)\mapsto ((1+iy)x).
\end{equation}
We claim that the image of $f$ is an exact Lagrangian sphere in $\C^n$ with one double point that
lifts to an embedded sphere in $\R^{2n+1}.$ When $n=1$ the image is
a figure eight in the plane with a double point at the origin.
\eex
The Reidemeister type theorem here is
\begin{lem}\label{1isotopy} 
If two Legendrian knots in $\R^{2n+1}, n>1,$ are Legendrian isotopic then their Lagrangian projections
are related by regular homotopies and isolated double point moves (see the left hand side of
Figure~\ref{fig:lagreid}).
\end{lem}


Though the Lagrangian projection of Legendrian knots has many nice properties, for example 
Lemma~\ref{1isotopy}, it has some drawbacks as well. Specifically, given a submanifold
$L^n$ of $\C^n$ it is not clear whether or not it is an exact Lagrangian submanifold
or lifts to an embedded Legendrian knot.
Thus it is sometimes useful to consider the front projection of a Legendrian knot.
The front projection $\Pi(L)$ is a codimension one subvariety of $\R^{n+1}.$ Living in a lower
dimensional space makes it somewhat easier to visualize. 
The drawback to the front projection
is that $\Pi(L)$ will have certain singularities. In fact, any singularity with a well defined tangent 
space one would expect from a projection
of $L^n\subset \R^{2n+1}$ to $\R^{n+1}$ can occur 
except tangential double points (these would correspond to 
double points of
$L$ in $\R^{2n+1}$). But any map from an $n$ dimensional space to
$\R^{n+1}$ with no vertical tangencies and 
the appropriate type of singularities
will lift to a Legendrian knot in $\R^{2n+1}.$

To see that any map into $\R^{n+1}$ with the appropriate singularities can be lifted to a Legendrian
knot we observe that at cusp edges (and all other allowable singularities) there is
a well defined tangent plane to the image of the map. Thus we can recover the $y_i$ coordinate
of the Legendrian knot by looking at the slope of the tangent plane in the $x_iz$-plane 
(since for the lift to be Legendrian its tangent planes have to satisfy $dz=\sum_{i=1}^n y_i dx_i$).

\begin{figure}[ht]
        {\epsfxsize=4in \centerline{\epsfbox{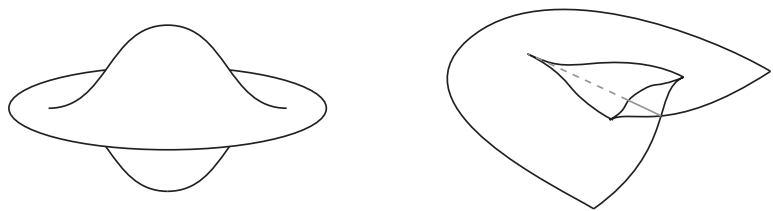}}}
        \caption{Front projection with a cusp edge (right) and a ``swallow-tail'' singularity (left).}
        \label{1fig:front}
\end{figure}

If $\phi:L\to \R^{2n+1}$ is a Legendrian embedding then (after a $C^\infty$ small perturbation
to another Legendrian) off of a codimension one subset
$\Sigma\subset L$ the map $\Pi \circ\phi:L\to\R^{n+1}$ will be an immersion. 
For a generic point $p$ in $\Sigma$ we can find a neighborhood $N$ of $p$ in $L,$ with coordinates
$(x_1,\ldots, x_n),$ and a coordinate neighborhood $N'$ of $\Pi\circ\phi(p)$ 
such that $\Pi_F\circ\phi$ is expressed in these coordinates
by $(x_1,\ldots, x_n)\mapsto(3x_1^2,x_2, \ldots, x_n, \pm 2x_1^3).$ We say such a point is on a 
\dfn{cusp edge}. See Figure~\ref{1fig:front}.
Note that the image of
$\Pi\circ \phi$ cannot have tangent planes containing the $z$-direction (since $\phi(L)$ is Legendrian).
It is the cusp edges that allow our front projections to avoid these vertical tangencies.  For a 
generic one dimensional Legendrian knot only cusp edges occur in the front projection.  In 
higher dimensions there can be other singularities as well, however they all occur in codimension
larger than one. So the top dimensional strata of $\Sigma$ consists entirely of cusp edges. 
In particular, later it will be important that any path on $L$ can be perturbed so 
that it only intersects
$\Sigma$ in cusp edges. For a thorough discussion of allowable singularities see \cite{AG}.

For Legendrian knots $L$ in $\R^{2n+1}$ we have two ``classical'' invariants. First,
Tabachnikov \cite{Tabachnikov88} has shown how to generalize the Thurston--Bennequin invariant to 
higher dimension. Specifically, let $L'$ be a copy of $L$ obtained by pushing $L$ slightly in the
$z$-direction, then 
\beq
tb(L)=\text{linking}(L,L').
\eeq
It is interesting to note that when $n$ is even $tb(L)=\frac12 \chi(L)$ and thus is a trivial invariant
of the Legendrian isotopy class of $L.$ When $n$ is odd $tb(L)$ is a non trivial invariant (see 
Section~\ref{highex}).

It is more difficult to generalize the rotation number a Legendrian knot in three manifolds to
higher dimensions. To each point $x\in L$ we get a Lagrangian plane $d\pi_x(T_xL)$
in $\C^n.$ Thus we get a map $L$ to $\text{Lag}(\C^n)$ the space of Lagrangians in $\C^n.$
The homotopy class of this map is an invariant of $L.$ In dimension three, $\pi_1(\text{Lag}(\C))=\Z$
and thus to Legendrian knots in three manifolds we have a rotation number. In higher dimensions and
when $L$ is not just a sphere this homotopy class can be much more complicated than an integer.
Again it is interesting to note that when $n$ is even and $L$ is an $n$ sphere then this rotation
class is always trivial. Thus when $n$ is even there are no classical invariants of Legendrian knots.

\subsection{Generalizations of the Chekanov-Eliashberg DGA}\label{highDGA}
Contact homology naturally generalizes to higher dimensions. Here we define the
DGA associated to a Legendrian knot in $\R^{2n+1}.$ 

\noindent
{\bf The Algebra.} Denote the double points of the Lagrangian projection, $\pi(L),$ 
by $\mathcal{C}.$ We assume $\CC$ is a finite set 
of transverse double points.
Let $\mathcal{A}$ be the free associative unital algebra over $\Z_2$ generated
by $\mathcal{C}.$\hfill\break
{\bf The Grading.} To each crossing $c\in\CC$ there are two points $c^+$ and $c^-$ in 
$L\subset\R^{2n+1}$ that
project to $c.$ We denote by $c^+$ the point with larger $z$-coordinate.  
Choose a map $\gamma_c:[0,1]\to L$ that parametrizes an arc running from $c^+$ to 
$c^-.$ (Note there could be more that one path.)
For each point $\gamma(\theta)\in L$ we have a Lagrangian plane 
$d\pi_{\gamma(\theta)}(T_{\gamma(\theta)}L)$ in $\C^n.$ Thus $\gamma$ give us a path 
$\widehat{\gamma}$ in
$\text{Lag}(\C^n).$ Since $c$ is a transverse double point   
$\widehat{\gamma}(0)$ is transverse to $\widehat{\gamma}(1).$ Thus we can find a complex structure
$J$ on $\R^{2n}$ that (1) induces the same orientation on $\C^n$ as the standard complex structure
and (2) $J(\widehat{\gamma}(1))=\widehat{\gamma}(0).$ Now set $\overline{\gamma}(\theta)=
e^{J\theta}\widehat{\gamma}(1).$ Note that $\widehat{\gamma}$ followed by $\overline{\gamma}$ is 
a closed loop in $\text{Lag}(\C^n).$ Moreover $\pi_1(\text{Lag}(\C^n)=\Z$ thus to this closed loop
we get an integer, $cz(c),$ the Conley-Zehnder invariant of $c.$  The grading on $c$ is
\[|c|=cz(c)-1.\]
We note that $c$ in general depends on the path $\gamma$ with which we started. To take care of
this ambiguity note that to each circle immersed in $L$ we get an integer through the above procedure.
Let $n$ be the greatest common divisor of all these integers.  It is easy to convince oneself that
$|c|$ is well defined modulo $n.$ \hfill\break
{\bf The Differential.} We will define the differential $\partial$ on $\A$ by defining it on
the generators of $\A$ and then extending by the signed Leibniz rule:
\[\partial ab= (\partial a)b+(-1)^{|a|}a\partial b.\]
Let $a\in \CC$ be a generator of $\A$ and let $b_1\ldots b_k$ be a word in the ``letters'' \CC.
Let $P_{k+1}$ be a $k+1$ sided polygon in $\C$ with vertices labeled counterclockwise $v_0,\ldots, v_k.$
We will consider maps $u:(P_{k+1},\partial P_{k+1})\to (\C^n, \pi(L))$ such that 
$u|_{\partial P_{k+1}\setminus \{v_i\}}$ lifts to a map to $L\subset \R^{2n+1}.$  
Call a vertex $v_i$  mapping to the double point $c$ is \dfn{positive} (resp. \dfn{negative})
if the lift of the arc just clockwise of $v_i$ 
in $\partial P_{k+1}$ lifts to an arc approaching $c^+$ (resp. $c^-$) and the arc just counterclockwise
of $v_i$ lifts to an arc approaching $c^-$ (resp. $c^+$), where $c^\pm$ are as in the definition of
grading.
Set 
\[\M^a_{b_1\ldots b_k}=
\{u:(P_{k+1},\partial P_{k+1})\to (\C^n, \pi(L)) \text{ such that $u$ satisfies 1.-- 4. below}\}/ \sim\]
where $\sim$ is holomorphic reparameterization and the conditions are
\begin{itemize}
\item[1.] $u|_{\partial P_{k+1}\setminus \{v_i\}}$ lifts to a map to $L\subset \R^{2n+1},$  
\item[2.] $u(v_0)=a$ and $v_0$ is positive,
\item[3.] $u(v_i)=b_i, i=1,\ldots, k,$ and $v_i$ is negative,
\item[4.] $u$ is holomorphic.
\end{itemize}
We can now define
\[\partial a=\sum_{b_1\ldots b_k} (\#_2\M) b_1b_2\ldots b_k,\]
where the sum is taken over all words in the letters \CC  for which 
$\dim(\M^a_{b_1\ldots b_n})=0$ and $\#_2$ denotes the modulo two count of
elements in $\M.$ 
\begin{thm}[Ekholm, Etnyre and Sullivan, \cite{EES1, EES2}]
With the notation above:
\begin{enumerate}
\item The map $\partial$ is a well defined differential that reduces the grading by 1.
\item The stable tame isomorphism class of $(\A,\partial)$ is an invariant of $L.$
\item The homology of $(\A,\partial)$ is an invariant of $L.$
\end{enumerate}
\end{thm}
We will compute some simple examples below.

\subsection{Examples}\label{highex}
The DGA of a Legendrian knot $L$ is defined in terms of the Lagrangian projection, but
it will be much easier to construct examples in the front projection. Thus we would like
to be able to recognize the double points in the Lagrangian projection by looking
at the front projection. We begin by noting that generically we can assume
that the double points in the Lagrangian projection do not occur along singularities
in the front projection. So we consider a point $p\in \Pi(L)$  in the front projection
that is not singular. Since $\Pi(L)$ is embedded near $p$ and the tangent plane does not
contain the vector $\frac{\partial}{\partial z}$ we can find an open set $U$ in $\R^n$ 
and a function $f:U\to \R$ such that near $p,$ $\Pi(L)$ is given by the graph of $f.$ 

\begin{figure}[ht]
        {\epsfxsize=4in \centerline{\epsfbox{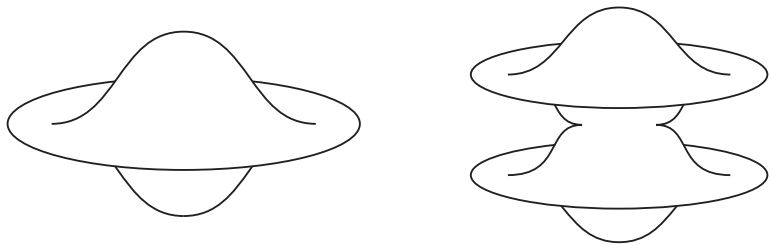}}}
        \caption{The front projection of two Legendrian spheres.}
        \label{fig:2h}
\end{figure}

Let $c^+$ and $c^-$ be two points in $\Pi(L).$ If they correspond to double point in the Lagrangian
projection then they must have the same $x_i$ coordinates, which we denote $x,$
and different $z$ coordinates. We 
assume the $z$ coordinate of $c^+$ is larger. Let $f^\pm$ be the functions described 
above for $c^\pm.$ We can assume that $f^+$ and $f^-$ have the same domain. Let $f=f^+-f^-.$
The following are equivalent:
\begin{itemize}
\item $c^+$ and $c^-$ correspond to a double point in $\pi(L)$
\item $T_{c^+}\Pi(L)=T_{c^-}\Pi(L)$
\item $df^+_x=df^-_x$ 
\item $df_x=0$
\end{itemize}
Thus we have a double point if $x$ is a critical point of $f.$ Moreover, the double point
is transverse if and only if $x$ is a nondegenerate double point. Thus to a transverse
double point we have a Morse index, $\text{Ind}_x(f).$ Suppose $c^+, c^-$ correspond to
a double point. Let $\gamma$ be a path from $c^+$ to $c^-$ that is transverse to the singular
set of $\Pi(L).$ Thus it intersects the singular set along the cusp edges in a finite number
of points. Let $D$ be the number of times $\gamma$ intersects the cusps while its $z$ coordinate
is decreasing and $U$ the number of times with its $z$ coordinate increasing.
\begin{lem}[Ekholm, Etnyre and Sullivan, \cite{EES1}]
With the notation as above. If $c^+$ and $c^-$ correspond to a transverse double point $c$ then
\beq
cz(c)= \text{Ind}_x(f)+D-U.
\eeq
\end{lem}
We are now ready for some examples.
\bex
In Figure~\ref{fig:2h} the front projections of two Legendrian knots are shown. Though
the pictures are of Legendrian $S^2$'s in $\R^5$ there are clearly analogous Legendrian
$S^n$'s in $\R^{2n+1}.$ We use contact homology to distinguish these examples
in all dimensions. When $n$ is even they cannot be distinguished by any classical invariants.
When $n$ is odd they have different Thurston-Bennequin invariants.
Note the Legendrian $L$ on the left has only one double point (the ``axis'' of the flying saucer) $c.$
Using the above lemma we see that $cz(c)=n+1$ so $|c|=n.$ Since $\partial$ is degree $-1$ 
we see that $\partial c=0.$ So the contact homology is generated by one element in grading $n.$

For $L'$ on the right of the figure, there is a double point $c$ corresponding the the ``axis'' again, but
now one must go down $3$ cusp edges from $c^+$ to $c^-$ thus $|c|=n+2.$ 
If the picture is drawn symmetrically, there are an $S^{n-1}$'s worth of double point. By deforming
the picture slightly so the cusps are tilted we can break the $S^{n-1}$'s worth of double points
into two double point $a$ and $b.$ They will have gradings $|a|=n$ and $|b|=1.$  We can again
conclude that $\partial c=0$ but we do not know if $\partial b=0$ or 1, and when $n=2$ if $\partial a
=b$ or $0.$ But in any case the contact homology of $L'$ will be different from that of
$L.$ 
\eex

We now describe a procedure for changing a given Legendrian knot.
Suppose $L$ is a Legendrian knot in $\R^{2n+1}$ containing two disks
$L_u$ and $L_l$ that project to the same disk under the Lagrangian projection.
This is a degenerate situation, but it is easy to isotop any Legendrian knot so
that it has such disks. Assume $L_u$ is above $L_l.$ Let $M$ be a $k$-manifold embedded in 
$L_l.$ There is an ambient isotopy of $\R^{n+1}$ that is supported near $\Pi(M)\times \R,$ where the
$\R$ factor is the $z$-axis, that moves $\Pi(M)$ up just past $\Pi(L_u).$ 
In the front projection replace $\Pi(L_l)$ with the disk obtained from $\Pi(L_l)$ by applying
the ambient isotopy. This new front diagram describes a new Legendrian knot $L_M$ called the
\dfn{stabilization of $L$ along $M$}. One may show
\begin{prop}[Ekholm, Etnyre and Sullivan, \cite{EES1}]
If $L_M$ is the stabilization of $L$ along $M$ and the notation is as above then
\begin{enumerate}
\item The rotation class of $L_M$ is the same as the rotation class of $L.$
\item If $D,U$ are the number of down, up, cusps along a generic path from $\Pi(L_u)$ to $\Pi(L_l)$ 
in $\Pi(L)$ then
\[
\tb(L_M)=
        \begin{cases}
        \tb(L), &  n \mbox{  even}, \\
        \tb(L)+ (-1)^{(D-U)}2\chi(M), & n \hbox{ odd,}
    \end{cases}
\]
\item If the stabilization of the Legendrian knot $L$ takes place in a small neighborhood
        of a cusp edge of $\Pi(L)$ then $CH_*(L)=0.$ 
\end{enumerate}
\end{prop}

\bbr
Note using this proposition it is easy to see that when $n$ is odd and not 1 then we can realize any
any integer, of the appropriate parity, as the $tb$ of some Legendrian knot.
\eer

\begin{figure}[ht]
  \relabelbox 
  \small{\epsfxsize=4.2in\centerline{\epsfbox{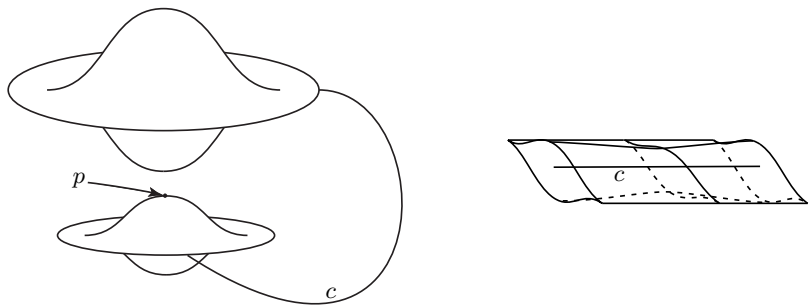}}}
  \relabel {p}{$p$}
  \relabel {c}{$c$}
  \relabel {1}{$c$}
  \endrelabelbox
  \caption{The Legendrian knot $L$ (left). To form $L$ a neighborhood of the arc $c$ should
    be replaced with the tube shown on the right.}
  \label{fig:bs}
\end{figure}

It is important that stabilization does not always force the contact homology to be trivial.
\bex
Let $L$ be the Legendrian knot shown in Figure~\ref{fig:bs}.
 It is clearly Legendrian isotopic to
the standard ``flying saucer'' in Figure~\ref{fig:2h}. Let $L_1$ be the stabilization of 
$L$ along the point $p$ indicated in the figure. Let $L_k=L_{k-1}\# L_1$  where $\#$ is defined
by replacing an arc connecting cusp edges of the two front projection by the tube shown on the right of
Figure~\ref{fig:bs}. One can show \cite{EES1} that the grading $-1$ linearized contact homology of $L_k$ has
dimension $k$:
\[\dim(L_1CH_{-1}(L_k))=k.\]
Thus none of the $L_k$'s are Legendrian isotopic. These examples generalize to all dimensions
but when $n$ is odd the knots $L_k$ are already distinguished by their $tb.$ Thus when $n$ is
odd to get an infinite family of examples of distinct Legendrian knots with the same classical
invariants one must alter the above construction. For details see \cite{EES1}.
\eex

For other examples of Legendrian knots in high dimensions and other constructions see \cite{EES1, EES3}.

\section{Applications}
There are numerous applications of Legendrian knot theory to topology and contact geometry. We
describe a few of these applications below.
\subsection{Legendrian surgery}
Let $(M,\xi)$ be a contact manifold and $L$ a Legendrian knot in $M.$ 
There is a standard neighborhood $N$ of $L$ with convex boundary. 
The dividing curves $\Gamma_{\partial N}$ consist of two parallel curves.
We can choose a framing on $N$ so that the dividing curves have slope $\infty$
and of course the meridian has slope $0.$ Let $M'$ be the manifold obtained
by gluing $S^1\times D^2$ to $\overline{M\setminus N}$ by the map $A:\partial(S^1\times D^2)\to \partial
\overline{M\setminus N}$ given by
\[A=\begin{pmatrix} 1 & 1\\ -1 &0
\end{pmatrix}.\]
Here we are using the framing on $\partial N$ and product structure on $\partial S^1\times D^2$ 
to write the map as a matrix.

The manifold $M'$ has a contact structure defined on the complement of $S^1\times D^2,$
but $\partial (S^1\times D^2)$ has an induced convex characteristic foliation. Using the 
classification of tight contact structures on the solid torus 
there is a unique tight contact structure on $S^1\times D^2$ inducing this characteristic
foliation. Thus we can get a well defined contact structure $\xi'$ on $M'.$ We say the contact
manifold $(M',\xi')$ is obtained form $(M,\xi)$ by \dfn{Legendrian surgery} on $L.$ 
There is an obvious generalization to Legendrian surgery on a Legendrian link in $(M,\xi).$
The fundamental result concerning Legendrian surgery is
\begin{thm}
If $(M,\xi)$ is a fillable contact structure then any contact manifold obtained by Legendrian 
surgery from $(M,\xi)$ is also fillable.
\end{thm}
Recall that a fillable contact structure is tight \cite{Eliashberg90a}. Moreover, the standard contact
structure on $S^3$ is fillable. Thus we can construct many tight contact structures
on three manifolds by doing Legendrian surgery on links in $S^3.$
\bex
Legendrian surgery on a Legendrian unknot in $S^3$ with $tb=-k$ yields tight contact
structures on the lens spaces $L(k+1,1).$

\begin{figure}[ht]
  \relabelbox 
  \small{\epsfxsize=4in \centerline{\epsfbox{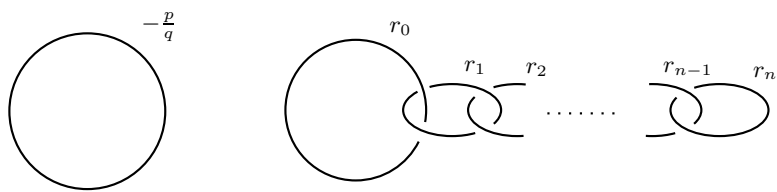}}}
  \relabel {0}{$r_0$}
  \relabel {1}{$r_1$}
  \relabel {3}{$r_{n-1}$}
  \relabel {2}{$r_2$}
  \relabel {4}{$r_n$}
  \relabel {p}{$-\frac{p}{q}$}
  \endrelabelbox
  \caption{Two surgery presentation of the lens space $L(p,q)$.}
  \label{lex}
\end{figure}

To get other lens spaces $L(p,q)$ consider the continued fractions expansion of $-\frac{p}{q}$:
\[-\frac{p}{q}=r_0-\frac{1}{r_1- \ldots \frac{1}{r_n}}.\]
Then the lens space $L(p,q)$ as a surgery presentation as shown in Figure~\ref{lex}, \cite{Rolfsen}.
It is easy to see that all the $r_i<-1$ thus we can realize the right hand side of Figure~\ref{lex}
by Legendrian unknots with $tb=r_i+1.$ Thus Legendrian surgery on this link will yield a
tight contact structure on $L(p,q).$ In is interesting to note that all tight contact
structures on $L(p,q)$ come from Legendrian surgery on a link in $S^3$ \cite{Giroux00, Honda1}. 
\eex
For an extensive discussion of the use of Legendrian surgery  to construct tight contact
structures see \cite{Gompf98}

One of the main questions in contact geometry is the following.
\begin{quest}
Does Legendrian surgery on a tight contact structure produce a tight contact structure?
\end{quest}
The answer is NO if the manifold is allowed to have boundary \cite{HondaDuke}
({\em cf.} \cite{Colin01a})
, but it is still possible
that for closed contact three manifolds this answer is YES.

\subsection{Invariants of contact structures}
Now that we have seen how to construct contact structures using Legendrian knots
we examine how to use
Legendrian knots to distinguishing contact structures on
three manifolds. For examples consider the contact structures
\[\xi_n=\ker(\sin(nz)dx + \cos (nz) dy),\]
on $T^3=T^2\times S^1$ where $x,y$ are coordinates on $T^2$,
and $z$ is
the coordinate on $S^1$. It is easy to see all these contact structures
on $T^3$ are tight (when pulled back to the universal cover of $T^3$ the
contact structures become the standard contact structure on $\R^3,$ \cite{Giroux00, Honda1, Kanda95})
and certainly look different. But how can they actually be distinguished?
They are in the same homotopy class of plane field so there is no ``algebraic'' way
to distinguish them.
Consider $\gamma=\{(x,y,z)| x=y=0\}.$ This is a Legendrian knot in $(T^3,\xi_n)$ and
the twisting of the contact planes relative to the natural product framing $\mathcal{F}$ 
is $tw(\gamma, \mathcal{F})=-n.$ The twisting with respect to $\mathcal{F}$ we will denote
$tb(\gamma),$ this is a slight abuse of notation but should not be confusing.
It turns out that this is the maximal twisting 
for any Legendrian knot in this knot type.
\begin{thm}
Let $\K$ be the topological knot type containing $\gamma.$ Then for the contact structure
$\xi_n$ we have
\[\overline{tb}(\K)=-n.\]
Thus all the contact structures $\xi_n$ on $T^3$ are distinct.
\end{thm}

Similar arguments have been used in distinguish some contact structures on other $T^2$ bundles
over $S^1$ and on some $S^1$ bundles over surfaces \cite{Giroux00, Honda1, Honda2}. 
This illustrates that
Legendrian knots are a very subtle invariant of contact structures. In fact, it is
possible that contact structures are completely determined by the Legendrian knots that 
exist within them.

\subsection{Plane curves}\label{sec:cusp}
Chekanov and Pushkar have recently solved Arnold's famous ``four vertex conjecture''\cite{Arnold93} using
invariants of Legendrian knots \cite{ChekanovPushkar}. Arnold conjecture is about generic wave fronts. A
\dfn{generic wave front} is a curve in $\R^2$ that is immersed at all but a finite number of
points at which there is a semi cubic cusp. Thus a curve looks very much like the front projection
of a Legendrian knot discussed in Section~\ref{sec:lk} except there might be vertical tangents.
Consider the manifold $M=\R^2\times S^1$ with the contact structure
\[\xi=\ker((cos \theta) dx + (sin \theta) dy),\]
where $x,y$ are coordinates on $\R^2$ and $\theta (\text{mod } 2\pi)$ is the coordinate on $S^1.$
(If we think of $M$ as the unit cotangent bundle of $\R^2$ then 
$\xi$ is a natural contact structure induced
from the standard complex structure on $T^*\R^2.$) Now if $\gamma$ is a generic (oriented) wave front in $\R^2$
then we can lift $\gamma$ to $\gamma_l$ in $M$ 
by sending $p\in \gamma$ to $(p,\theta)$ in $M$ where $\theta$ is the
$\frac{\pi}{2}$ plus the angle $T_p\gamma $ makes with the $x$-axis in $\R^2.$  
Clearly $\gamma_l$ is a Legendrian curve in $M.$
In addition generic Legendrian curves in $M$ will project to generic wave fronts in $\R^2.$  Thus
we may clearly study wave fronts by studying their corresponding Legendrian lifts.

Suppose $\gamma_t, t\in[0,1]$ is a family of generic (oriented)
wave fronts in $\R^2$ with $\gamma_0$ and $\gamma_1$ oppositely oriented smooth circular wave fronts
and for all $t,$ $\gamma_t$ having no self tangencies at which the orientations on the two
intersecting strands agree. Such a family is usually called an ``eversion of a smooth circular 
wave front with no
dangerous self tangencies''. Note this last condition guarantees that the Legendrian lift of each $\gamma_t$
is an embedded Legendrian knot in $M.$  Thus such a family of wave fronts gives a Legendrian isotopy of
the Legendrian lifts.
\begin{conj}[Arnold's Four Cusp Conjecture, \cite{Arnold93}]
An eversion of a smooth circular wave front with no
dangerous self tangencies must contain some curves with at least four cusps.
\end{conj}
In \cite{Arnold93} Arnold verified this conjecture in certain cases. Using an invariant for Legendrian knots
in $M$ similar to the decomposition invariant described in Section~\ref{sec:di}, Chekanov and 
Pushkar have verified 
the conjecture in all cases \cite{ChekanovPushkar}.

\subsection{Knot concordance}\label{app:concord}
Using the inequality
\[\frac12(tb(L)+|r(L)|+1)\leq g_s(\K)\]
for all $L\in\L(\K)$ 
from Section~\ref{sec:sg} we can easily make many observations about (non) slice knots
and, in particular, find many examples of elements of infinite order in the knot 
concordance group. We be begin with the simple observation
\begin{lem}
If there is a Legendrian knot $L\in\L(\K)$ with $tb(L)+|r(L)|+1\geq 0$ then
$\K$ is not smoothly slice. Moreover $\K$ is of infinite order in the knot concordance
group.
\end{lem}
The first observation follows directly from the above inequality. The second follows from
the inequality and the fact that $tb(L+L')=tb(L)+tb(L')+1.$ 

\begin{figure}[ht]
  \relabelbox 
  \small{\epsfxsize=4in \centerline{\epsfbox{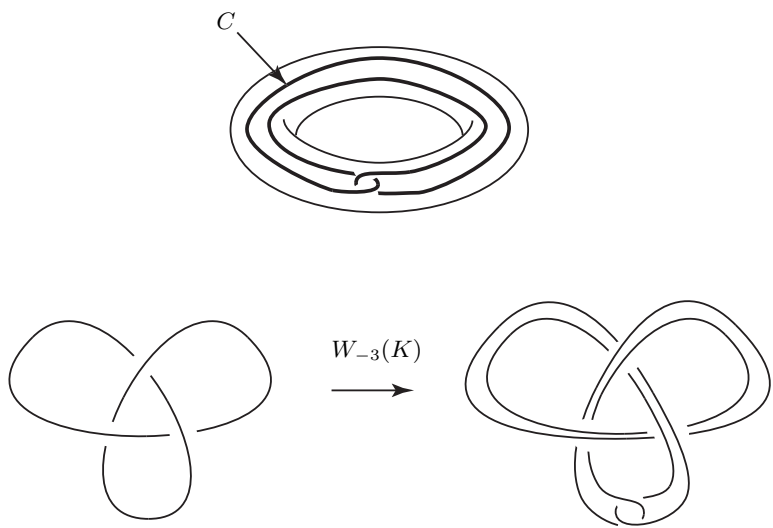}}}
  \relabel {C}{$C$}
  \relabel {W}{$W_{-3}(K)$}
  \endrelabelbox
  \caption{The solid torus used in the Whitehead double (top) and the $-3$ twisted double
    of the left handed trefoil knot.}
  \label{double}
\end{figure}

\begin{figure}[ht]
  \relabelbox 
  \small{\epsfxsize=4in \centerline{\epsfbox{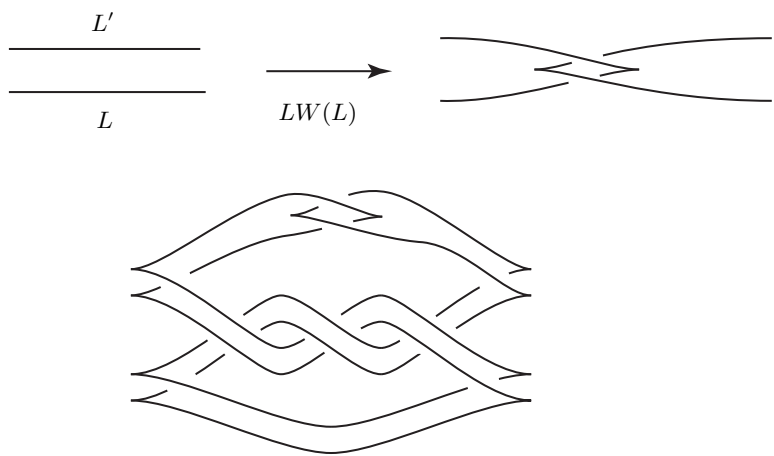}}}
  \relabel {L'}{$L'$}
  \relabel {L}{$L$}
  \relabel {W}{$LW(L)$}
  \endrelabelbox
  \caption{The change to $L\cup L'$ used to make $LW(L).$ Legendrian 
    Whitehead double of the trefoil (bottom).}
  \label{Ldouble}
\end{figure}

We now apply these ideas to study Whitehead doubles of knots.
If the neighborhood $N$ of a knot $K$ is replaced by the solid torus 
$D^2\times S^1$ shown if Figure~\ref{double}, 
so that $\{pt\}\times S^1,$ 
with $pt\in \partial D^2,$ is identified with a longitude for $K$ on $\partial N$ then the
image of $C$ is called the \dfn{$0$-Whitehead double} of a knot $K$ and is denoted $W_0(K).$
The \dfn{$n$-twisted Whitehead double} $W_n(K)$ is obtained in an analogous way after
putting $n$-full right handed twist in the solid torus $D^2\times S^1.$
One may compute that the Alexander polynomial of $W_n(K)$ is 0 for any know type. 
Freedman \cite{Freedman} has used this to show that the zero twisted Whitehead double of any
knot is {\em topologically locally flatly} slice. 
There is a well known conjecture:
\begin{conj}
A knot is smoothly slice if and only if its zero twisted Whitehead double is.
\end{conj}
In relation to this conjecture we can show the following theorem is true.
\begin{thm}[Rudolph 1995, \cite{Rudolph}]
Given any knot type $\K$ for which $\overline{tb}(\K)\geq 0$ then
all iterated zero twisted Whitehead doubles of $\K$ are not smoothly slice.
\end{thm}
To prove this theorem we consider the a
Legendrian version of the Whitehead double. Let $L$ be a Legendrian knot in
$(\R^3,\xi_{std})$  and $L'$ a parallel copy pushed up (in the $z$-direction)
slightly from $L.$ The Legendrian Whitehead double $LW(L)$ of $L$ is obtained by replacing
a horizontal stand of $L$ and $L'$ by a clasp as shown in Figure~\ref{Ldouble}.
One may easily
compute that as a topological knot $LW(L)$ is $W_{tb(L)}(L)$ and that
\[tb(LW(L))=1\quad \text{ and }\quad r(LW(L))=0.\]
Thus given \K let $L$ be an Legendrian knot in $\L(\K)$ realizing maximal $tb.$ Forming
$LW(L)$ we see that $W_{\overline{tb}(\K)}(\K)$ is of infinite order in the concordance 
group. By considering the Legendrian doubles of stabilizations of $L$ we see that $W_{n}(\K)$ is
of infinite order for all $n\leq \overline{tb}(\K).$ In particular if $\overline{tb}(\K)\geq 0$ then
the zero twisted Whitehead double is not slice.


\subsection{Invariants of classical knots}
It is quite easy to construct invariants of topological knots out of Legendrian knot theory.
For example the maximal Thurston--Bennequin invariant is a topological knot invariant
\[\overline{tb}(\K)=\text{max}\{tb(L)|L\in\L(\K)\}.\]
As we have seen in Section~\ref{sec:sg} this invariant is closely related to the slice genus
of a knot. Moreover it is very sensitive to mirroring. 
\bex
Let $\K$ be a $(p,q)$ torus knot, with $p,q>0.$ It's mirror $\overline{\K}$ is a $(-p,q)$ torus
knot. From Section~\ref{sec:torus} we know $\overline{tb}(\K)=pq-p-q$ and $\overline{tb}(\overline{\K})=-pq.$ 
\eex
It also seems that $\overline{tb}$ is sensitive to many topological operations on a knot,
like mutation. This prospect seems very interesting but has not been pursued.

Analogous to the definition of finite type invariants of topological knots one can define
finite type invariants of Legendrian knots.
\begin{thm}[Fuchs and Tabachnikov 1997, \cite{FuchsTabachnikov97}]
Any finite type Legendrian knot invariant for Legendrian knots in $(\R^3, \xi_{std})$ 
with the same $tb$ and $r$ is
also finite type invariant of the underlying framed topological knot type. 
\end{thm}
Put another way, two distinct Legendrian knots with the same topological knot type, $tb$ and $r$
cannot be distinguished by a finite order invariant.
This theorem has been generalized to other contact manifold, \cite{Tchernov}.

One might try to use this theorem to try to construct infinite order invariants of
topological knots out of Legendrian knot theory. For example let $\L_m(\K)$ be all
the Legendrian knots with maximal $tb$ and set $CH(\K)=\{CH_*(L)|L\in \L\}.$ So
$CH(\K)$ is the set of all contact homologies for Legendrian knots in $\L(\K)$ with
maximal $tb.$ Since $CH_*(L)$ is not a finite type invariant of $L$ it would be surprising
if $CH(\K)$ is a finite type invariant of $\K.$ More generally we ask the following question.
\begin{quest}
Can one extract non finite type invariants of topological knots $\K$ out of the DGA associated
to Legendrian knots in $\L(\K)$?
\end{quest}

\subsection{Contact homology and topological knot invariants}
In this section we discuss invariants of topological knots in $\R^3$ that come from from 
high dimensional Legendrian knot theory. Specifically, to a topological knot in $\R^3$ we will 
construct a Legendrian torus
in $\R^3\times S^2.$  The Legendrian isotopy class of this torus will be an
invariant of the knot. Thus any invariant of the Legendrian torus will be an invariant
of the topological knot.

Let $W=\{(x,v)\in T^*\R^3| v\in T^*_xM \text{ and } |v|=1\}$ be the unit cotangent bundle of $\R^3.$ 
It is easy to see there is a natural contact structure $\xi$ on $W.$ Specifically let 
$x_i$ be coordinates on $\R^3$ and $y_i$ be coordinates on the fibers of $T^*\R^3.$ The
restriction of $\lambda=\sum y_i dx_i$ to $W$ is a contact form on $W.$ Let $\xi=\ker\alpha.$ 

Given a topological knot $K$ in $\R^3$ let $T_K=\{(x,v)\in W| v(T_xK)=0 \}$ be the unit 
conormal bundle. One may easily check that $T_K$ is a Legendrian torus in $W$ and that
a topological isotopy of $K$ will produce a Legendrian isotopy of the associated torus.
Thus $T_K$ is an invariant of $K.$ 

As discussed above there is no ``classical invariants'' of $T_K,$ but one can use the contact homology 
of $T_K$ to come up with an invariant of $K.$ It is actually somewhat difficult to compute the
contact homology (not to mention contact homology has not been shown to be well defined in this
case yet), but following this idea Ng \cite{Ng1} has combinatorially defined the invariants of
topological knots and braids. 
Moreover, in \cite{Ng1} it is shown that these contact homology
invariants can distinguish knots with the same Alexander polynomial and signature.
In \cite{Ng2} a topological interpretation is given to a small piece of the
contact homology invariant.
Here we describe some of the geometric idea behind the braid invariants
and refer the reader to \cite{Ng1, Ng2} for details on knot invariant and the combinatorial
proofs of invariance. 

Let $U$ be the unknot in $\R^3$ represented as a 1-braid. So $T_U$ is a Legendrian torus in
$W.$ A Legendrian submanifold always has a neighborhood contactomorphic to the 1-jet space
of the submanifold \cite{GeigesIntro, McDuffSalamon}. 
That is, given a manifold $M$ its one jet space is $J^1(M)=T^*M\times \R$ with
the contact form $\alpha=dz+\lambda$ (here $z$ is the coordinate in the $\R$ direction and
$\lambda$ is the tautological 1-form on $T^*M$). It is easy to see that the zero section $Z$ in
$J^1(M)$ is a Legendrian submanifold. Thus the Legendrian torus $T_U$ has a neighborhood 
contactomorphic to $J^1(T^2).$ 

Any braid $B$ in $\R^3$ can be isotoped into an arbitrarily small neighborhood of the unknot $U.$
Moreover, any isotopy of the braid can be done in this small neighborhood. Thus the Legendrian
manifold $T_B$ in $W$ can be thought to sit inside $J^1(T^2)$ and the Legendrian isotopy class of
$T_B$ as a subset of $J^1(T^2)$ depends only on the isotopy class of $B$ as a braid.

There is a projection $\pi:J^1(T^2)\to T^*T^2$ that has properties identical to the
Lagrangian projection discussed in Sections~\ref{highknot}. Thus $\pi(T_B)$ is a Lagrangian submanifold
of $T^*T^2$ and generically has a finite number of transverse double points. We can now define the
contact homology of $T_B$ as in Section~\ref{highDGA}. It is the contact homology of $T_B$ that is
the invariant of $B.$

We now describe (at least heuristically) how to compute $CH_*(T_B).$ The 1-jet space also has a 
``front projection'' $\Pi:J^1(T^2)\to T^2\times \R.$ This front projection has all the properties
described in Section~\ref{highknot} for the standard front projection. If $B$ is an $n$-braid we can find 
$n$ functions $f_i:[0,1]\to R^2$ that parameterize $B$ in the sense that their
graphs in $[0,1]\times\R^2/(0,x)\sim(1,x)$ give $B.$ With such a representation of $B$ we can
express $\Pi(T_B)$ as the union of the graphs of the functions 
$F_i(t,\theta)=f_i(t)\cdot(\cos 2\pi\theta,
\sin 2\pi\theta):[0,1]^2\to \R,$  where $[0,1]^2$ glued to form $T^2.$ As described in Section~\ref{highex}
the double point the the Lagrangian projection of $T_B$ correspond to the critical points of
$F_i-F_j,$ and these correspond to the critical points of $|f_i-f_j|.$ The function
$|f_i-f_j|$ must have a maximum and a minimum, we label these $b_{ij}, a_{ij},$ respectively.
Though there may be other critical points for $|f_i-f_j|$ we can eliminate them up to stable
tame isomorphism. There are no cusps in the front projection of $T_B$ so the grading on the
double points is one less than the Morse index of $F_i-F_j,$ thus $|a_{ij}|=0$ and $|b_{ij}|=1.$ 

We now describe the differential. Clearly we have
\[\partial a_{ij}=0.\]
The differential on the $b_{ij}$'s is somewhat more complicated.
\begin{equation} 
    \partial b_{ij} = a_{ij} -  \phi_B(a_{ij}).
\end{equation}
Where $\phi_B$ is defined as follows: let for a generator $\sigma_k$ of the braid
group we define
\begin{equation}
\phi_{\sigma_k}:   \begin{cases}
        a_{ki} \mapsto -a_{k+1,i}-a_{k+1,k}a_{ki} & i \neq k, k+1 \\
        a_{ik} \mapsto -a_{i,k+1}-a_{ik}a_{k,k+1} & i\neq k, k+1\\
        a_{k+1,i} \mapsto a_{ki} & i\neq k, k+1\\
        a_{i,k+1} \mapsto a_{ik} & i\neq k, k+q\\
        a_{k,k+1} \mapsto a_{k+1,k} &\\
        a_{k+1,k} \mapsto a_{k,k+1} &\\
        a_{ij} \mapsto a_{ij} & i,j \neq k, k+1.
    \end{cases}
  \end{equation}
If $B=\sigma_{i_1}\ldots \sigma_{i_l}$ then define $\phi_B=\phi_{\sigma_{i_l}}\circ \ldots\circ 
\phi_{\sigma_{i_1}}.$ This differential $\partial$ is computed by looking for holomorphic disks in 
$T^*(T^2)$ with boundary on the projection of $T_B.$ This is done by considering ``gradient 
flow trees'' for all the functions $F_i-F_j.$ See \cite{Ng1}. In addition, see \cite{Ng1},
for a description of the knot invariant.

\end{document}